\documentclass[a4paper,10pt]{amsart}
\usepackage{hyperref}
\usepackage{amsmath,amstext,amsthm,amscd,amsopn,verbatim,amssymb}

\usepackage[margin=1in]{geometry}
\usepackage{tikz}
\usetikzlibrary{matrix}
\usetikzlibrary{shapes}
\usetikzlibrary{arrows}
\tikzset{ext/.style={circle, draw,inner sep=1pt},int/.style={circle,draw,fill,inner sep=1pt},nil/.style={inner sep=1pt}}
\tikzset{exte/.style={circle, draw,inner sep=3pt},inte/.style={circle,draw,fill,inner sep=3pt}}
\tikzset{diagram/.style={matrix of math nodes, row sep=3em, column sep=2.5em, text height=1.5ex, text depth=0.25ex}}
\tikzset{diagram2/.style={matrix of math nodes, row sep=0.5em, column sep=0.5em, text height=1.5ex, text depth=0.25ex}}
\tikzset{descr/.style={}}

\theoremstyle{plain}
  \newtheorem{thm}{Theorem}[section]
  \newtheorem{defi}[thm]{Definition}
  \newtheorem{deflemma}[thm]{Definition/Lemma}
  \newtheorem{prop}[thm]{Proposition}
  \newtheorem{cor}[thm]{Corollary}
  \newtheorem{lemma}[thm]{Lemma}
\theoremstyle{definition}
  \newtheorem{ex}[thm]{Example}
  \newtheorem{rem}[thm]{Remark}

\newcommand{\alg}[1]{\mathfrak{{#1}}}
\newcommand{\co}[2]{\left[{#1},{#2}\right]} % commutator
\newcommand{\pderi}[2]{ { \frac{\partial {#1} }{\partial {#2} } } }

\newcommand{\ad}{{\text{ad}}}

\newcommand{\p}{\partial}
\newcommand{\Hom}{\mathrm{Hom}}

\newcommand{\C}{{\mathbb{C}}}
\newcommand{\R}{{\mathbb{R}}}

\newcommand{\CG}{{\mathsf{ICG}}}
\newcommand{\ICG}{{\mathsf{ICG}}}
\newcommand{\TCG}{{\mathsf{TCG}}}
\newcommand{\Ne}{{\mathcal{N}}} % Nerve
\newcommand{\Graphs}{{\mathsf{Graphs}}}
\newcommand{\Gra}{{\mathsf{Gra}}}

\newcommand{\fullGraphs}{{\mathsf{fGraphs}}}
\newcommand{\fGraphs}{{\mathsf{fGraphs}}}

\newcommand{\Br}{{\mathsf{Br}}}
\newcommand{\op}{\mathcal}

\newcommand{\bpm}{\begin{pmatrix}}
\newcommand{\epm}{\end{pmatrix}}
\DeclareMathOperator{\bs}{\mathbf{s}}
\newcommand{\tr}{\mathit{tr}}
\newcommand{\FreeLie}{{\F_{\Lie}}}

\newcommand{\fC}{\mathit{fC}}

\newcommand{\GF}{\mathbb{K}}
\newcommand{\id}{{\mathit{id}}}

\renewcommand{\S}{\mathbb{S}}
\newcommand{\bbS}{\mathbb{S}}

\newcommand{\Tpoly}{ T_{\rm poly} }
\newcommand{\Dpoly}{ D_{\rm poly} }
\newcommand{\PCD}{ \mathit{PaCD} }
\newcommand{\PaCD}{ \PCD }
\newcommand{\PaP}{ \mathit{PaP} }

\newcommand{\mV}{{\mathcal{V}}}
\newcommand{\mU}{{\mathcal{U}}}
\newcommand{\mF}{{\mathcal{F}}}

\DeclareMathOperator{\End}{End}
\DeclareMathOperator{\Aut}{End}
\DeclareMathOperator{\Der}{Der}

\newcommand{\whl}{\circlearrowleft}

\newcommand{\sder}{\alg{sder}}

\newcommand{\grt}{\alg{grt}}
\newcommand{\GRT}{\mathrm{GRT}}
\newcommand{\GC}{\mathsf{GC}}
\newcommand{\dGC}{\mathsf{dGC}}

\newcommand{\fullGC}{\mathsf{fGC}}
\newcommand{\fGC}{\mathsf{fGC}}

\newcommand{\dgra}{\mathsf{dgra}}

\newcommand{\Def}{\mathsf{Def}}
\newcommand{\hoe}{\mathit{hoe}}
\newcommand{\hoLie}{\mathit{hoLie}}
\newcommand{\Lie}{\mathit{Lie}}
\newcommand{\Ass}{\mathit{Ass}}
\newcommand{\coAss}{\Ass^*}%{\mathit{coAss}}
\newcommand{\coCom}{\Com^*}%{\mathit{coCom}}
\newcommand{\coLie}{\Lie^*}%{\mathit{coLie}}
\newcommand{\coe}{e^*}%{\mathit{coe}}
\newcommand{\Com}{\mathit{Com}}
\newcommand{\hoCom}{\mathit{Com}_\infty}
\newcommand{\conn}{{\mathit{conn}}}
\newcommand{\F}{\mathbb{F}}
\newcommand{\Q}{\mathbb{Q}}
\newcommand{\Z}{\mathbb{Z}}
\newcommand{\N}{\mathbb{N}}
\newcommand{\K}{\mathbb{K}}
\newcommand{\im}{\mathit{Im}}

\newcommand{\sgn}{\mathit{sgn}}
\newcommand{\coker}{\mathit{coker}}
\newcommand{\bracket}{{\co{}{}}}
\newcommand{\Harr}{\mathit{Harr}}
\newcommand{\poly}{ {\rm poly} }
\newcommand{\Tw}{{\mathit{Tw}}}

\newcommand{\gr}{{\mathrm{gr}}}

\begin{document}
\title[M. Kontsevich's graph complex and $\alg{grt}$]{M. Kontsevich's graph complex and the Grothendieck-Teichm\"uller Lie algebra}
%\shorttitle{ee}
\author{Thomas Willwacher}
\address{Department of Mathematics\\ University of Zurich\\ Winterthurerstrasse 190 \\ 8057 Zurich, Switzerland}
\email{thomas.willwacher@math.uzh.ch}

\thanks{The author was partially supported by the Swiss National Science Foundation (grant 200020-105450).}
% \subjclass[2000]{16E45; 53D55; 53C15; 18G55}
% \date{}
\keywords{Formality, Deformation Quantization, Operads}

\begin{abstract}
We show that the zeroth cohomology of M. Kontsevich's graph complex is isomorphic to the Grothendieck-Teichm\"uller Lie algebra $\grt_1$. The map is explicitly described. This result has applications to deformation quantization and Duflo theory. 
We also compute the homotopy derivations of the Gerstenhaber operad. They are parameterized by $\grt_1$, up to one class (or two, depending on the definitions). More generally, the homotopy derivations of the (non-unital) $E_n$ operads may be expressed through the cohomology of a suitable graph complex.
Our methods also give a second proof of a result of H. Furusho, stating that the pentagon equation for $\grt_1$-elements implies the hexagon equation.
\end{abstract}
\maketitle

\section{Introduction}
The Grothendieck-Teichm\"uller group $\GRT_1$ is a pro-unipotent group introduced by V. Drinfeld \cite{drinfeld}, based on ideas of A. Grothendieck. 
This group is of great interest because it plays a central role in a variety of constructions in mathematics. 
The Grothendieck-Teichm\"uller group acts freely transitively on the set of Drinfeld associators. These associators are the central algebraic objects in Drinfeld's construction of quasi-Hopf algebras \cite{drinfeld} and the quantization of Lie bialgeras \cite{ek1, tamenquant}, in deformation quantization \cite{tamdefq}, in showing formality of the little disks operad \cite{tamarkin}, in the study of multiple zeta values \cite{lemurakami, brown, furushodshuffle} or may be used to solve the Kashiwara-Vergne conjecture in Lie theory \cite{AET, schneps}. 
The original (profinite) version of this group was proposed by Grothendieck to study the Galois group $\mathrm{Gal}(\bar \Q/\Q)$. For a more detailed survey we refer the reader to \cite{schnepssurvey, furushosurvey}.

In this paper we will add two ``stories'' to the above list. 
The first is related to M. Kontsevich's graph complex $\GC_2$, which is a combinatorially defined complex (roughly linear combinations of isomorphism classes of graphs) introduced in the study of the formality conjecture \cite{K3}. This complex in fact carries the structure of a differential graded Lie algebra, which acts on M. Kontsevich's formality morphisms. Little is known about its cohomology $H(\GC_2)$.
The first main result of this paper is the following.
\begin{thm}
\label{thm:main}
The zeroth cohomology of the Kontsevich graph complex, considered as Lie algebra, is isomorphic to the Grothendieck-Teichm\"uller Lie algebra.
\[
 H^0(\GC_2) \cong \grt_1
\]
Furthermore the cohomology in negative degrees vanishes.
\[
 H^{<0}(\GC_2) =0
\] 
\end{thm}
The cohomology in positive degrees is still unknown. In particular, it is an open conjecture that $H^1(\GC_2)=0$. On the other hand computer experiments \cite{barnatanmckay} show that  $H^3(\GC_2)\neq 0$. 

The second main result is about the homotopy theory of the Gerstenhaber operad $e_2$, which is quasi-isomorphic to the operad of chains of the little disks operad.
D. Tamarkin has shown \cite{tamarkingrtaction} that $\grt_1$ acts faithfully on some (and hence any) cofibrant resolution of $e_2$ by operadic derivations.
I. e., there is an injective map of Lie algebras
\[
\grt_1 \hookrightarrow H^0(\Der(\hoe_2))
\] 
where $\hoe_2$ is the minimal quasi-free resolution of $e_2$ and $\Der(\dots)$ denotes the differential graded Lie algebra of derivations. Precise definitions will be given below.
The second main result of this paper is the computation of $H^0(\Der(\hoe_2))$. In particular we will show that the above inclusion is an isomorphism, up to one ``trivial'' class acting by scaling of the generators.

\begin{thm}
\label{thm:hoe2}
 The automorphism of the Gerstenhaber operad up to homotopy are given by $\grt_1$ and one class, i.e., there is an isomorphisms of Lie algebras
\[
 H^0(\Der(\hoe_2))\cong \grt_1 \rtimes \K =: \grt  %\cong H^0(\GC_2).
\]
where $\K$ acts on $\grt_1$ by multiplication with the degree with respect to the grading on $\grt_1$.
\end{thm}

Let $e_n$ be the operad governing $n$-algebras and let $\hoe_n$ be its minimal cofibrant resolution.
The third main contribution of this paper is to show that the cohomology of $\Der(\hoe_n)$ can be expressed through the cohomology of a suitable graph complex $\GC_n$, which was also introduced (in some form) by M. Kontsevich \cite{K4}. %In fact, the definitions of $\GC_n$ and $\GC_2$ differ only in signs and gradings. 

\begin{thm}
\label{thm:GCdef}
 Let $e_n$ be the operad governing $n$-algebras, $n=2,3,\dots$. Then
\begin{equation}
\label{equ:GCdef}
 H(\Der(\hoe_n)) \cong S^+(H(\GC_n)[-n-1] \oplus \K[-n-1] \oplus V_n[-n-1] )[n+1]
\end{equation}
as graded vector spaces where
\[
V_n = \bigoplus_{\substack{j\geq 1 \\ j\equiv 2n+1 \mod 4 } } \K[n-j]
\]
and $S^+(\cdots)$ denotes the completed symmetric tensor product, without the zeroth term (i.e., without $\K$). Furthermore the inclusion $H(\GC_n)\to  H(\Der(\hoe_n))$ is a map of Lie algebras.
\end{thm}

\begin{rem}
One may put a natural graded Lie algebra structure on $H(\GC_n) \oplus \K \oplus V_n$ as follows. $H(\GC_n) \oplus V_n$ is the cohomology of a larger graph complex $\fGC_{n,\conn}$ (to be defined below), which is naturally a differential graded Lie algebra. The extra generator acts by multiplication by the Euler characteristic minus one. Then Theorem \ref{thm:GCdef} may be strengthened by saying that the inclusion of $H(\GC_n) \oplus \K \oplus V_n$ is a map of graded Lie algebras. Note however that we do not make claims about the Lie algebra structure on the whole of $H(\Der(\hoe_n))$.
\end{rem}

\begin{rem}
It is a well known fact due to D. Tamarkin that an element of $H^0(\Der(\hoe_2))$ determines an $L_\infty$ derivation of the Lie algebra of polyvector fields $\Tpoly[1]$ up to homotopy. The map from $H^0(\Der(\hoe_2))$ to $H^0(\GC_2)$ encoded in Theorem \ref{thm:GCdef} is a ``stable'' version of this map. %To my knowledge, it has previously not been known whether the 
\end{rem}

\begin{rem}
 A result similar to Theorem \ref{thm:hoe2} has recently been obtained by B. Fresse \cite{fresse}.
The differences between his result and ours are: (i) B. Fresse considers $e_2$ as a Hopf operad and studies automorphisms of a resolution as Hopf operad and (ii) he computes the full automorphism group (up to homotopy), not just the Lie algebra.
The second difference is minor, since a pro-unipotent group is isomorphic to its Lie algebra. See Appendix \ref{app:unipotent} for a more detailed discussion.
The first difference is more severe and makes the results of \cite{fresse} and the present paper distinct. There is one ``link'': We construct the $\grt$ action below by acting on some Hopf operad quasi-isomorphic to $e_n$ (as Hopf operad). It will be apparent that the action respects the Hopf operad structure, see Remark \ref{rem:Hopfaction}. 
\end{rem}

\begin{rem}
 The author was made aware by V. Turchin recently that the deformation complexes of operad maps $e_m\to e_n$, but with the bracket being sent to zero, have also been considered and related to a version of graph cohomology in \cite{Turchin1, Turchin2, Turchin3, LambrechtsTurchin}.
\end{rem}

% 
% For completeness, we also state a generalization of Theorem \ref{thm:main}.
% \begin{thm}
% \label{thm:maingen}
% Let  $\alg{t}$ be the Drinfeld-Kohno (simplicial) Lie algebra, and let $H(\alg{t})$ be its simplicial cohomology. Then 
% \[
%  H(\alg{t})[3] \cong H(\GC_2)\oplus \R[1].
% \]
% \end{thm}
% \begin{rem}
% Since $H^3(\alg{t})\cong{\grt}$ Theorem \ref{thm:main} follows, except for the compatibility with the Lie algebra structure.
% The additional summand $\R[1]$ can be thought of as representing the graph with one vertex and one edge, forming a loop\footnote{the ``divergence operator graph''}. This graph is not included in the graph complex by convention.  
% %Here the right hand side could be more naturally written as $H(\GC_2^\whl)$, where $\GC_2^\whl$ is the wheeled version of the graph complex, 
% \end{rem}

\subsection{Structure of the paper}
We tried to keep the main body of this paper short, at the cost of a larger number of Appendices which contain some technical results.
In section 2 we fix our notation and recall some standard definitions.
 Section \ref{sec:grcomplexes} contains the definition of the graph complexes we work with. Section \ref{sec:themap} contains a proof of  Theorem \ref{thm:GCdef}. In section \ref{sec:grt} the definition of the Grothendieck-Teichm\"uller Lie algebra $\grt_1$ is recalled and some of its properties are derived.
The proof of the main Theorems \ref{thm:main} and \ref{thm:hoe2} follows in section \ref{sec:theproof}. Section \ref{sec:appdefq} contains some applications of this Theorem. Section \ref{sec:explicitmap} contains a largely independent ``pedestrian'' description of the map between $\GC_2$ and $\grt_1$. The important result is a computation of the leading order terms of the graph cohomology classes corresponding to the conjectural generators $\sigma_3, \sigma_5,\dots$ of $\grt_1$.
Finally the numerous appendices provide background for some constructions and notations used in the main text.
Of particular importance is Appendix \ref{sec:optwists}, in which a technical construction we call ``operadic twisting'' is introduced.

\subsection*{Acknowledgements}
I am very grateful to Anton Alekseev, Vasily Dolgushev, Pavel Etingof, Giovanni Felder, Benoit Fresse, David Kazhdan, Anton Khoroshkin, Sergei Merkulov, Pavol \v Severa and V. Turchin for many helpful discussions and their encouragement. In particular Anton Alekseev and Pavol \v Severa helped me a lot in shaping my understanding of $\grt$ and the graph complex. I am highly indebted to Vasily Dolgushev who carefully read the manuscript and pointed out several mistakes in the original version, one of which was relatively severe and made some changes to this manuscript necessary.

%I also thank the colleagues who pointed out to me (correctly) that the original version of this manuscript was badly written and hardly understandable. 
For the present revised version I added some more details, streamlined the presentation, fixed the mistakes and changed the notation a bit to adhere more to standard conventions. I apologize for potential referential inconsistency with the older version.

I am very grateful for support by the Harvard Society of Fellows during most of the writing of this work. 

As final note, let me also point out the more detailed account \cite{vasilynotes} on some of the results of this paper given by V. Dolgushev and C. Rogers.

\section{Notation and Conventions}
\label{sec:defcomplexes}
We work over a ground field $\K$ of characteristic zero. For $V$ a graded or differential graded (dg) $\K$-vector space, we denote its $r$-fold desuspension by $V[r]$. 
The operator on elements shifting degrees by one we denote by $\bs$, so that for $x\in V$, $\bs x\in V[1]$.
For $x\in V$ a homogeneous element, we denote by $|x|$ its degree, so that for example $|\bs x|=|x|-1$. 
In general we will use cohomological conventions, i.e., differentials will have degree +1. 
For a complex $V$ with differential $d=d_1+d_2$ such that $d_1^2=0$, we denote by $(V, d_1)$ the graded vector space $V$ endowed with the differential $d_1$. 

A ($\Z$-)grading on a vector space $V$ is a decompostion $V\cong \oplus_{n\in \Z}V_n$ into a direct sum of subspaces $V_n$. We similarly call a decomposition 
\[
 V\cong \prod_{n\in \Z}V_n
\]
into a direct \emph{product} of subspaces a \emph{complete grading} or just \emph{grading} (abusing notation) of the vector space $V$.

The symmetric groups will be denoted by $S_n$, $n=1,2,\dots$.
We denote the (completed) symmetric product space of the (dg) vector space $V$ by
\begin{equation}
\label{equ:Splus}
S^+(V) := \prod_{j\geq 1}(V^{\otimes j})^{S_j}
\end{equation}
where the symmetric groups act by permutations of the factors.
Sometimes we will also use the version with a $j=0$-term
\[
S(V) := \GF\oplus S^+(V). 
\]

 Our conventions about operads will mostly follow the textbook \cite{lodayval} by J.-L. Loday and B. Vallette. For $\op P$ an operad, we denote its space of $N$-ary operations by $\op P(N)$. The operadic $r$-fold desuspension is an operad $\op P\{r\}$ such that 
\[
\op P\{r\}(N) = \op P(N)\otimes \sgn_N^{\otimes r} [(N-1)r]
\]
where $\sgn_N$ is the sign representation of $S_N$. Identical notation is used for cooperads. 
All cooperads we will encounter arise as Koszul duals (see \cite[7.2]{lodayval}) of standard quadratic operads. In particular the Koszul dual cooperad of the commutative operad $\Com$ is the Lie cooperad (up to suspension)
\[
\Com^\vee = \coLie\{1\}
\] 
and similarly (see \cite[chapter 13]{lodayval})
\begin{align*}
\Lie^\vee &= \coCom\{1\} \\
\Ass^\vee &= \coAss\{1\}.
\end{align*}

A central role will be played by the operads $e_n$ governing $n$-algebras for $n=2,3,\dots$, see \cite[13.3.22]{lodayval}. Their Koszul duals are
\[
e_n^\vee = \coe_n\{n\}.
\]

If $\op C$ is a coaugmented cooperad, we will denote the quasi-free operad obtained by cobar construction (see \cite[7.3.3]{lodayval}) by $\Omega(\op C)$. If $\op P$ is a Koszul operad, then there is a canonical quasi-isomorphism
\[
\Omega(\op P^\vee) \to \op P
\]
and $\Omega(\op P^\vee)$ is the minimal resolution of $\op P$. The operads mentioned above and their desuspensions are Koszul and we use the following abbreviations:
\begin{align*}
\hoe_n &= \Omega(e_n^\vee) &  \hoLie_n &= \Omega((\Lie\{n-1\})^\vee) = \Omega(\coCom\{n\}) 
\end{align*}

The subscript on the right is chosen so that we have an embedding $\hoLie_n\to \hoe_n$.
Using a formality morphism of the little $n$-cubes operad (see \cite{tamarkin} for $n=2$ or \cite{K2}, \cite{LV}), $\hoe_n$ can be seen as a model for the $E_n$ operad, without zero-ary operations. 

For $\op C$ a cooperad and $\op P$ an operad, we denote the convolution dg Lie algebra by 
\begin{equation}
\label{equ:convalgebra}
\Hom_{\S}(\op C, \op P) = \prod_{N\geq 0} \Hom_{S_N}(\op C(N), \op P(N))
\end{equation}
as in \cite[section 6.4.4]{lodayval}. In the cases relevant to this work $\op C$ will always be coaugmented, with $\op C(1)$ one dimensional, $\op C(0)=0$ and we are given a map of operads $\Omega(\op C)\to \op P$. Such a map determines a Maurer-Cartan element $\alpha$ in the differential graded (dg) Lie algebra $\Hom_{\S}(\op C, \op P)$. We may twist by this Maurer-Cartan element to obtain a Lie algebra which we denote
\[
\Def(\Omega(\op C)\to \op P)
\] 
and call the deformation complex of the map $\Omega(\op C)\to \op P$. The corresponding notation in  
\cite[section 6.4.9]{lodayval} is $\Hom_{\S}^\alpha(\op C, \op P)$.

\begin{rem}
Our notation here is slightly non-standard. In some respect (and in some not) it is more natural to call 
\[
\prod_{N\geq 2} \Hom_{S_N}(\op C(N), \op P(N))
\]
the deformation complex. 
However, our convention will significantly streamline later proofs since we will not have to treat the $N=1$ case as special. In the relevant examples, $\op P(1)$ is one dimensional and the cohomology of both candidate deformation complexes differs by a one dimensional space.
\end{rem}

We will often use the following Lemma, which is a form of the left lifting (up to homotopy) property of $\Omega(\op C)$. 
\begin{lemma}
\label{lem:defpresqiso}
Let $\op C$ be a coaugmented cooperad, let $\op P$ and $\op P'$ be operads, and let 
\[
\Omega(\op C) \to \op P \to \op P'
\]
be operad maps, with the right hand arrow being a quasi-isomorphism. Then the induced map of differential graded Lie algebras
\[
\Def(\Omega(\op C)\to \op P) \to \Def(\Omega(\op C)\to \op P')
\]
is a quasi-isomorphism.
\end{lemma}
\begin{proof}
We have to show that the mapping cone is acyclic. Put a descending complete filtration on the mapping cone as follows:
\[
\mF^p := \prod_{N\geq p} \Hom_{S_N}(\op C(N), \op P(N))[1] \oplus \Hom_{S_N}(\op C(N), \op P'(N)).
\]
The associated graded is acyclic since $P \to \op P'$ is a quasi-isomorphism. Hence it follows by standard spectral sequence arguments that the mapping cone is acyclic as well.
\end{proof}

We recall the following definitions from \cite[section 5.1]{tamarkingrtaction}:
Let $\op P$ be a (differential graded) operad and let $\op P'$ be the same operad with zero differential. Define the dg commutative algebra $a_n=\K[\epsilon]/(\epsilon^2)$ with $|\epsilon|=-n$. A map of $\S$-modules $g:\op P'\to \op P'$ of degree $n$ is called a derivation of degree $n$ if the map
\[
\id + \epsilon g : \op P'\to \op P'\otimes a_n
\] 
is a map of operads.
%\footnote{The word non-unital is non-standard, but will simplify our proofs. }
 The derivations of degree $-n$ form a vector space and the space of derivations%\footnote{The word pseudo-derivation is not used in \cite{tamarkingrtaction}.} of $\op P$
\[
\Der(\op P)'
\]
is defined to be the direct sum of these spaces for all $n$. $\Der(\op P)'$ is endowed with a differential derived from the differential on $\op P$. Note that the ``true'' derivations of $\op P$ are the closed degree 0 elements of $\Der(\op P)'$, so the notation is slightly abusive. 
Furthermore there is a Lie bracket on $\Der(\op P)'$ by the commutator:
\[
\co{D}{D'} := D\circ D' -(-1)^{n n'} D'\circ D
\] 
where $D$ and $D'$ are derivations of degrees $n$ and $n'$. It may be checked that $\co{D}{D'}$ is indeed derivation of degree $n+n'$.
In the homotopy theory of operads, one is interested in studying $\Der(\tilde{ \op P})'$ and in particular its cohomolgy for some cofibrant replacement $\tilde{ \op P}$ of $\op P$. Up to homotopy it does not matter which cofibrant replacement is used.

The important examples for us are the operads $\op P=\hoe_n$. One may work out the definitions to see that the complex $\Der(\hoe_n)'$ is isomorphic to the codimension one subcomplex
\[
 \prod_{N\geq 2} \Hom_{S_N}(e_n^\vee(N), \op P(N))[1] \subset \Def(\hoe_n\stackrel{\id}{\longrightarrow} \hoe_n)[1].
\]
Note however that the Lie brackets are different and mind the degree shift. In order to streamline later proofs we will define the complex 
\[
\Der(\hoe_n) := \prod_{N\geq 1} \Hom_{S_N}(e_n^\vee(N), \op P(N))[1] \cong \Def(\hoe_n\stackrel{\id}{\to} \hoe_n)[1]
\]
and also call it the complex of derivations of $\hoe_n$. The Lie bracket naturally extends to this larger complex. 

\begin{rem}
Note that scaling by arity, i. e. the map
\begin{gather*}
A\colon \op P \to \op P \\
\op P(N) \ni x \mapsto (N-1) x 
\end{gather*}
is a derivation for any operad $\op P$, which commutes with all other derivations. Our definition of $\Der(\cdots)$ is such that this trivial derivation is rendered exact. For example:
\[
\Der(\hoe_n) = \Der(\hoe_n)' \oplus \K B
\]
where the coboundary of $B$ is $A$.
\end{rem}

For $\op C$ a coaugmented cooperad, $\op P$ an operad and $f: \Omega(\op C)\to \op P$ an operad map, there are canonical maps
\[
\Der(\Omega(\op C)) \to \Def(\Omega(\op C)\stackrel{f}{\to} \op P)[1] \leftarrow \Der(\op P)' 
\]  
by post- or pre-composition. For homogeneous $D\in \Der(\Omega(\op C))$, $D'\in \Der(\op P)'$ we will use the following notation for these maps
\begin{equation}
\label{equ:circnotation}
\begin{aligned}
D &\mapsto (-1)^{|D|} \bs f\circ D \\
D' &\mapsto (-1)^{|D'|} \bs D'\circ f.
\end{aligned}
\end{equation}

\begin{rem}
In the following, we will often encounter deformation complexes of the form 
$\Def(\hoe_n\to \op P)$, where $\op P$ is some operad.
The differential on this complex then has the form
\begin{equation}
\label{equ:diffsplit1}
d = \delta + \co{\alpha}{\cdot}
\end{equation}
where $\delta$ is in induced by the differential on $\op P$ and $\alpha$ is the Maurer-Cartan element defined by the map $\hoe_n\to\op P$. Most of the time, the latter map will in addition factor as 
\[
\hoe_n\to e_n\to \op P
\]
where the first map is the canonical projection. In this case the only nonvanishing part of $\alpha$ in the direct product decomposition \eqref{equ:convalgebra} is that for $N=2$. $e_n^\vee(2)$ is two dimensional and we may split $\alpha=\alpha_\wedge+\alpha_{\co{}{}}$ into one part for each of the two co-generators. We will furthermore abbreviate
\begin{align}
\label{equ:dwedge}
%\begin{aligned}
d_\wedge &:= \co{\alpha_\wedge}{\cdot} & d_{\co{}{}} &:= \co{\alpha_{\co{}{}}}{\cdot}. 
%\end{aligned}
\end{align}
Hence the differential becomes
\begin{equation}
\label{equ:diffsplit2}
d = \delta + d_\wedge + d_{\co{}{}} %\co{\alpha_\wedge}{\cdot} + \co{\alpha_{\co{}{}}}{\cdot} 
\end{equation}
\end{rem}

\begin{rem}
Although we will often follow the notation and terminology of \cite{lodayval}, there is one notable conflict of notation: We introduce in Appendix \ref{sec:optwists} the (novel) notion of \emph{operadic twisting} and will denote the twisted version of an operad $\op P$ by $\Tw\op P$. This has nothing to do with the notation $\Tw$ (for twisting cochains) as used in \cite[section 6.4.8]{lodayval}.
\end{rem}
\section{Graph operads and graph complexes}
\label{sec:grcomplexes}
\begin{figure}
\centering
\[
\begin{tikzpicture}[
vert/.style={circle,draw, minimum size=5pt, inner sep=0}, invvert/.style={inner sep=-1,minimum size=-1}, ext/.style={circle,draw, minimum size=5pt, inner sep=0}, scale=2 ]
\node (v0) at (0.3,1.5) [vert,label=90:2] {};
\node (v1) at (0.3,1) [vert,label=270:1] {};
\node (v2) at (1.3,1.5) [vert,label=90:1] {};
\node (v3) at (1,1) [vert,label=270:2] {};
\node (v4) at (1.6,1) [vert,label=270:3] {};
\node (v5) at (0.7,1.3) [invvert] {$\circ_1$};
\node (v6) at (1.9,1.3) [invvert] {=};
\node (v7) at (2.4,1) [vert,label=270:2] {};
\node (v8) at (2.6,1.3) [vert,label=180:1] {};
\node (v9) at (2.8,1) [vert,label=270:3] {};
\node (v10) at (2.6,1.6) [vert,label=90:4] {};
\node (v11) at (3.4,1) [vert] {};
\node (v12) at (3.6,1.3) [vert] {};
\node (v13) at (3.8,1) [vert] {};
\node (v14) at (3.3,1.5) [vert] {};
\node (v15) at (4.2,1) [vert] {};
\node (v16) at (4.6,1) [vert] {};
\node (v17) at (4.4,1.3) [vert] {};
\node (v18) at (4.7,1.5) [vert] {};
\node (v19) at (3,1.3) [invvert] {+};
\node (v20) at (4,1.3) [invvert] {+};
\draw (v0)--(v1) (v3)--(v2) (v2)--(v4) (v4)--(v3) (v7)--(v9) (v9)--(v8) (v8)--(v7) (v8)--(v10) (v14)--(v11) (v12)--(v11) (v11)--(v13) (v13)--(v12) (v15)--(v17) (v17)--(v16) (v16)--(v15) (v16)--(v18) ;
\end{tikzpicture}
\]
 \caption{\label{fig:grainsert} The operadic insertion in the operads $\Gra_n$.}
\end{figure}
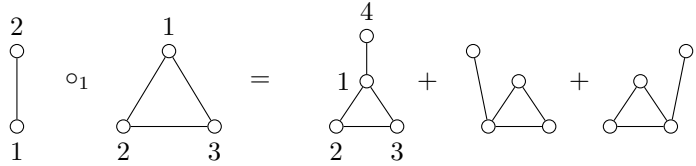

Several versions of graph complexes have been introduced by M. Kontsevich \cite{K3, K4}.\footnote{Strictly speaking, we consider the cohomological version here, while M. Kontsevich mainly considers the pre-dual, homological version in \cite{K4}. It does not matter much.}
Let us define them here in a combinatorial way. Let $\dgra_{N,k}$ be the set of directed graphs with $N$ numbered vertices and $k$ ordered directed edges. Concretely, a graph in $\dgra_{N,k}$ is given by a $k$-tuple of pairs $(i,j)$, $1\leq i,j \leq N$, a pair $(i,j)$ representing an edge from vertex $i$ to vertex $j$. We explicitly allow multiple edges and tadpoles, i. e., edges of the form $(i,i)$.
There is a natural action of the permutation group $P_k=S_k\ltimes S_2^{\times k}$ on $\dgra_{n,k}$ by permuting the order and flipping the directions of the edges. %Denote by $\sgn$ the sign representation of the symmetric group.
Denote by $\sgn_k$ the one-dimensional sign representations of $S_k$.
We will understand $\sgn_k$ (respectively $\sgn_2^{\otimes k}$) as a representation of $P_k$ by letting $S_2^{\times k}$ (respectively $S_k$) act trivially.
Define graded vector spaces as follows.
\[
 \Gra_n^\whl(N) = 
\begin{cases} 
 \bigoplus_{k\geq 0} \left( \GF\langle \dgra_{N,k} \rangle \otimes_{P_k} \sgn_k \right)[k(n-1)]
 &\quad \text{$n$ even} \\
\bigoplus_{k\geq 0} \left( \GF\langle \dgra_{N,k} \rangle \otimes_{P_k} \sgn_2^{\otimes k} \right)[k(n-1)] 
  &\quad \text{$n$ odd.}
 \end{cases}
\]
In words, we give each edge the degree $1-n$, and introduce the appropriate signs. For $n$ odd we additionally fix an orientation on the graph by prescribing directions on edges, identifying the graph with an edge direction flipped with minus the original graph.
The spaces $\Gra_n^\whl(N)$ naturally assemble to operads $\Gra_n^\whl$. The $S_N$ action is given by permuting the labels on the vertices. The operadic compositions $\Gamma_1 \circ_j \Gamma_2$ are given by `inserting'' the graph $\Gamma_2$ at vertex $j$ of graph $\Gamma_1$ and summing over all possible ways of reconnecting the edges incident to vertex $j$ in $\Gamma_1$ to vertices of $\Gamma_2$, see Figure \ref{fig:grainsert}.
 In case $n$ is even one needs to put an ordering (up to signed permutation) on the edges of the newly formed graphs. 
The natural choice is to consider the edges of $\Gamma_1$ to stand on the left of those of $\Gamma_2$.

We furthermore define 
\[
\Gra_n\subset \Gra_n^\whl
\]
to be the sub-operads spanned by graphs without tadpoles.

\begin{rem}
Note that by symmetry reasons, there can be tadpoles in graphs occurring in $\Gra_n^\whl$ only for $n$ even, and multiple edges only for $n$ odd. More precisely, for $n$ odd the operation of flipping a tadpole edge is odd, hence a graph containing a tadpole yields a zero element of $\Gra_n^\whl$ for $n$ odd. Similarly, for $n$ even the operation of permuting two edges of a multiple edge is odd. In particular $\Gra_n^\whl=\Gra_n$ for $n$ odd.
\end{rem}

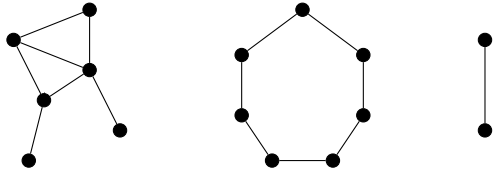
\begin{figure}
\centering
\[
\begin{tikzpicture}[
vert/.style={circle,draw,fill, minimum size=5pt, inner sep=0}, invvert/.style={inner sep=-1,minimum size=-1}, ext/.style={circle,draw, minimum size=5pt, inner sep=0}, scale=2 ]
\node (v0) at (3.4,1.3) [vert] {};
\node (v1) at (3.4,0.7) [vert] {};
\node (v2) at (0.3,1.3) [vert] {};
\node (v3) at (0.8,1.5) [vert] {};
\node (v4) at (0.5,0.9) [vert] {};
\node (v5) at (0.8,1.1) [vert] {};
\node (v6) at (1,0.7) [vert] {};
\node (v7) at (1.8,1.2) [vert] {};
\node (v8) at (1.8,0.8) [vert] {};
\node (v9) at (2,0.5) [vert] {};
\node (v10) at (2.4,0.5) [vert] {};
\node (v11) at (2.6,0.8) [vert] {};
\node (v12) at (2.6,1.2) [vert] {};
\node (v13) at (2.2,1.5) [vert] {};
\node (v14) at (0.4,0.5) [vert] {};
\draw (v13)--(v7) (v7)--(v8) (v8)--(v9) (v9)--(v10) (v10)--(v11) (v11)--(v12) (v12)--(v13) (v4)--(v2) (v4)--(v5) (v5)--(v6) (v3)--(v2) (v2)--(v5) (v5)--(v3) (v0)--(v1) (v4)--(v4) (v4)--(v14) ;
\end{tikzpicture}
\]
 \caption{\label{fig:GCgraphs} Three graphs in $\fullGC_n$. The differential in $\fullGC_n$ is the Lie bracket with the right hand graph.
 % The right hand graph is a Maurer-Cartan element.
 }
\end{figure}

The signs in the definition of the operad $\Gra_n$ are chosen so that there is an injective map 
\[
e_n \to \Gra_n 
\]
from the operad $e_n$ governing $n$ algebras to $\Gra_n$. The operad $e_n$ is generated by two binary operations, the product operation $\cdot \wedge \cdot$ of degree 0 and the bracket operation $\co{\cdot}{\cdot}$ of degree $1-n$. These operations are mapped to graphs as follows.
\begin{align}
\label{equ:enGramapdef}
\cdot \wedge \cdot &\mapsto
\begin{tikzpicture}[baseline=-0.65ex]
\node[ext, label=-90:1] at(0,0) {};
\node[ext, label=-90:2] at(.5,0) {};
\end{tikzpicture}
&
\co{\cdot}{\cdot} &\mapsto
\begin{tikzpicture}[baseline=-0.65ex]
\node[ext, label=-90:1](v) at(0,0) {};
\node[ext, label=-90:2](w) at(.5,0) {};
\draw (v) edge (w);
\end{tikzpicture}
\end{align}

By restriction, one obtains maps $\Lie\{n-1\} \to \Gra_n$ from the degree shifted Lie operad $\Lie\{n-1\}$, and hence also maps $\hoLie_{n} \to \Gra_n$ from its minimal resolution $\hoLie_{n}$ to $\Gra_n$.
The full graph complex $\fGC_n$ is by definition the deformation complex of the latter map.
\[
 \fullGC_n := \Def(\hoLie_{n} \to \Gra_n)
\]
As a graded vector space, $\fGC_n$ is just the space of (anti-)invariant elements in $\Gra_n$ under the action of the permutation group.
\[
 \fullGC_n \cong 
 \begin{cases}
 \prod_N  (\Gra_n(N)[n(1-N)])^{S_N} & \text{$n$ even} \\
  \prod_N  (\Gra_n(N)\otimes \sgn_N [n(1-N)])^{S_N} & \text{$n$ odd}.
 \end{cases}
\]
%Here the symmetric group action shall incorporate signs if $n$ is odd. 

\begin{rem}
In pictures, we shall draw elements of $\fullGC_n$ as undirected graphs with unlabelled black vertices. To actually obtain an element of $\fullGC_n$ from such a picture, one has to choose an ordering and directions of the edges, sum over all possible ways of assigning labels $1,2,\dots$ to the vertices and divide by the order of the symmetry group of that graph. The overall sign is left undetermined. See Figure \ref{fig:GCgraphs} for examples of elements in $\fGC_n$.
\end{rem}

By the description as a deformation complex, we know that $\fGC_n$ has the structure of a differential graded Lie algebra. Similarly, we define a dg Lie algebra $\fGC_n^\whl \supset \fGC_n$ as
\[
\fullGC_n^\whl := \Def(\hoLie_{n} \to \Gra_n^\whl).
\]
It differs from $\fullGC_n$ only in so far that tadpoles are permitted in graphs.

\begin{figure}
 \centering
\[
\begin{tikzpicture}[scale=1,
vert/.style={draw,outer sep=0,inner sep=0,minimum size=5,shape=circle,fill},
helper/.style={outer sep=0,inner sep=0,minimum size=5,shape=coordinate},
default_edge/.style={draw},
ext/.style={draw,outer sep=0,inner sep=0,minimum size=5,shape=circle},
every loop/.style={}]

\node [vert] (v21) at (8,8) {};
\node (v2) at (1.5,9) [vert] {};
\node (v12) at (1.5,8) [vert] {};
\node (v13) at (2.5,8) [vert] {};
\node [vert] (v15) at (3,9) {};
\node [vert] (v16) at (4,9) {};
\node [vert] (v18) at (7,9) {};
\node [vert] (v19) at (7,8) {};
\node [vert] (v23) at (8,9) {};
\node [vert] (v25) at (2.4,7.2) {};
\node [vert] (v26) at (1.331,6.545) {};
\node [vert] (v28) at (1.695,5.029) {};
\node [vert] (v30) at (3.106,5.029) {};
\node [vert] (v32) at (3.47,6.545) {};
\node [vert] (v34) at (2.4,6) {};
\node [vert] (v44) at (5.9,7.2) {};
\node [vert] (v45) at (4.831,6.545) {};
\node [vert] (v47) at (5.195,5.029) {};
\node [vert] (v49) at (6.606,5.029) {};
\node [vert] (v51) at (6.97,6.545) {};
\node [vert] (v53) at (5.695,6.029) {};
\node [vert] (v60) at (8.9,7.2) {};
\node [vert] (v61) at (7.831,6.545) {};
\node [vert] (v63) at (8.195,5.029) {};
\node [vert] (v65) at (9.606,5.029) {};
\node [vert] (v67) at (9.965,6.545) {};
\node [vert] (v69) at (8.695,6.029) {};
\node (v76) at (1.2,8.8) [helper,label=90:{$d$}] {};
\node (v77) at (1.2,7.8) [helper,label=90:{$d$}] {};
\node [helper,label=90:{$d$}] (v78) at (6.5,8.2) {};
\node [helper,label=90:{$d$}] (v79) at (0.895,5.629) {};
\node [helper,label=90:{$=-$}] (v80) at (2.3584,8.7731) {}; %\quad \frac{1}{2}
\node (v81) at (3.2,7.8) [helper,label=90:{$=0$}] {};
\node [helper,label=90:{$=0$}] (v82) at (8.7,8.2) {};
\node [helper,label=90:{$=$}] (v83) at (4.095,5.629) {};
\node [helper,label=90:{$+$}] (v85) at (7.295,5.629) {};
\node [vert] (v86) at (6.195,6.029) {};
\node [vert] (v87) at (9.195,6.029) {};

\draw[default_edge] (v19) to (v21);
\draw[default_edge] (v12) to (v13);
\draw[default_edge] (v15) to (v16);
\draw[default_edge] (v18) to (v19);
\draw[default_edge] (v21) to (v23);
\draw[default_edge] (v25) to (v26);
\draw[default_edge] (v26) to (v28);
\draw[default_edge] (v28) to (v30);
\draw[default_edge] (v30) to (v32);
\draw[default_edge] (v32) to (v34);
\draw[default_edge] (v32) to (v25);
\draw[default_edge] (v34) to (v26);
\draw[default_edge] (v34) to (v28);
\draw[default_edge] (v34) to (v30);
\draw[default_edge] (v34) to (v25);
\draw[default_edge] (v18) to (v23);
\draw[default_edge] (v18) to (v21);
\draw[default_edge] (v23) to (v19);
\draw[default_edge] (v44) to (v45);
\draw[default_edge] (v45) to (v47);
\draw[default_edge] (v47) to (v49);
\draw[default_edge] (v49) to (v51);
\draw[default_edge] (v53) to (v86);
\draw[default_edge] (v51) to (v44);
\draw[default_edge] (v53) to (v45);
\draw[default_edge] (v53) to (v47);
\draw[default_edge] (v86) to (v51);
\draw[default_edge] (v53) to (v44);
\draw[default_edge] (v60) to (v61);
\draw[default_edge] (v61) to (v63);
\draw[default_edge] (v63) to (v65);
\draw[default_edge] (v65) to (v67);
\draw[default_edge] (v67) to (v69);
\draw[default_edge] (v67) to (v60);
\draw[default_edge] (v69) to (v61);
\draw[default_edge] (v69) to (v63);
\draw[default_edge] (v87) to (v65);
\draw[default_edge] (v86) to (v49);
\draw[default_edge] (v69) to (v87);
\draw[default_edge] (v87) to (v60);
\end{tikzpicture}
\]
\caption{\label{fig:graphdex} Several example computations of the graph differential. Note that in drawing these pictures we are cheating a bit since the signs are ambiguous.
%we did not say what element of the graph complex the picture of a graph stands for. However, this is ``merely'' a matter of prefactors. The prefactors shown here can be recovered up to sign by interpreting a graph with unnumbered vertices as the sum over all numberings of vertices, divided by the order of the automorphism group.
}
\end{figure}
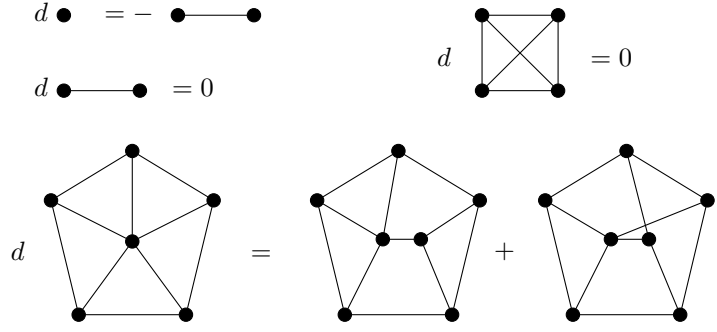

%
%Next, we want to twist the operads $\Gra_n^\pm$, see Appendix \ref{sec:optwists}. First, one obtains a graded Lie algebra
%\[
% \fullGC_{m,n}^\pm = \Def(\hoLie_m \stackrel{0}{\to} \Gra_n^\pm)
%\]
%of derivations of the zero map of operads. Elements of this Lie algebra are series of graphs as before, but with ``unidentifiable'' vertices. We will draw those vertices black in pictures, see Figure \ref{fig:GCgraphs}. In order to twist, we need a Maurer-Cartan element. We require that the graph with two vertices and one edge (see Figure \ref{fig:GCgraphs}) is such a Maurer-Cartan element. This has several consequences:
%\begin{enumerate}
% \item It must have degree 1. Hence we are forced to take $m=n$, i.e., consider $\fullGC_{n,n}^\pm$.
% \item It should not be zero by symmetry. Hence, we are forced to consider the case ``+'' for $n$ even and ``-'' for ``n'' odd.
%\end{enumerate}
%We will focus on these cases. The Maurer-Cartan element provides us with a differential $d$. We define the dg Lie algebras
%\[
% \fullGC_n = (\fullGC_{n,n}^{(\pm)^n}, d).
%\]
\begin{rem}
 Combinatorially, the differential applied to some graph $\Gamma$ has the form
\[
 d\Gamma \text{``}=\text{''} \sum_{v\in V(\Gamma)} \frac{1}{2}\text{(splitting of $v$)}-\text{(adding an edge at $v$)}.
\]
Here the sum runs over all vertices of $\Gamma$. The ``$\text{(splitting of $v$)}$'' means the vertex $v$ is replaced by a pair of vertices connected by an edge, and one sums over all possible ways of reconnecting the incoming edges at $v$ to the two newly created vertices. The term ``$\text{(adding an edge at $v$)}$'' stands for a graph obtained by adding a new vertex and connecting it to $v$. If there are any incoming edges at $v$, the second term cancels those graphs from the first term, in which all incoming edges at $v$ had been connected to one of the newly added vertices. Several examples can be found in Figure \ref{fig:graphdex}.
\end{rem}

Every graph can be seen as a union of its connected components. Hence we can write
\[
 \fullGC_n = S^+(\fullGC_{n,\conn }[-n])[n]
\]
where $\fullGC_{n,\conn}$ is the subcomplex of $\fullGC_n$ spanned by the connected graphs. Similarly,
\[
 \fullGC_n^\whl = S^+(\fullGC^\whl_{n,\conn }[-n])[n]
\]
where $\fullGC^\whl_{n,\conn}$ is the subcomplex of $\fullGC_n^\whl$ spanned by the connected graphs.
 We define M. Kontsevich's graph complex $\GC_n$ to be the subspace of $\fullGC_{n,\mathit{conn}}$ spanned by those graphs with all vertices of valence at least 3.
 %, with no tadpoles, but we \emph{allow} multiple edges.
 A significant part of the following proposition has been proven by M. Kontsevich.
 % \cite{K4,K5}.

\begin{prop}[partially contained in \cite{K4,K5}]
\label{prop:GCred}
 $\GC_n\subset \fullGC_n^\whl$ is a sub-dg Lie algebra. The cohomology satisfies 
\begin{align*}
H(\fullGC^\whl_{n,\mathit{conn}}) &= H(\GC_n) \oplus \bigoplus_{ \substack{j\geq 1 \\ j\equiv 2n+1 \mod 4} } \GF[n-j] \\
H(\fullGC_{n,\mathit{conn}}) &= H(\GC_n) \oplus \bigoplus_{ \substack{j\geq 3 \\ j\equiv 2n+1 \mod 4} } \GF[n-j].
\end{align*}
Here, the class corresponding to $\GF[n-j]$ is represented by a loop with $j$ edges, as in Figure \ref{fig:loops}.
\end{prop}
\begin{proof}[Sketch of proof]
A part of the proof is already contained in \cite{K4,K5}. First note that the differential does not produce any vertices of valence $\leq 1$, nor can it produce tadpoles or multiple edges if there were none in the graph before. Assume a graph $\Gamma$ with all vertices at least trivalent is given. Then the differential $d\Gamma$ contains graphs with one bivalent vertex, one such graph for every edge. However, each such graph comes twice with opposing signs, from splitting either one of the adjacent vertices. Pictorially:
\[
d
\begin{tikzpicture}[baseline=-.65ex]
\begin{scope}[]
\node[int] (v1) at (0,0) {};
\node[int] (v3) at (.5,0) {};
\draw (v1) edge (v3);
\foreach \h in {-.3,0,.3}
{
\draw (v1) -- +(-.3,\h);
\draw (v3) -- +(.3,\h);
}
\end{scope}
\end{tikzpicture}
=
\begin{tikzpicture}[baseline=-.65ex]
\begin{scope}
\node[int] (v1) at (0,0) {};
\node[int] (v2) at (0.3,0) {};
\node[int] (v3) at (1,0) {};
\draw (v2) edge (v1) edge (v3);
\foreach \h in {-.3,0,.3}
{
\draw (v1) -- +(-.3,\h);
\draw (v3) -- +(.3,\h);
}
\end{scope}
\end{tikzpicture}
-
\begin{tikzpicture}[baseline=-.65ex]
\node[int] (v1) at (0,0) {};
\node[int] (v2) at (0.7,0) {};
\node[int] (v3) at (1,0) {};
\draw (v2) edge (v1) edge (v3);
\foreach \h in {-.3,0,.3}
{
\draw (v1) -- +(-.3,\h);
\draw (v3) -- +(.3,\h);
}
\end{tikzpicture}
=0
\]

 Hence we can conclude that $d\GC_n\subset \GC_n$. This shows that $\GC_n$ is indeed a subcomplex.
Next Let $\GC_n^1$ be the space spanned by graphs having at least one vertex of valence $1$, and let $\GC_n^2$ be the space spanned by graphs having no vertex of valence one, but at least one vertex of valence two. 
Such vertices cannot be killed by the differential and hence we have a decomposition (of complexes)
\[
 \fullGC_{n,\mathit{conn}} \cong \GC_n^{\geq 3} \oplus \GC_n^2 \oplus \GC_n^1
\]
where $\GC_n^{\geq 3}$ is spanned by graphs containing only at least trivalent vertices. 

We claim that $\GC_n^1$ is acyclic. 
Indeed for the subcomplex of graphs not containing a trivalent vertex (i.e., ``linear graphs''), this is easily shown. Assume next that a graph, say $\Gamma$, has at least one trivalent vertex. We call an ``antenna'' a maximal connected subgraph consisting of one- and two-valent vertices in $\Gamma$. Then full graph $\Gamma$ can be seen as some ``core graph'' (the complement of the union of all antennas) with antennas of various lengths attached. See Figure \ref{fig:antennas} for a graphical illustration of those terms. One can set up a spectral sequence such that the first differential is the one increasing the sum of the lengths of the antennas. It is easily seen that this complex is acyclic and hence the claim is shown.\footnote{Worries about the convergence of this spectral sequence are adressed in Appendix \ref{sec:specs}.}

\begin{figure}
\centering
\[
\begin{tikzpicture}[
vert/.style={circle,draw,fill, minimum size=5pt, inner sep=0}, invvert/.style={inner sep=-1,minimum size=-1}, ext/.style={circle,draw, minimum size=5pt, inner sep=0},redvert/.style={circle,draw=gray,fill=gray, minimum size=5pt, inner sep=0}, invvert/.style={inner sep=-1,minimum size=-1}, scale=2 ]
\node (v0) at (0.9,1) [vert] {};
\node (v1) at (1.3,0.7) [vert] {};
\node (v2) at (1.9,0.9) [vert] {};
\node (v3) at (1.5,1.4) [vert] {};
\node (v4) at (2,1.6) [vert] {};
\node (v5) at (2.4,1.9)[ redvert] {};
\node (v6) at (2.6,1.4) [redvert] {};
\node (v7) at (1.2,1.7) [redvert] {};
\node (v8) at (0.8,1.8) [redvert] {};
\node (v9) at (2.7,0.9) [redvert] {};
\node (v10) at (1.3,0.4) [redvert] {};
\node (v11) at (1.1,0.2) [redvert] {};
\node (v12) at (0.4,1.6) [redvert] {};
\draw (v0)--(v1) (v1)--(v10) (v1)--(v2) (v0)--(v3) (v3)--(v1) (v3)--(v2) (v2)--(v4) (v4)--(v5) (v4)--(v6) (v7)--(v3);
\draw[gray] (v6)--(v9) (v7)--(v8) (v12)--(v8) (v10)--(v11);
\end{tikzpicture}
\]
 \caption{\label{fig:antennas} The ``antennas'' in the graph are drawn in gray.}
\end{figure}
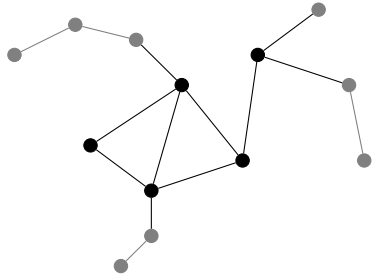

\begin{figure}
\centering
\[
\begin{tikzpicture}[
vert/.style={circle,draw,fill, minimum size=5pt, inner sep=0}, invvert/.style={inner sep=-1,minimum size=-1}, ext/.style={circle,draw, minimum size=5pt, inner sep=0}, scale=2, every loop/.style={} ]
\node (v0) at (0.5,1) [vert] {};
\draw (v0) to [out=45,in=135,loop] ();
\node (v1) at (1.1,0.8) [vert] {};
\node (v2) at (1.1,1.2) [vert] {};
\draw (v1) edge [bend left] (v2);
\draw (v1) edge [bend right] (v2);
\node (v3) at (1.6,0.8) [vert] {};
\node (v4) at (1.6,1.2) [vert] {};
\node (v5) at (1.9,1) [vert] {};
\node (v6) at (2.3,0.8) [vert] {};
\node (v7) at (2.3,1.2) [vert] {};
\node (v8) at (2.7,1.2) [vert] {};
\node (v9) at (2.7,0.8) [vert] {};
\node (v10) at (3.1,1) [vert] {};
\node (v11) at (3.3,1.3) [vert] {};
\node (v12) at (3.3,0.7) [vert] {};
\node (v13) at (3.6,0.8) [vert] {};
\node (v14) at (3.6,1.2) [vert] {};
\draw (v3)--(v4) (v4)--(v5) (v5)--(v3) (v6)--(v7) (v7)--(v8) (v8)--(v9) (v9)--(v6) (v10)--(v11) (v11)--(v14) (v14)--(v13) (v13)--(v12) (v12)--(v10) ;
\end{tikzpicture}
\]
 \caption{\label{fig:loops} Some graphs from the subcomplex of ``loops'', as occuring in the proof of Proposition \ref{prop:GCred}. Note: Some of these graphs (depending on $n$) are zero by symmetry.}
\end{figure}
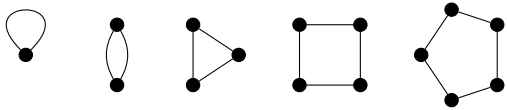

Next we claim that the cohomology of $\GC_n^2$ is 
\[
H(\GC_n^2)=\bigoplus_{ \substack{j\geq 1 \\ j\equiv 2n+1 \mod 4} } \GF[n-j].
\]
Indeed one checks that this is the cohomology of the  subcomplex of $\GC_n^2$ spanned by graphs without trivalent vertices, i. e., by the ``loop'' graphs (see Figure \ref{fig:loops}). Hence the claim reduces to showing that the subcomplex of $\GC_n^2$ spanned by graphs having at least one trivalent vertex is acyclic. For such a graph $\Gamma$, define its \emph{core} to be the graph obtained by deleting all bivalent vertices and, for each deleted vertex, connecting the two adjacent edges. Each edge in the core can be labelled by a natural number $k$, recording how many edges of $\Gamma$ were joined in forming that edge of the core. From the labelled core graph, the original graph $\Gamma$ may be reconstructed. For example:
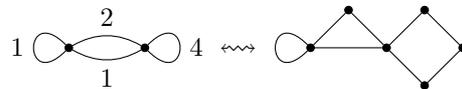
\[
\begin{tikzpicture}[every loop/.style={distance=.8cm}, baseline=-.65ex]
\node[int] (v1) at (0,0) {};
\node[int] (v2) at (1,0) {};
\draw (v1) 	edge[bend left] node[above]{2} (v2)
			edge[bend right] node[below]{1}(v2);
\draw (v1) edge[ out=135, in=-135, loop] node[left]{1} (v1);
\draw (v2) edge[ out=45, in=-45, loop] node[right]{4} (v2);
\end{tikzpicture}
\leftrightsquigarrow
\begin{tikzpicture}[every loop/.style={distance=.8cm}, baseline=-.65ex]
\node[int] (v1) at (0,0) {};
\node[int] (v2) at (1,0) {};
\node[int] (v12) at (.5,.5) {};
\node[int] (v2a) at (1.5,0.5) {};
\node[int] (v2b) at (1.5,-0.5) {};
\node[int] (v2c) at (2,0) {};
\draw (v1) 	edge (v2)
	  (v12) edge (v1) edge (v2)
	  (v2a) edge (v2) edge (v2c)
	  (v2b) edge (v2) edge (v2c);
\draw (v1) edge[ out=135, in=-135, loop] (v1);
\end{tikzpicture}
\]

% Any such graph can be written as a ``core'' with only trivalent vertices, and the edges labelled by natural numbers. An edge labelled by $k$ represents a ``string'' of $k-1$ bivalent vertices. 
One can set up a spectral sequence such that the first differential increases one of the labels on the core by one, i. e., it creates a vertex of valence 2. 
%The resulting complex than has the form
%\[
%\prod_{\Gamma} (\bigotimes_{e\in E(\Gamma)} C )^{G_\Gamma}
%\]
%where the product ranges over all possible core graphs, $G_\Gamma$ is the automorphism group of the core graph $\Gamma$, the tensor product ranges over the edges of $\Gamma$ and the complex $C$ has the form
%\[
%0\leftarrow \GF \stackrel{0}{\leftarrow} \GF
%\stackrel{\id}{\leftarrow} \GF
%\stackrel{0}{\leftarrow} \GF
%\stackrel{\id}{\leftarrow} \GF
%\leftarrow \cdots
%.
%\]
%The group $G_\Gamma$ acts by permutations (with appropriate signs) on the factors of the tensor product, in the same manner it permutes the edges of $\Gamma$. Since the cohomology of $C$ is concentrated in lowest degree 
 More precisely, the differential will increase only any even label and map oddly labelled edges to zero. It is easily checked that the resulting complex is indeed acyclic. Here some care has to be taken in case some labelled edge forms a tadpole, but the argument still works.

\begin{figure}
\centering
\[
\begin{tikzpicture}[scale=1,
vert/.style={draw,outer sep=0,inner sep=0,minimum size=5,shape=circle,fill},
helper/.style={outer sep=0,inner sep=0,minimum size=5,shape=coordinate},
default edge/.style={draw},
 every loop/.style={out=140, in=50, looseness=.8, distance=.8cm }]

\node (v0) at (0.600000023841858,0.799999952316284) [vert] {};
\node (v1) at (2.20000004768372,0.799999952316284) [vert] {};
\node (v2) at (0.600000023841858,0.399999976158142) [helper] {};
\node (v3) at (0.400000005960464,0.399999976158142) [helper] {};
\node (v4) at (0.800000011920929,0.399999976158142) [helper] {};
\node (v5) at (2.20000004768372,0.399999976158142) [helper] {};
\node (v6) at (2,0.399999976158142) [helper] {};
\node (v7) at (2.40000009536743,0.399999976158142) [helper] {};
\node (v14) at (2.20000004768372,1.30000001192093) [vert] {};
\node (v16) at (1.20000004768372,0.799999952316284) [helper,label=0:$\mapsto$] {};

\draw[default edge] (v1)--(v6);
\draw[default edge] (v5)--(v1);
\draw[default edge] (v1)--(v7);
\draw[default edge] (v0)--(v3);
\draw[default edge] (v0)--(v2);
\draw[default edge] (v0)--(v4);
\draw[default edge] (v14)--(v1);
%\draw (v0) to [out=0,in=40,loop] ();
\draw (v0) to [out=70,in=110,loop] ();
%\draw (v0) to [out=140,in=180,loop] ();

%\draw (v1) to [out=0,in=40,loop] ();
\draw (v14) to [out=70,in=110,loop] ();
%\draw (v1) to [out=140,in=180,loop] ();

\end{tikzpicture}
\]
 \caption{\label{fig:tadpolediff} Picture of the part of the differential on $\GC_n^{\geq 3}$ producing a valence three tadpole vertex, as occuring in the proof of Proposition \ref{prop:GCred}.}
\end{figure}
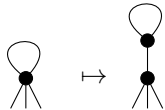

Finally consider $\GC_n^{\geq 3}$. We still need to show that we can omit the tadpoles. One can set up a spectral sequence coming from the filtration according to the number of tadpole vertices of valence three. The first differential hence produces such a vertex, see Figure \ref{fig:tadpolediff}. Similar arguments as for $\GC_n^1$ then show that the cohomology is given by the tadpole-free graphs. This proves the proposition.

\end{proof}

The complexes $\GC_n$ are Kontsevich's (cohomological) graph complexes, see \cite{K4,K5,K3}. We note that the complexes $\GC_n$ for even $n$ are isomorphic, if one reduces the $\mathbb{Z}$ grading to a $\mathbb{Z}_2$ grading. The same holds true for the complexes $\GC_n$ for odd $n$. However, we do not know a way of relating the cohomology of the even $n$ to that of the odd $n$ complexes.

\begin{rem}
\label{rem:onevi}
% As probably\footnote{He announced a proof that probably contains this fact.} noted by M. Kontsevich \cite{K3}, 
As noted by Conant et al. \cite{conant_cut_2005}, %(TODO: CHECK)
the subcomplex of $\GC_n$ spanned by 1-vertex irreducible graphs is quasi-isomorphic to $\GC_n$. We give a short sketch of a proof in Appendix \ref{sec:onevi}.
\end{rem}

%Often in the literature, tadpoles are excluded from the graph complex from the start. We will have use for the tadpole-free version of the graph complex later and introduce the following notation: For $n$ even, $\Gra_n$ is the suboperad of $\Gra_n^-$ spanned by tadpole-free graphs, while $\Gra_n^\whl=\Gra_n^-$ is the full operad. For $n$ odd, $\Gra_n=\Gra_n^\whl$ is simply $\Gra_n^+$. (There are no tadpoles by symmetry reasons anyways.) In general, whenever we draw a symbol ``${}^\whl$'' beneath some graph space, it means that tadpoles are allowed. 

% Similarly to Proposition \ref{prop:GCred} one then proves the following statements.
% 
% \begin{prop}
% \label{prop:GCred2}
% \begin{align*}
% H(\fullGC_n) &\cong H(\Def(\hoLie_n\to \Gra_n))  & \text{for $n$ odd} \\
% H(\fullGC_n) &\cong H(\Def(\hoLie_n\to \Gra_n))\oplus \R[n-1] & \text{for $n$ even} \\
% H(\Def(\hoLie_n\to \Gra_n)\conn) &\cong H(\GC_n)\oplus \bigoplus_{ \substack{j\geq 2 \\ j\equiv 2n+1 \mod 4} } \R[n-j]. &
% \end{align*}
% 
% \end{prop}
% The ``missing'' cohomology class in $H(\Def(\hoLie_n\to \Gra_n))$ for $n$ even is given by the graph with one vertex and a tadpole.

\subsection{Action on polyvector fields and M. Kontsevich's motivation}
For $n\in \Z, r\in \N$ consider the graded commutative polynomial algebra 
\[
A = \K[x_1,\dots, x_r, \xi_1,\dots, \xi_r]
\]
where the $x_j$ have degree $0$ and the $\xi_j$ have degree $n-1$. This algebra becomes an $e_n$ algebra if we impose the relations 
\[
\co{\xi_i}{x_j} = \delta_{ij}.
\]
In fact, the $e_n$ algebra structure is obtained from a $\Gra_n$ algebra structure on $A$, via pull-back along the inclusion $e_n\to\Gra_n$ we encountered above.
The action of an element $\Gamma \in \Gra_n(N)$ determined by a graph which we also denote $\Gamma\in \dgra_{N,k}$ with edge set $E(\Gamma)$ is
\[
\Gamma \cdot(\gamma_1,\dots, \gamma_N)
=
\mu_N \left( \prod_{(i,j)\in E(\Gamma)} \pi_{ij} (\gamma_1\otimes\cdots \otimes \gamma_N) \right).
\] 
Here $\gamma_1,\dots, \gamma_N\in A$, $\mu_N$ is the $N$-fold product and 
\begin{equation}
\label{equ:Graaction}
\pi_{ij} := \sum_{l=1}^r \pderi{}{\xi_l^{(i)}} \pderi{}{x_l^{(j)}} + (-1)^n \pderi{}{\xi_l^{(j)}} \pderi{}{x_l^{(i)}}
\end{equation}
where the superscripts $(i)$ and $(j)$ shall indicate that the derivative acts on the $i$-th, respectively $j$-th factor in the tensor product.
We leave it to the reader to verify that this formula is well defined and that it defines an action of $\Gra_n$ on $A$.

It follows that there is a map of complexes 
\[
\fGC_n \to C(A)
\]
from the full graph complex to the Chevalley-Eilenberg complex of the $\Lie\{n-1\}$ algebra $A$. 
In fact, this map is injective and one may think of $\fGC_n$ as a ``universal version'' (independent of the dimension $r$) of $C(A)$. See \cite{mepolystable} for a more precise form of this statement.

Note also that for $\K=\R$ we may take instead of $A$ the larger algebra
\[
A' = C^\infty(\R^r)[\xi_1,\dots, \xi_r]
\]
and all statements above remain valid. In particular for $n=2$ this is the space of polyvector fields on $\R^r$. 
M. Kontsevich's motivation for considering the graph complex $\GC_2$ in \cite{K3} was the following:
For $n=2$, we obtain a map from $\GC_2$ into the Chevalley-Eilenberg complex of the space of polyvector fields. However, the first cohomology of this complex is the natural recipient for obstructions in the construction of a formality morphism.
If one only consider a special class of tentative formality morphisms, defined through certain graphs (``Kontsevich graphs'') then one may make the obstruction fall into $H^1(\GC_2)$. It is a difficult open conjecture that $H^1(\GC_2)=0$.

\subsection{The (twisted) operads \texorpdfstring{$\Graphs_n$}{Graphsn}}
\label{sec:graphsn}
\begin{figure}
\centering
\[
\begin{tikzpicture}[scale=1,
vert/.style={draw,outer sep=0,inner sep=0,minimum size=5,shape=circle,fill},
helper/.style={outer sep=0,inner sep=0,minimum size=5,shape=coordinate},
default edge/.style={draw},
ext/.style={draw,outer sep=0,inner sep=0,minimum size=5,shape=circle}]

\node (v0) at (0.400000005960464,0.899999976158142) [ext] {};
\node (v1) at (0.900000035762787,0.199999928474426) [ext] {};
\node (v2) at (1.80000007152557,0.899999976158142) [ext] {};
\node (v3) at (3.10000014305115,0.899999976158142) [ext] {};
\node (v4) at (3.79999995231628,0.399999976158142) [ext] {};
\node (v5) at (4.5,1.30000001192093) [ext] {};
\node (v6) at (3.70000004768372,1.30000001192093) [vert] {};
\node (v7) at (3.79999995231628,1) [vert] {};
\node (v8) at (1,1) [vert] {};
\node (v9) at (1.30000007152557,0.5) [vert] {};
\node (v10) at (1.30000007152557,1.39999997615814) [vert] {};
\node (v11) at (0.400000005960464,1.79999999701977) [vert] {};
\node (v12) at (0.600000023841858,1.39999997615814) [vert] {};
\node (v13) at (0.900000035762787,1.79999999701977) [vert] {};
\node (v29) at (4.5,0.699999928474426) [ext] {};

\draw[default edge] (v11)--(v12);
\draw[default edge] (v12)--(v13);
\draw[default edge] (v13)--(v11);
\draw[default edge] (v0)--(v1);
\draw[default edge] (v1)--(v8);
\draw[default edge] (v8)--(v9);
\draw[default edge] (v9)--(v2);
\draw[default edge] (v9)--(v1);
\draw[default edge] (v8)--(v10);
\draw[default edge] (v0)--(v8);
\draw[default edge] (v3)--(v6);
\draw[default edge] (v6)--(v7);
\draw[default edge] (v3)--(v7);
\draw[default edge] (v6)--(v5);
\draw[default edge] (v7)--(v5);
\draw[default edge] (v4)--(v29);
\end{tikzpicture}
\]
 \caption{\label{fig:graphsexample} A typical graph from $\fullGraphs_n(3)$ (left) and one from $\Graphs_n(4)\subset \fullGraphs_n(4)$ (right).}
\end{figure}
In his proof of the formality of the little cubes operads \cite{K2}, M. Kontsevich introduced operads $\Graphs_n$, whose elements are linear combinations of graphs with two kinds of vertices: (i) ``external'' vertices, which are numbered $1,\dots, N$ and (ii) ``internal'' vertices, which are ``indistinguishable'', i. e., unnumbered. 
For a combinatorial description of these operads we refer the reader to \cite[section 3.3.3]{K2}. 
Here we will take a slightly different standpoint and re-define the operads $\Graphs_n$ in a way that makes their relation to the operads $\Gra_n$ we encountered before and to the graph complexes $\GC_n$ more transparent.
Concretely, there is an operation we call \emph{operadic twisting}, which, from an operad $\op P$ with an arrow $\Lie\{k\}\to \op P$, produces in a functorial way another operad $\Tw \op P$, with an arrow  $\Lie\{k\}\to \Tw \op P$. The general theory of operadic twisting is described in Appendix \ref{sec:optwists}, and in more detail in \cite{vasilymetwisting}.
If we apply the twisting functor to the operad $\Gra_n$ (with $k=n-1$) we obtain an operad 
\[
\fullGraphs_n:=\Tw \Gra_n.
\] 
Generators of $\fullGraphs_n(N)$ (as a vector space) can be depicted by graphs with two kinds of vertices: (i) ``external'' vertices, which are numbered $1,\dots, N$ and (ii) ``internal'' vertices, which are ``indistinguishable''. In pictures, we draw external vertices white and internal vertices black. For an example, see Figure \ref{fig:graphsexample}. The operadic composition is obtained by insertion at external vertices, similarly to the composition in $\Gra_n$. One can repeat the construction allowing tadpoles to obtain operads
\[
\fullGraphs_n^\whl:=\Tw \Gra_n^\whl.
\]

\begin{deflemma}
We define $\Graphs_n^\whl\subset\fullGraphs_n^\whl$ be the sub-operad spanned by graphs with all internal vertices at least trivalent and with no connected component consisting entirely of internal vertices. We define $\Graphs_n = \Graphs_n^\whl \cap \fullGraphs_n$ to be the suboperad spanned by graphs without tadpoles.
\end{deflemma}
\begin{proof}
 One needs to show that the subspaces $\Graphs_n^\whl$ and $\Graphs_n$ are closed under the differential and under the operadic compositions.
 However, using the explicit description of the differential and the composition from Appendix \ref{sec:optwists} it is straightforward to check that neither of the forbidden features can be created by insertions or the differential.
\end{proof}

The operad $\Graphs_n$ as in the above definition agrees with the operad $\Graphs_n$ as defined by M. Kontsevich \cite[section 3.3.3]{K2}. 
The theory of operadic twisting guarantees that there is an action of the deformation complex $\Def(\hoLie_{k+1}\to \op P)$ on the twisted operad $\Tw \op P$ by derivations.
In our setting, this means that there is an action of the full graph complex $\fGC_n=\Def(\hoLie_{n}\to \Gra_n)$ on the operad $\fullGraphs_n:=\Tw \Gra_n$.
%Again by the general theory there is a dg action of the dg Lie algebra $\fullGC_n$ on the operad $\fullGraphs_n$. 
%Loosely speaking, the action is the same as the adjoint action on $\fullGC_n$, except that some vertices happen to be colored white.

\begin{lemma}
\label{lem:gcnongraphsn}
The suboperad $\Graphs_n\subset \fullGraphs_n$ is closed under the action of the sub-dg Lie algebra $\GC_n\subset \fGC_n$. 
\end{lemma}
\begin{proof}
A straightforward verification, given the explicit description of the action as in Appendix \ref{sec:optwists_graphs}.
\end{proof}
The action of $\GC_n$ on $\Graphs_n$ will be important later.
%Similarly to the case of graph complexes, we want to simplify the operad $\fullGraphs_n$ a bit. Let $\Graphs_n^\whl$ be the suboperad spanned by graphs with all internal vertices at least trivalent and with no connected component consisting entirely of internal vertices. Let $\Graphs_n$ be the same thing, except that graphs with tadpoles are also excluded.
%\begin{lemma}
% $\Graphs_n^\whl$ and $\Graphs_n$ are dg suboperads.
%\end{lemma}
%\begin{proof}
% One notes that neither of the forbidden things can be created by insertions or the differential.
%\end{proof}
Let us discuss the cohomology of the operads $\Graphs_n$, $\Graphs_n^\whl$, $\fullGraphs_n$ and $\fullGraphs^\whl_n$.
Note that there is a splitting
\[
 \fullGraphs_n(N) = \fullGraphs_{n,c}(N)\otimes (\GF \oplus \fullGC_n[-n])
\]
where $\fullGraphs_{n,c}$ is the suboperad consisting of graphs with no connected components consisting entirely of internal vertices.
For the version with tadpoles there is a similar splitting
\[
 \fullGraphs_n^\whl(N) = \fullGraphs_{n,c}^\whl(N)\otimes (\GF \oplus \fullGC_n^\whl[-n]).
\]
Copying the proof of Proposition \ref{prop:GCred}, one can show the following.
\begin{prop}
\label{prop:Graphsincl}
The inclusions 
\[
\begin{tikzpicture}
\matrix(m)[matrix of math nodes, column sep=.5em]{
& \Graphs_n^\whl & \\
\Graphs_n & & \fullGraphs_{n,c}^\whl \\
& \fullGraphs_{n,c} & \\
};
\path (m-2-1) to node[sloped] {$\subset$} (m-1-2) 
      (m-2-1) to node[sloped] {$\subset$} (m-3-2)
	  (m-1-2) to node[sloped] {$\subset$} (m-2-3)
	  (m-3-2) to node[sloped] {$\subset$} (m-2-3);
\end{tikzpicture}
\]
are quasi-isomorphisms.
\end{prop}

The cohomology of $\Graphs_n$ has been computed by M. Kontsevich \cite{K2} and also P. Lambrechts and I. Volic \cite{LV}.
\begin{prop}[\cite{K2}, \cite{LV}]
\label{prop:Graphscohom}
 The map $e_n\to \Graphs_n$ (defined by the assignments \eqref{equ:enGramapdef}) is a quasi-isomorphism. In particular
\[
 H(\Graphs_n) \cong e_n.
\]
\end{prop}

\begin{rem}
We note in particular that the cohomology of $\Graphs_2^\whl$ is the Gerstenhaber operad $e_2$ and not the Batalin-Vilkovisky (BV) operad. However, if one quotients out from $\Graphs_2^\whl$ the spaces spanned by graphs with a tadpole at an internal vertex, the cohomology is the BV operad. %(TODO: $\exists$ reference?)
\end{rem}

Following \cite{pavol}, one notes that each graph in $\Graphs_n$ decomposes into a (co)product of internally connected components. Here ``internally connected'' means connected after deleting all external vertices. One can hence write 
\begin{equation}
\label{equ:GraphsICG}
 \Graphs_n = S(\ICG_n[1] )
\end{equation}
and similarly
\begin{equation*}
 \Graphs_n^\whl = S(\ICG_n^\whl[1] )
\end{equation*}
where $\CG_n$ (respectively $\ICG_n^\whl$)is spanned by internally connected graphs, shifted in degree by $1$.\footnote{In \cite{pavol} the notation $\mathsf{CG}$ has been used instead of $\CG$. We use the latter notation to prevent confusion with $\GC$.}
One checks that the differential is compatible with the coproduct and hence the spaces $\CG_n$ (respectively $\ICG_n^\whl$) form (operads of) $L_\infty$ algebras.

The following proposition was shown in \cite{pavol} for $n=2$.
\begin{prop}
\label{prop:ICG}
 The cohomologies of the $\CG_n$ and $\ICG_n^\whl$ are the completions of the (operads of) graded Lie algebras $\alg{t}^{(n)}$, where $\alg{t}^{(n)}(N)$ is generated by symbols $t_{ij}=(-1)^{n-1} t_{ji}$, $1\leq i\neq j\leq N$, of degree $2-n$, with relations $\co{ t_{ij} }{ t_{ik}+t_{jk} }=0$ for $\#\{i,j,k\}=3$ and $\co{ t_{ij} }{ t_{kl} }=0$ for $\#\{i,j,k,l\}=4$.
\end{prop}
\begin{proof}
 Copy the proof of the Appendix of \cite{pavol} and change the gradings.
\end{proof}
In the following, we will mostly use the $n=2$ case, which we abbreviate by $\alg{t} := \alg{t}^{(2)}$.

\begin{rem}
 Here it is crucial that we allow multiple edges for $n$ odd. Otherwise it is not true that $\Graphs_n = S(\CG_n[1] )$.
\end{rem}

\section{A spectral sequence for \texorpdfstring{$\Def(\hoe_n\to \op P)$}{Def(hoen -> P)}}
We will often deal with deformation complexes of the form $\Def(\hoe_n\to \op P)$. In this section we introduce a spectral sequence to compute their cohomology. We apply this tool to compute $H(\Def(\hoe_n\to \op P))$ or reduce the computation to a simpler one for several $\op P$. 

\subsection{A grading on \texorpdfstring{$e_n$}{en} and a filtration on \texorpdfstring{$\Def(\hoe_n\to \op P)$}{Def(hoen-> P)}}
\label{sec:enfiltration}
It is a well known fact that the symmetric sequence $e_n$ may be written as 
\[
e_n = \Com \circ \Lie\{n-1\}
\]
where $\circ$ is the monoidal product in the category of symmetric sequences, i.e., the ``plethysm'' operation.
More concretely, this means that $e_n(N)$ is a direct sum of spaces of the form
\begin{equation*}
%\label{}
%e_n(N) \cong \oplus_{j=1}^N \oplus_{\substack{ i_1,\dots, i_j \\ i_1+\cdots +i_j=N }} 
\Com(k) \otimes \Lie\{n-1\}(i_1)\otimes \cdots \otimes \Lie\{n-1\}(i_k)
\end{equation*}
where $k\in \{1,\dots, N\}$ and $i_1,\dots, i_k$ are such that $i_1+\cdots +i_k=N$. 
There is a grading on $e_n$ by the number $k-1$. 
It can be checked to be compatible with the operad structure.

\begin{rem}
On any graded operad (i. e., one for which the differential is zero) there exist two ``trivial'' gradings, one by the arity minus one, and one by the cohomological degree. The grading we consider here has the form
\[
(\text{arity})-1 + (\text{cohomological degree})/(n-1).
\]
\end{rem}

Dually we obtain a grading on the cooperads $e_n^*$ and furthermore on the cooperads $e_n^*\{n\}=e_n^\vee$ Koszul dual to $e_n$. 
This grading also transfers to the operads $\hoe_n =\Omega(e_n^\vee)$.
Now let $\op P$ be some other operad and
\[
f: \hoe_n \to \op P
\]
an operad map. Let us consider $\op P$ to be concentrated in degree zero with respect to our additional grading. Of course $f$ will then not respect the additional grading in general. However, since $\Omega(e_n^\vee)$ is concentrated in non-negative degrees we still inherit a descending filtration $\mF$ on the complex
\[
\Def(\hoe_n \to \op P).
\]
We will use the associated spectral sequence as a tool to compute $H(\Def(\hoe_n \to \op P))$ in some cases below.

\subsection{The special case \texorpdfstring{$\op P=\Graphs_n^\whl$}{P=Graphs n whl}}
Consider $\op P=\Graphs_{n}^\whl$. 
The goal of this subsection is to identify a certain subcomplex of $\Def(\hoe_n\to\Graphs_{n}^\whl)$
which is quasi-isomorphic to the full complex.
\begin{deflemma}
\label{deflem:fC}
We define $\fC\subset \Def(\hoe_n\to\Graphs_{n}^\whl)$ to be the subcomplex spanned by elements $x\in \Def(\hoe_n\to\Graphs_{n}^\whl)$ such that
\begin{enumerate}
\item The image of all components of $x$ (which are maps in $\Hom_{S_N}(e_n^\vee(N), \Graphs_{n}^\whl(N))$, $N=1,2,\dots$) is in the space of linear combinations of graphs all of whose external vertices have valence exactly 1.
\item The components of $x$  vanish on all of  $e_n^\vee(N)$ except possibly on the subspace of bottom cohomological degree $(1-N)n$. (Or put differently, the components of $x$ descend to the quotient $\Com^*\{n\}$.)
\end{enumerate}
\end{deflemma}

Interpreting elements of $\Def(\hoe_n\to\Graphs_{n}^\whl)$ combinatorially as linear combinations of certain graphs as explained in Appendix \ref{sec:defengra}, $\fC$ may alternatively be defined as the subspace spanned by graphs with all clusters of size one, and all external vertices of valence one. The differential has two pieces, one splits internal vertices, and one attaches a new external vertex of valence one to some internal vertex.
%TODO: picture?

\begin{proof}[Proof of Definition/Lemma \ref{deflem:fC}.]
We need to show that the subspace $\fC$ defined above is indeed a subcomplex.
The differential on $\Def(\hoe_n\to\Graphs_{n}^\whl)$ splits as in \eqref{equ:diffsplit2}. It is clear that $\delta$ may not invalidate the second condition of the definition. Furthermore, $\delta$ at most decreases valences of vertices, and cannot produce new valence one vertices, so it also cannot invalidate the first condition. 
The effect of the other terms in the differential, i.~e., of the bracket with the Maurer-Cartan element $\alpha$ is most easily discussed using the graphical interpretation of elements of $\Def(\hoe_n\to\Graphs_{n}^\whl)$ as in Appendix \ref{sec:defengra}. 
Also note that $\co{\alpha}{\cdot}=\alpha\circ\cdot \pm \cdot \circ \alpha$ where $\circ$ is the pre-Lie product in $\Def(\hoe_n\to\Graphs_{n}^\whl)$. Both operations $\alpha_\wedge\circ\cdot$ and $\cdot \circ \alpha_\wedge$ will create graphs with one cluster of size two from a graph with all clusters of size one. However, since the external vertices have valence one, the terms produced are identical and cancel each other. Hence $\co{\alpha_\wedge}{\cdot}$ acts trivially on $\fC$. The operation $\co{\alpha_{\co{}{}}}{\cdot}$ adds one cluster of size one, and does not change cluster sizes otherwise, hence it automatically preserves the second property demanded for in the definition. However, it might a priori invalidate the first. The violating terms are graphs with one external vertex of valence two, connected to some other external vertex. But again, these terms are produced in pairs, one by $\alpha_{\co{}{}}\circ\cdot$ and one by $\cdot \circ \alpha_{\co{}{}}$, and cancel.
\end{proof}

The main result of this subsection is the following.
\begin{prop}
\label{prop:fCqiso}
The inclusion $\fC\to \Def(\hoe_n\to\Graphs_{n}^\whl)$ is a quasi-isomorphism.
\end{prop}

The proposition will follow from the following Lemma.
\begin{lemma}
\label{lem:fCqiso}
The operation $d_\wedge$ squares to zero and the inclusion
\[
(\fC, 0) \to (\Def(\hoe_n\to\Graphs_{n}^\whl), d_\wedge)
\]
is a quasi-isomorphism, where the notation means that the left hand side is considered a complex with trivial differential and the right hand side is considered a complex with differential $d_\wedge$. Here $d_\wedge$ is as in \eqref{equ:dwedge}.
\end{lemma}
\begin{proof}
The first assertion follows from the Maurer-Cartan equation
\[
\delta \alpha + \frac 1 2 \co{\alpha}{\alpha}=0
\]
with $\alpha=\alpha_\wedge+\alpha_{\co{}{}}$ (as in \eqref{equ:diffsplit2}) by restriction.
For the second assertion we have to compute the cohomology of $V:=(\Def(\hoe_n\to\Graphs_{n}^\whl), d_\wedge)$. 
Let the \emph{character} of a graph $\Gamma$ with external and internal vertices be the isomorphism class of the graph obtained by deleting the external vertices but keeping the dangling edges. The differential on $V$ does not alter the character and hence $V$ splits into a direct product of subcomplexes, one for each character. Call the subcomplex $V_c$ for the fixed character $c$, say with $k$ dangling edges. To $c$ is associated some automorphism group $G$ acting by permutations on the dangling edges. Concretely, $V_c$ is isomorphic to the complex
\[
\oplus_{p=1}^k V_{p,k,n}^G
\]
up to an overall degree shift depending on $c$. Here $V_{p,k,n}$ is as in Appendix \ref{sec:harrison}. It then follows from Lemma \ref{lem:harrisoncohom} that $V_c$ has one dimensional cohomology, and tracing the representative, one sees that this is an element in the (one dimensional) intersection of $V_c$ with the image of the inclusion $(\fC, 0) \to (\Def(\hoe_n\to\Graphs_{n}^\whl), d_\wedge)$. The statement of the lemma immediately follows.
\end{proof}

In the following proofs we will use some spectral sequence arguments. The following lemma will settle convergence.
\begin{lemma}
\label{lem:findimdecomp}
$\Def(\hoe_n\to\Graphs_{n}^\whl)$ decomposes into a direct product of subcomplexes which are finite dimensional in each degree.
\end{lemma}
\begin{proof}
We put an additional filtration on $\Graphs_{n}^\whl$ by assigning to a graph with $v$ internal vertices and $e$ edges the degree $v-e$. The differential on $\Graphs_{n}^\whl$ preserves this degree.
Furthermore we assign a co-Gerstenhaber word in $e_n^\vee(N)$ with $b$ cobrackets the degree $b-N$. This yields a grading on each factor in the direct product \eqref{equ:convalgebra}. 
One checks that the differential on $\Def(\hoe_n\to\Graphs_{n}^\whl)$ does not alter the additional degree, and hence the direct product of the (additional) degree $k$ components of all factors yields a subcomplex. In each such subcomplex, fixing the cohomological degree fixes the total number of vertices and bounds the number of edges in graphs occurring. This leaves a finite dimensional subspace.
Note also that using the graphical interpretation of elements of $\Def(\hoe_n\to\Graphs_{n}^\whl)$ from Appendix \ref{sec:defengra} the additional degree is the total number of vertices minus the total number of edges, i. e., the Euler characteristic. 
\end{proof}

The following technical proposition will be used in later sections.
\begin{prop}
\label{prop:fCqiso2}
The operation %$\delta + \co{\alpha_\wedge}{\cdot}$ 
$\delta + d_\wedge$ 
squares to zero and the inclusion
\[
(\fC, \delta) \to (\Def(\hoe_n\to\Graphs_{n}^\whl), \delta + d_\wedge)
\]
is a quasi-isomorphism.
\end{prop}
\begin{proof}
Take a filtration on the mapping cone by the cohomological degree in $\Graphs_{n}^\whl$. The associated graded is the mapping cone of the map in Lemma \ref{lem:fCqiso} and hence acyclic.
The Lemma then follows from Lemma \ref{lem:findimdecomp}. More explicitly: The mapping cone is a direct product of subcomplexes as in Lemma \ref{lem:findimdecomp}. On each subcomplex the filtration is bounded and hence the spectral sequence converges to cohomology. Hence each subcomplex is acyclic. Hence their direct product, i.~e., the mapping cone is acyclic.
\end{proof}

\begin{rem}
\label{rem:tildef}
Let $\tilde f:e_n\to \Graphs_{n}^\whl$ be the operad map mapping the product generator to the empty graph as in \eqref{equ:enGramapdef}, but mapping the bracket to zero. Let further $\tilde f'$ be the composition $\hoe_n\to e_n \stackrel{\tilde f} \longrightarrow \Graphs_{n}^\whl$.
Then
\[
(\Def(\hoe_n\to\Graphs_{n}^\whl), \delta + d_\wedge)
\cong 
\Def(\hoe_n\stackrel{\tilde f} \longrightarrow \Graphs_{n}^\whl)
\]
may be understood as a deformation complex as well.
\end{rem}

\begin{proof}[Proof of Proposition \ref{prop:fCqiso}.]
We have to show that the mapping cone is acyclic. To do this consider the spectral sequence coming from the filtration of section \ref{sec:enfiltration}. The differential on the associated graded of $\Def(\hoe_n\to\Graphs_{n}^\whl)$ has the form $\delta + d_\wedge$. 
The differential on the associated graded of $\fC$ is just $\delta$. 
Proposition \ref{prop:fCqiso2} then says that the associated graded of the mapping cone is acyclic. It again follows by Lemma \ref{lem:findimdecomp} that the mapping cone is acyclic.
%
%The differential on the associated graded of the mapping cone has the form 
%\[
%\delta + \co{\alpha_\wedge}{\cdot} + d_\iota
%\]
%where $\delta$ is supposed to act on both pieces $\fC$ and $\Def(\hoe_n\to\Graphs_{n}^\whl)$, $\co{\alpha_\wedge}{\cdot}$ is zero on $\fC$ (cf. the proof of Definition/Lemma \ref{deflem:fC}) and $d_\iota$ is the part of the differential on the mapping cone coming from the inclusion map.
%We claim that the associated graded of the mapping cone is acyclic. To show this, consider yet another filtration by the cohomological degree of elements of $\Graphs_{n}^\whl$. On the associated graded, this will reduce the differential to 
%\[
%\co{\alpha_\wedge}{\cdot} + d_\iota.
%\]
%We claim that this differential renders the associated graded (of the associated graded) of the mapping cone acyclic. But this is the statement of Lemma \ref{lem:fCqiso}.
%TODO: discuss convergence of spec. seq.
\end{proof}

\subsubsection{The connected part}
Let us again use the graphical interpretation of elements $\Def(\hoe_n\to\Graphs_{n}^\whl)$ from Appendix \ref{sec:defengra}. 
Let 
\[
\Def(\hoe_n\to\Graphs_{n}^\whl)_\conn \subset \Def(\hoe_n\to\Graphs_{n}^\whl)
\]
be the subcomplex generated by the connected graphs.
Any graph decomposes into its connected components and the differential may combinatorially neither glue nor disconnect components. Hence
\[
 \Def(\hoe_n\to\Graphs_{n}^\whl) \cong S^+(\Def(\hoe_n\to\Graphs_{n}^\whl)_\conn)
\]
as complexes. Similarly let $\fC_\conn\subset \fC$ be the 
intersection of $\fC$ and $\Def(\hoe_n\to\Graphs_{n}^\whl)_\conn$.
Alternatively, it may be described as the space of maps which take values in $\Graphs_n^\whl$ which are linear combinations of connected graphs.
The complex $\fC$ splits into a symmetric product of its connected subcomplex
\[
\fC = S^+(\fC_\conn)
\]
Furthermore, checking the proof of Propositions \ref{prop:fCqiso} and \ref{prop:fCqiso2} we see that the result restricts to the connected parts.
\begin{prop}
\label{prop:fCconnqiso}
The inclusions
\[
(\fC_\conn, \delta) \to (\Def(\hoe_n\to\Graphs_{n}^\whl)_\conn, \delta + d_\wedge)
\]
and
\[
\fC_\conn \to \Def(\hoe_n\to\Graphs_{n}^\whl)_\conn
\]
are quasi-isomorphisms.
\end{prop}
%
%and accordingly
%\[
%H(\fC) = S^+(H(\fC_\conn)).
%\]
%
%$\Def(\hoe_n\to\Graphs_{n}^\whl)$
%There is a subcomplex 

%\begin{rem}
%\label{rem:fCsplit}
%Let $\fC_\conn\subset \fC$ be the subcomplex of $\fC$ consisting of maps which takes values in $\Graphs_n^\whl$ which are linear combinations of connected graphs.
%Note that the complex $\fC$ splits into a symmetric product of its connected subcomplex
%\[
%\fC = S^+(\fC_\conn)
%\]
%and accordingly
%\[
%H(\fC) = S^+(H(\fC_\conn)).
%\]
%\end{rem}

\subsubsection{A closer look at $\fC_\conn$}
\label{sec:closerlook}
This subsection contains some technical calculations, the motivation for which will be given in section \ref{sec:derhoe2grt} below.
Let $\alg{g}$ be an operad in $L_\infty$ algebras. Its (completed) Chevalley complexes 
\[
C(\alg g)(N) := C(\alg{g}(N))\cong \prod_{k\geq 0} S^k(\alg{g}(N)[1])
\]
form an operad $C(\alg g)$ of dg vector spaces. There is always a map
\[
\Com\to C(\alg g)
\]
by sending the product operation to the basis element of $\GF \cong S^0(\alg{g}(2)[1]) \subset C(\alg g)(2)$. Also, there is a map 
\[
e_n\to C(\alg g)
\]
by pre-composing with the map $e_n\to \Com$ that sends the bracket operation to zero and the product operation to the product operation.  
Denote the composition $\Com_\infty\to \Com \to C(\alg g)$ by $f$ and the composition $\hoe_n\to e_n\to C(\alg g)$ by $g$. Then we may form the complexes
\begin{align*}
\Def(\Com_\infty\stackrel{ f} \longrightarrow C(\alg{g}))
&\cong
\prod_{N\geq 1} \prod_{k\geq 0} \Hom_{S_N}(\Com^\vee(N), S^k(\alg{g}(N)[1])) \\
\Def(\hoe_n\stackrel{ g} \longrightarrow C(\alg{g}))
&\cong
\prod_{N\geq 1} \prod_{k\geq 0} \Hom_{S_N}(e_n^\vee(N), S^k(\alg{g}(N)[1])).
\end{align*}
One checks that the subspaces
\[
\prod_{N\geq 1}  \Hom_{S_N}(\Com^\vee(N), \alg{g}(N)[1])
\] 
and 
\[
\prod_{N\geq 1}  \Hom_{S_N}(e_n^\vee(N), \alg{g}(N)[1])
\] 
are subcomplexes, which we denote (by abuse of notation)
\[
\Def(\Com_\infty \rightarrow \alg{g}[1])
\]
and
\[
\Def(\hoe_n \rightarrow \alg{g}[1]).
\]

One example is $\alg{g}=\ICG_n^\whl$, for which 
\[
\Graphs_n^\whl(N) \cong C(\ICG_n^\whl(N)).
\] 
 Note that the complex $\fC_\conn$ from above is a subcomplex of $\Def(\hoe_n \rightarrow \ICG_n^\whl[1])$.
\begin{prop}
\label{prop:fCICG}
The inclusions
\[
(\fC_\conn, \delta) \subset (\Def(\hoe_n \rightarrow \ICG_n^\whl[1]), \delta + d_\wedge)
\]
and
\[
\fC_\conn\subset \Def(\hoe_n \rightarrow \ICG_n^\whl[1])
\]
are quasi-isomorphisms.
\end{prop}
\begin{proof}
The proof is a copy of the proofs of Propositions \ref{prop:fCqiso}, \ref{prop:fCqiso2}.
\end{proof}

Another interesting example is the (completed) Drinfeld-Kohno operad of Lie algebras $\hat{\alg{t}}$. Recall from \cite{pavol} that it is quasi-isomorphic to $\ICG_2$, and hence to $\ICG_2^\whl$. Using this fact one may show the following result.

\begin{prop}
\label{prop:icgandt}
$H(\Def(\hoe_2 \rightarrow \ICG_2^\whl[1]))\cong H(\Def(\hoe_2 \rightarrow \alg{t}[1]))$
and 
$H(\Def(\hoe_2 \rightarrow \ICG_2^\whl[1]), \delta+\co{\alpha_\wedge}{\cdot})\cong H(\Def(\hoe_2 \rightarrow \alg{t}[1]),d_\wedge )$.
\end{prop}
\begin{proof}
There is a zig-zag of quasi-isomorphisms 
\[
\hat{\alg{t}} \leftarrow \TCG^\whl \to \ICG_2^\whl.
\]
where $\TCG^\whl$ (the truncated version of $\ICG_2^\whl$) is identical to $\ICG_2^\whl$ in degrees $<0$, is zero in degrees $>0$ and the closed subspaces of $\ICG_2^\whl$ in degree 0. The right hand arrow is the inclusion, while the left hand arrow is the projection of the degree 0 part to the cohomology.
We may hence form a zig-zag of complexes
\[
\Def(\hoe_2 \rightarrow \hat{\alg{t}}[1])\leftarrow
\Def(\hoe_2 \rightarrow \TCG^\whl[1])\to
\Def(\hoe_2 \rightarrow \ICG_2^\whl[1])
\]
Put a filtration on the mapping cones of each map by the total cohomological degree minus the cohomological degrees of $\TCG^\whl$ and $\ICG_2^\whl$. The associated graded is acyclic, hence by arguments similar to those in the proof of Lemma \ref{prop:fCqiso2} the mapping cones are acyclic.
\end{proof}

%
%
%We restrict now to the case $n=2$.\footnote{This is not for technical reasons, but because we want to use some results of \cite{pavol}, and only the $n=2$ case is considered there. }
%Recall from \cite{pavol} that the operad of $L_\infty$ algebras $\ICG_2$ is quasi-isomorphic to the Drinfeld-Kohno operad of Lie algebras $\alg{t}$.
%
%
%
%
% via a zig-zag of quasi-isomorphisms 
%\[
%\alg{t} \leftarrow \TCG \to \ICG.
%\]

%There is an increasing filtration on the operad $\Graphs_n^\whl$ such that 
%\[
%\mF^p \Graphs_n^\whl(N) \cong 
%\prod_{k= 0}^p S^k(\ICG_n^\whl(N)[1]).
%\]

\subsection{The special case \texorpdfstring{$\op P=\Gra_n^\whl$}{P=Gra n whl}}
The operad $\Gra_n^\whl$ is a quotient of $\Graphs_n^\whl$ by sending to zero all graphs with internal vertices. 
The corresponding quotient of $\fC$ is rather trivial, being either one-dimensional ($n$ even) or a symmetric product of a one dimensional space.
One may copy verbatim the proofs of Propositions \ref{prop:fCqiso} and \ref{prop:fCconnqiso} to arrive at the following result.
\begin{prop}
\label{prop:hoengracohomfull}
The cohomologies of $\Def(\hoe_n\to \Gra_n^\whl)_\conn$ and $\Def(\hoe_n\to \Gra_n^\whl)$ satisfy
 %The cohomology of $\Def(\hoe_n\to \Gra_n^\whl)_\conn$ is 
\begin{align*}
 H(\Def(\hoe_n\to \Gra_n^\whl)_\conn) &\cong \GF[-1] \\
 H(\Def(\hoe_n\to \Gra_n^\whl)) &\cong S^+(\GF[-n-1])[n].
\end{align*}
%and \emph{for even} $n$, the cohomology of its tadpole free part is
%\[
% H(\Def(\hoe_n\to \Gra_n)_\conn) \cong \R[-1]\oplus\R[n-2]
%\]
\end{prop}
\begin{rem}
 The cohomology class $\GF[-1]$ corresponds to a relative rescaling of the Lie bracket and product. 
It corresponds to a graph with two vertices (each in its own cluster) and one edge. Elements of the symmetric product correspond to unions of such graphs. This is only possible for $n$ odd by symmetry reasons, otherwise the cohomology is just one-dimensional.
 %The class $\R[n-2]$ encodes a map sending the Lie co-bracket to the Lie bracket. The corresponding graph has two vertices in a cluster and an edge connecting them. In the complex with tadpoles this class is the coboundary of the ``divergence'' operation, i.e., the graph with one vertex and a tadpole.
\end{rem}
\begin{rem}
 For $n=2$ this result is the ``universal'' version of the well-known statement that the Gerstenhaber algebra $T_\poly$ of polyvector fields on $\R^n$ is non-deformable, up to homotopy.
The operads $\Gra_2$ or $\Gra_2^\whl$ can be seen as universal versions of the operad $\End(T_\poly)$.
\end{rem}

\section{The map between 
\texorpdfstring{$\GC_n$}{GCn} 
and 
\texorpdfstring{$\Der(\hoe_n)$}{Der(hoen)}
}
\label{sec:themap}
In this section we will define the map between $H(\GC_n)$ and $H(\Der(\hoe_n))$. In fact, we know two distinct ways to describe that map, and it is advantageous to know both.
We will give the first definition in section \ref{sec:firstdef} and the second definition in section \ref{sec:seconddef} and show that the two maps thus defined are identical in section \ref{sec:mapsagree}.

\subsection{Reduction to the connected part}
Any graph may be split into its connected components. In section \ref{sec:grcomplexes} we saw that the differential on the full graph complex $\fGC_n$ respects this splitting, and $\fGC_n$ can be written as a symmetric product of the subcomplex $\fGC_{n,\conn}$ spanned by connected graphs.
The goal of this subsection is to describe an analogous, but less obvious notion of ``connectedness'' for some basis elements of the complex $\Def(\hoe_n \to e_n)$. 
Concretely, we will define a subcomplex 
\[
\Def(\hoe_n \to e_n)_\conn \subset \Def(\hoe_n \to e_n)
\]
which we call the ``connected part'' of $\Def(\hoe_n \to e_n)$  such that 
\[
 \Def(\hoe_n \to e_n) \cong S^+(\Def(\hoe_n \to e_n)_\conn [-n])[n]
\]
as complexes. Here $S^+(\cdots)$ denotes the completed symmetric product space, without the term $\GF$, as in \eqref{equ:Splus}.
To define the connected part, recall the well known fact that 
\[
e_n \cong \Com \circ \Lie\{n-1\}
\]
as $\bbS$-modules, where ``$\circ$'' is the usual monoidal product in the category of $\bbS$-modules, see \cite[section 5.1.6]{lodayval}.
This may be re-stated by saying that operations in $e_n(N)$ can be identified with linear combinations of expressions of the form 
\begin{equation}
%\label{equ:coenform}
 L_1(X_{1},\cdots,X_{N})\wedge \cdots \wedge L_k(X_{1},\cdots, X_{N})
\end{equation}
where $X_1,\dots,X_N$ are formal variables, $L_j$ are $\Lie\{n-1\}$-words and each $X_i$ occurs exactly once in the expression. The action of the symmetric group $S_N$ on $e_n(N)$ is given by permuting the indexes of the formal variables $X_1,\dots, X_N$. Picking some set of $\Lie\{n-1\}$-words forming a basis for $\Lie\{n-1\}(m)$ for each $m$ (for example the one from \cite[section 13.2.5]{lodayval}) a basis of $e_n(N)$ by a set of expressions of the above form may be written down.

%recall the well known fact that 
%\[
%e_n \cong \Com \circ \Lie\{n-1\}
%\]
%as symmetric sequences, where ``$\circ$'' is the usual monoidal product in the category of symmetric sequences. More concretely, this means that
%\[
%e_n(N) \cong \bigoplus_{1\leq k \leq N}\left( \Com(k)\otimes \bigoplus_{\substack{k_1,\dots, k_k \\ k_1+\cdots+k_k=N}} \Lie\{n-1\}(k_1)\otimes \cdots \otimes  \Lie\{n-1\}(k_k) \otimes \GF [S_N] \right)_{G_{N, k, k_1,\dots, k_k}}
%\]
%where $\GF[S_N]$ is the group algebra of the symmetric group and $G_{N, k, k_1,\dots, k_k}$ is a permutation group.
%Alternatively, we may identify $e_n(N)$ with the space of linear combinations of two level trees
%
%
%
% The only important fact for us is that $e_n$ inherits a filtration by the number $k$ in the above sum, i. e.,
%\[
%\mF^p e_n(N) \cong \bigoplus_{1\leq k \leq p}\left( \Com(k)\otimes  \Lie\{n-1\}(k_1)\otimes \cdots \otimes  \Lie\{n-1\}(k_k) \otimes \GF[S_N] \right)_{G_{N, k, k_1,\dots, k_k}}.
%\] 

Next consider again
\[
\Def(\hoe_n\to e_n) = \prod_{N\geq 1} \Hom_{S_N}( e_n^*\{n\}(N), e_n(N) )
\cong
\prod_{N\geq 1} \left( e_n(N)\otimes e_n(N)\otimes \sgn^{\otimes n}[(1-N)n] \right)^{S_N}.
\]
Similarly to what was said above, we can define a basis of $e_n(N)\otimes e_n(N)$ given by a set of expressions of the form 
\begin{equation}
\label{equ:enen}
\phi:=(L_1(X_{1},\cdots,X_{N})\wedge \cdots \wedge L_k(X_{1},\cdots, X_{N}))
\otimes 
(L_1'(Y_{1},\cdots,Y_{N})\wedge \cdots \wedge L_{k'}'(Y_{1},\cdots, Y_{N}))
\end{equation}
where $X_1,\dots,X_N, Y_1,\dots, Y_N$ are formal variables and $L_1, \dots , L_{k'}'$ are $\Lie\{n-1\}$-words such that each formal variable occurs exactly once in the expression. 
We say that two $\Lie\{n-1\}$-words $L(X_{1},\cdots,X_{N})$ and $L'(Y_{1},\cdots,Y_{N})$ \emph{share a variable} if there is some integer $j$ such that $X_j$ occurs in $L$ and $Y_j$ occurs in $L'$. For example $\co{X_1}{X_2}$ and $\co{Y_5}{Y_2}$ share a variable while $\co{X_1}{X_2}$ and $\co{Y_5}{Y_4}$ do not. Given some formal expression $\phi$ of the form \eqref{equ:enen} we say that Lie words $L_i$ and $L_j'$ are \emph{connected} in this expression if there is an ordered set of integers $j_1, i_1,\dots, j_r, i_r$ ($r$ is some other integer) such that $L_i$ and $L_{j_1}'$ share a variable, $L_{i_1}$ and $L_{j_1}'$ share a variable, $L_{i_1}$ and $L_{j_2}'$ share a variable etc., and finally $L_{i_r}$ and $L_j'$ share a variable. Similarly, we say that $L_{i}$ and $L_j$ (resp. $L_{i}'$ and $L_j'$) are connected if $i=j$ or if $L_{i}$ and $L_j$ (resp. $L_{i}'$ and $L_j'$) if there is a sequence of integers $j_1, i_1,\dots, j_r$ (resp. $i_1, j_1,\dots, i_r$) such that $L_{i}$ shares a variable with $L_{j_1}'$, which shares a variable with $L_{i_1}$ etc. 
Clearly connectedness defines an equivalence relation on the set $\{L_1,\dots, L_{k'}' \}$ of  $\Lie\{n-1\}$-words occurring in the expression $\phi$.  We define the \emph{number of connected components} of $\phi$ as the number of equivalence classes under this equivalence relation.
We say that $\phi$ is \emph{connected} if its number of connected components is one.
\begin{ex}
The expression 
\[
(\co{X_1}{X_2}\wedge \co{X_3}{X_4})
\otimes 
(Y_1\wedge\co{Y_2}{Y_3}\wedge Y_4)
\]
is connected while the expression 
\[
(\co{X_1}{X_2}\wedge \co{X_3}{X_4})
\otimes 
(\co{Y_1}{Y_2}\wedge Y_3 \wedge Y_4)
\]
has two connected components.
\end{ex}
Picking a basis of $e_n(N)\otimes e_n(N)$ by expressions of the form \eqref{equ:enen} we may define the subspace 
\[
V_N \subset e_n(N)\otimes e_n(N)
\]
spanned by the connected basis elements.
 $V_N$ does not depend on the particular choice of basis and it is closed under the action of $S_N$.
We now define 
\[
\Def(\hoe_n \to e_n)_\conn = 
\prod_{N\geq 1} \left( V_N\otimes \sgn^{\otimes n}[(1-N)n] \right)^{S_N}
\subset \Def(\hoe_n \to e_n).
\]

\begin{prop}
\label{prop:defenbyconnectedpart}
$\Def(\hoe_n \to e_n)_\conn \subset \Def(\hoe_n \to e_n)$ is a sub-complex and furthermore 
\[
 \Def(\hoe_n \to e_n) \cong S^+(\Def(\hoe_n \to e_n)_\conn [-n])[n]
\]
as complexes, where $S^+(\cdots)$ is the completed symmetric product space, see \eqref{equ:Splus}. In particular
\[
 H(\Def(\hoe_n \to e_n)) \cong S^+(H(\Def(\hoe_n \to e_n)_\conn) [-n])[n].
\]
\end{prop}
\begin{proof}
A straightforward verification.
\end{proof}

In order to prove Theorem \ref{thm:GCdef} it suffices to show below the following result
\begin{prop}
\label{prop:fgcconndefenconn}
There is an injective map $H(\fGC_{n,\conn})\to H(\Def(\hoe_n \to e_n)_\conn)$ with one-dimensional co-kernel, concentrated in degree 0.
\end{prop}

%\subsubsection{A generalization}

\subsection{The first description of the map}
\label{sec:firstdef}
There are canonical quasi-isomorphisms 
\[
\hoe_n\to e_n \to \Graphs_n\to \fGraphs_{n,c}^\whl,  
\]
and hence also quasi-isomorphisms of complexes\footnote{Our conventions regarding deformation complexes are stated in section \ref{sec:defcomplexes}.}
\begin{equation}
\label{equ:defentodefgraphs}
 \Der(\hoe_n)[-1] \to \Def(\hoe_n\to e_n) \to \Def(\hoe_n\to \Graphs_n)\to \Def(\hoe_n\to \fGraphs_{n,c}^\whl).
\end{equation}
Elements of the complex on the right can be written as linear combinations of graphs, with certain symmetry properties. See Appendix \ref{sec:defengra} for a more detailed description. Any such graph splits into a union of its connected components, and the differential acts on each component separately and preserves connectedness. Let $\Def(\hoe_n\to \fGraphs_{n,c}^\whl)_\conn\subset \Def(\hoe_n\to \fGraphs_{n,c}^\whl)$ be the subspace spanned by the connected graphs.
Then 
\[
\Def(\hoe_n\to \fGraphs_{n,c}^\whl)
\cong
S^+(\Def(\hoe_n\to \fGraphs_{n,c}^\whl)_\conn[-n])[n]
\]
as complexes.
% furthermore the map \eqref{equ:defentodefgraphs} is induced by a map
%\[
%\Def(e_n)_\conn [-1] \to \Def(\hoe_n\to \fGraphs_{n,c}^\whl)_\conn.
%\]
%If follows that this map is a quasi-isomorphism as well.
The dg Lie algebra $\fGC_{n, \conn}^\whl$ acts canonically on $\fGraphs_{n,c}^\whl$ by operadic derivations, see Lemma \ref{lem:gcnongraphsn}.
By composition with the map $\hoe_n\to \fGraphs_{n,c}^\whl$ one obtains from any derivation of $\fGraphs_{n,c}^\whl$ an element in $\Def(\hoe_n\to \fGraphs_{n,c}^\whl)$.
Derivations coming from elements of $\fGC_{n, \conn}^\whl$ in fact are mapped into the connected part $\Def(\hoe_n\to \fGraphs_{n,c}^\whl)_\conn$.
Hence there is a map of complexes
\[
\Phi\colon \fGC_{n, \conn}^\whl \to \Def(\hoe_n\to \fGraphs_{n,c}^\whl)_\conn[1].
\]
%Theorem \ref{thm:GCdef} (to be proven) states that this map is a quasi-isomorphism, up to one term.
%For the following proof it will be slightly more economic to allow tadpoles in our graphs, i.e., to work with $\Graphs_n^\whl$ instead of $\Graphs_n$. Since the inclusion $\Graphs_n\to \Graphs_n^\whl$ is a quasi-isomorphism, this is not a severe change.

Let us denote by $\bar \Phi$ the composition 
\[
H(\fGC_{n, \conn}^\whl) \to
H(\Def(\hoe_n\to \fGraphs_{n,c}^\whl)_\conn)[1] \to %\stackrel{\cong}{\to}
H(\Der(\hoe_n))
\]
where the first map is induced by $\Phi$. We will see below that $\bar\Phi$ is injective with one dimensional cokernel, in degree 0.
For now, let us show the following result.
\begin{prop}
\label{prop:firstmaplie}
$\bar \Phi$ is a map of Lie algebras. 
\end{prop}
\begin{proof}
%By the action of $\fGC_{n, \conn}^\whl$ on $\fGraphs_{n,c}^\whl$ we have a map of dg Lie algebras
%\[
%\fGC_{n, \conn}^\whl \to \Der'(\fGraphs_{n,c}^\whl)
%\]
%using the notation of section\ref{sec:defcomplexes} and $\bar \Phi$ factors through the map induced on cohomology by this map. Here the map 
%\[
% H(\Der'(\fGraphs_{n,c}^\whl)) \to H(\Der(\hoe_n))
%\]
%is defined as follows. Let a cocycle $D\in \Der'(\fGraphs_{n,c}^\whl)$ be given. Then $D\circ f$ defines a cocycle in 
%
%
%It will hence be sufficient to show that the composition 
%\[
% H^0(\GC) \to H^1(\Def(G_\infty\to \Graphs)) \to H^0(\Def(G_\infty))
%\]
%respects Lie brackets. So 
Let homogeneous cocycles $\gamma, \nu \in \fGC_{n, \conn}^\whl$ be given. By the action of $\fGC_{n, \conn}^\whl$ on $\fGraphs_{n,c}^\whl$ we obtain derivations $D_\gamma, D_\nu\in \Der(\fGraphs_{n,c}^\whl)$, using the notation of section \ref{sec:defcomplexes}. They satisfy $\co{\gamma}{\nu}=D_{\co{\gamma}{\nu}}$. The images of $\gamma,\nu$ in $\Def(\hoe_n\to \fGraphs_{n,c}^\whl)_\conn[1]$ have the form $(-1)^{|\gamma|}\bs D_\gamma\circ f$, $(-1)^{|\nu|}\bs D_\nu\circ f$, where $f:\hoe_n\to \fGraphs_{n,c}^\whl$ is the composition $\hoe_n\to e_n \hookrightarrow \fGraphs_{n,c}^\whl$ and $\bs$ shall denote the degree shift operation, which is also responsible for the extra sign. 
Since the map 
\[
 \phi: \Der(\hoe_n) \to \Def(\hoe_n\to \fGraphs_{n,c}^\whl)[1],
\]
is a quasi-isomorphism of complexes, there are derivations $F_\gamma,F_\nu, F_{\co{\gamma}{\nu}}\in \Der(\hoe_n)$ such that 
\begin{align*}
\phi(F_\gamma) &:= (-1)^{|\gamma|}\bs 
f\circ F_\gamma 
= (-1)^{|\gamma|}\bs D_\gamma\circ f + \text{(coboundaries)} \\
\phi(F_\nu) &:= (-1)^{|\nu|}\bs f\circ F_\nu =(-1)^{|\nu|}\bs D_\nu\circ f + \text{(coboundaries)} \\
\phi(F_{\co{\gamma}{\nu}}) &:= (-1)^{|\gamma|+|\nu|}\bs f\circ F_{\co{\gamma}{\nu}} = (-1)^{|\gamma|+|\nu|}\bs D_{\co{\gamma}{\nu}}\circ f + \text{(coboundaries)}.
\end{align*}
 Our goal is to show that
\begin{equation}
\label{equ:Frespectbr}
 \co{F_\gamma}{F_\nu} = F_{\co{\gamma}{\nu}} +\text{(coboundaries)}.
\end{equation}
Since $\phi$ is a quasi-isomorphism, it is sufficient to show the equation 
\[
 \phi(\co{F_\gamma}{F_\nu}) := (-1)^{|\gamma|+|\nu|}\bs f\circ (\co{F_\gamma}{F_\nu}) =  D_{\co{\gamma}{\nu}}\circ f +\text{(coboundaries)}. 
\]
obtained by composing \eqref{equ:Frespectbr} with $\phi$ from the left and using the defining property of $F_{\co{\gamma}{\nu}}$.
The bracket in the middle is the commutator (of $\bbS$-module morphisms). Compute
\begin{align*}
f\circ F_\gamma \circ F_\nu 
&= 
 (D_\gamma \circ f +  \text{(coboundaries)} ) \circ F_\nu 
\\ &=
 D_\gamma \circ f \circ F_\nu +  \text{(coboundaries)} 
\\ &=
 D_\gamma \circ (D_\nu \circ f +  \text{(coboundaries)} )+  \text{(coboundaries)} 
\\ &=
 D_\gamma \circ D_\nu \circ f +  \text{(coboundaries)}
\end{align*}
and similarly for $\gamma$ and $\nu$ interchanged. Hence we obtain
\begin{align*}
 \phi(\co{F_\gamma}{F_\nu}) 
 &=
 D_\gamma \circ D_\nu \circ f -(-1)^{|\mu||\nu|} D_\nu \circ D_\gamma \circ f + \text{(coboundaries)}
 = \co{D_\gamma}{D_\nu}\circ f + \text{(coboundaries)}
 \\&= D_{\co{\gamma}{\nu}}\circ f + \text{(coboundaries)}
\end{align*}
This was to be shown.
\end{proof}

\subsubsection{Explicit form}
Let us work out the explicit form of the map 
\[
\psi : \fGC_{n, \conn}^\whl \to
\Def(\hoe_n\to \fGraphs_{n,c}^\whl)_\conn [1].
\]
First note that the differential on the right hand side splits according to eqns. \eqref{equ:diffsplit1}, \eqref{equ:diffsplit2}.

%by definition has the form
%\[
%d = \delta + \co{\alpha}{\cdot}
%\]
%where $\delta$ is induced by the differential on $\fGraphs_{n,c}^\whl$, and $\alpha$ is the Maurer-Cartan element defined by the map $\hoe_n\to \fGraphs_{n,c}^\whl$. This map in turn has two non-vanishing components, for the product and bracket generators of $e_n$. One may decompose $\alpha$ accordingly into two pieces $\alpha=\alpha_{\wedge}+\alpha_{\co{}{}}$. Then we define
%\begin{align*}
%d_{\wedge} &= \co{\alpha_{\wedge}}{\cdot} & d_{\co{}{}} &= \co{\alpha_{\co{}{}}}{\cdot}
%\end{align*}
%so that 
%\[
%d = \delta + d_{\wedge} + d_{\co{}{}}.
%\]
Now let $\gamma\in \fGC_{n, \conn}^\whl$ be a homogeneous element. In order to find its image in $\Def(\hoe_n\to \fGraphs_{n,c}^\whl)_\conn$, we have to compute the action of $\gamma$ on the two 2-vertex graphs to which the generators of $e_n$ are sent to, see \eqref{equ:enGramapdef}. Since both these graphs contain no internal vertices, the only relevant part of the action is 
\[
a \mapsto \gamma_1 \cdot a
\]
where $\gamma_1\in \fGraphs_{n,c}^\whl(1)$ is obtained from $\gamma$ by marking vertex 1 as external (cf. \eqref{equ:x1def}) and the action $\cdot$ is defined in \eqref{equ:cdotdef}.
Note also that $\gamma_1$ (or any element in $\fGraphs_{n,c}^\whl)_\conn$) may be considered as an element of $\Def(\hoe_n\to \fGraphs_{n,c}^\whl)_\conn$ as well (sending the counit to $\gamma_1$).
Also, the Lie bracket on the convolution dg Lie algebra $\Def(\hoe_n\to \fGraphs_{n,c}^\whl)$ with such elements is given by a formula formally identical to \eqref{equ:cdotdef}. One may hence check that $\gamma$ is sent to 
\begin{equation}
\label{equ:gammaimg}
(-1)^{|\gamma|} \bs \co{\gamma}{\alpha} = -\bs \co{\alpha}{\gamma} = -\bs (d_{\wedge}\gamma_1 + d_{\co{}{}}\gamma_1  )
\end{equation}
where $\bs$ shall denote the degree shift operation, which is also responsible for the extra sign on the left (since $\gamma_1$ appears to the right of $\bs$).

%One can show the following proposition by a direct combinatorial computation.
%\begin{prop}
% Denote the map on cohomology induced by $\Phi$ by
%\[
% \bar \Phi\colon H(\fGC_{n, \conn}^\whl) \to H(\Def(\hoe_n\to \fGraphs_{n,c}^\whl)_\conn)[1] %\cong H(\Def(e_n)_\conn).
%\]
%Then $\bar \Phi$ is injective with one dimensional cokernel, in degree 0.
%\end{prop}
%The proof is technical and will be given in Appendix \ref{sec:defengra}.
%Assuming the proposition, we immediately obtain Proposition \ref{prop:fgcconndefenconn} as a corollary, and hence, together with Proposition \ref{prop:GCred} a proof of Theorem \ref{thm:GCdef}.

%\begin{rem}
%$\Def(\hoe_n\to \fGraphs_{n,c}^\whl)_\conn$ may be equipped with a filtration by the number of clusters of size one and valence one (using the terminology from Appendix \ref{sec:defengra}). By a combinatorial argument, one may then compute the first convergent of the associated spectral sequence and find that it is isomorphic to $\fGC_{n,\conn}[-1]\oplus \GF[-1]$. This can be made into a quick proof of Theorem \ref{thm:GCdef}. We will however follow a different route below. 
%\end{rem}

\subsection{The second description of the map%, and a proof of Theorem \ref{thm:GCdef}
}
\label{sec:seconddef}
There are natural maps $\hoLie_n\to \hoe_n$ and $e_n\to \Graphs_n^\whl\to \Gra_n^\whl$. They allow one to write a sequence of maps
\begin{equation*}
 0\to \Def(\hoe_n\to e_n) \to \Def(\hoe_n\to \Gra_n^\whl) \to \Def(\hoLie_n\to \Gra_n^\whl) \to 0 .
\end{equation*}
All these spaces can be written as (suitably degree shifted) symmetric powers of their connected parts. In this section, we will only care about the connected parts.
\begin{equation}
\label{equ:defseq}
0\to\Def(\hoe_n\to e_n)_\conn  \to \Def(\hoe_n\to \Gra_n^\whl)_\conn \to \Def(\hoLie_n\to \Gra_n^\whl)_\conn%
 \to 0.
\end{equation}
It is shown in Proposition \ref{prop:hoengracohomfull} that 
\[
H(\Def(\hoe_n\to \Gra_n^\whl)_\conn) \cong \GF[-1].
\]
The single non-trivial cohomology class is represented by a graph with one edge and is easily checked to map to zero in $H(\Def(\hoLie_n\to \Gra_n^\whl)_\conn)$.
%The above sequence is exact on the left and right, but not in the middle. The cohomology of the middle term is calculated in Proposition \ref{prop:hoengracohom} to be $\R[-1]$. Furthermore this $\R[-1]$ is mapped to zero in $H^1(\Def(\hoLie_n\to \Gra_n^\whl)\conn)$.
Assume now that the sequence \eqref{equ:defseq} was exact. Then from the corresponding long exact sequence in homology together with the previous remarks one could conclude that
\[
 H^{\bullet+1}(\Def(\hoe_n\to e_n)_\conn) = H^\bullet(\Def(\hoLie_n\to \Gra_n^\whl)_\conn) \oplus \GF[-1].
\]
Hence, using that 
\[
H(\Def(\hoLie_n\to \Gra_n^\whl)_\conn) \cong H(\GC_n) \bigoplus_{ \substack{j\geq 1 \\ j\equiv 2n+1 \mod 4} } \GF[n-j]
\]
by Proposition \ref{prop:GCred}, Theorem \ref{thm:GCdef} would be proven. The problem is now that \eqref{equ:defseq} is \emph{not} exact in the middle. However, one can cure that defect:
\begin{prop}
\label{prop:homseq}
 There is a long exact sequence in homology
\begin{multline}\label{equ:homseq}
  \cdots \to H^k(\Def(\hoe_n\to \Gra_n)_\conn) \to H^k(\Def(\hoLie_n\to \Gra_n)_\conn) \to
\\
%\stackrel{\phi}{
\to H^{k+1}(\Def(\hoe_n\to e_n)_\conn) \to H^{k+1}(\Def(\hoe_n\to \Gra_n)_\conn) \to \cdots
\end{multline}
%On graphs $\Gamma\in \GC_n\subset \Def(\hoLie_n\to \Gra_n)$, the connecting morphism $\phi$ agrees with the map $\Phi$ defined above.\footnote{Here the identification $H(\Def(\hoe_n\to e_n)) \cong H(\Def(\hoe_n\to \Graphs_n))$ is used.}
\end{prop}

For the proof one can apply the following Lemma, which the author learned from D. Kazhdan.
\begin{lemma}
\label{lem:kazhdan}
 Let
\[
 0 \to A \to B \to C \to 0
\]
be a sequence of complexes such that the composition of consecutive arrows is 0. Assume that the total cohomology of the double complex $D=A[-1] \oplus B\oplus C[1]$, vanishes. Then there is a long exact sequence in cohomology
\[
 \cdots \to H^k(B) \to H^k(C) \to H^{k+1}(A) \to H^{k+1}(B) \to \cdots
\]
\end{lemma}
\begin{proof}[Proof (Sketch)]
Compute the cohomology of the double complex $D$ using the associated spectral sequence. The $E^1$-term is 
\[
 0 \to H(A) \to H(B) \to H(C) \to 0.
\]
The spectral sequence must collapse at the $E^3$ page by degree reasons, and hence $E^3=0$ by the assumption in the Lemma, and the fact that spectral sequences of bounded double complexes converge to the true cohomology. Concretely, the nontrivial parts of the complex $E^2$ are
\begin{gather*}
 0\to \ker(H(B) \to H(C)) / \im(H(A) \to H(B)) \to 0 \\
0\to \ker(H(A) \to H(B)) \to \coker(H(B) \to H(C)) \to 0.
\end{gather*}
From the vanishing of the cohomology of $E^3$ the exactness of the long sequence then follows.
\end{proof}

 Let us now turn to the proof of Proposition \ref{prop:homseq}. 
We would like to apply the Lemma to the sequence \eqref{equ:defseq}. However, since the middle two arrows do not compose to zero, we need to replace 
$\Def(\hoe_n\to e_n)_\conn$ in the sequence \eqref{equ:defseq} by a sub-complex. Let $\Xi\subset \Def(\hoe_n\to e_n)$ the subcomplex of Appendix \ref{sec:endef}. It is shown in Proposition \ref{prop:Xiqiso} that the inclusion is a quasi-isomorphism. Define $\Xi_\conn\subset \Xi$ to be the connected part. It is easily checked that $\Xi_\conn\to \Def(\hoe_n\to e_n)_\conn$ is also a quasi-isomorphism, and hence so is $\Xi_\conn\to \Def(\hoe_n\to \Gra_n^\whl)_\conn$.
We will apply the above lemma to the sequence
\begin{equation}
\label{equ:defseq2}
0\to \Xi_\conn  \to \Def(\hoe_n\to \Gra_n^\whl)_\conn \to \Def(\hoLie_n\to \Gra_n^\whl)_\conn%
 \to 0.
\end{equation}
In other words, $A=\Xi_\conn$, $B=\Def(\hoe_n\to \Gra_n^\whl)_\conn$ and $C=\Def(\hoLie_n\to \Gra_n^\whl)_\conn$.
 In order for the lemma to be applicable, we need to check the following result, which is proven in Appendix \ref{sec:prekazhdanproof}.
\begin{lemma}
\label{lem:prekazhdan}
The double complex defined by the sequence \eqref{equ:defseq2} is acyclic.
\end{lemma}
Let us believe Lemma \ref{lem:prekazhdan} for now and use it to prove the first assertion of Theorem \ref{thm:GCdef}. First we apply Lemma \ref{lem:kazhdan}, which proves Proposition \ref{prop:homseq}. Together with the remarks made prior to this proposition, Proposition \ref{prop:fgcconndefenconn} follows. Using the description of the cohomology of $\Def(\hoe_n\to e_n)$ provided by Proposition \ref{prop:defenbyconnectedpart}, together with the result of Proposition \ref{prop:GCred} then shows \eqref{equ:GCdef}.%Theorem \ref{thm:GCdef}.
\hfill \qed
 
\subsection{The two maps agree}
\label{sec:mapsagree}
The goal of this subsection is to show the following result:
\begin{prop}
\label{prop:mapsagree}
The two maps
\[
H(\fGC_{n, \conn}) \to H(\Def(\hoe_n\to e_n)_\conn) [-1]
\]
that were defined in sections \ref{sec:firstdef} and \ref{sec:seconddef} are identical.
\end{prop}

The first map came from the map 
\[
\fGC_{n, \conn} \to \Def(\hoe_n\to \fGraphs_{n,c})_\conn [-1] % TODO: change notation
\]
given by the action of $\fGC_{n, \conn}$ on $\fGraphs_{n, c}$, which is worked out explicitly in Appendix \ref{sec:optwists_graphs}.
Let $\Gamma \in\fGC_{n, \conn}$ be a cocycle. Then its image is the cocyle
\[
X = (d_{\wedge} + d_{\co{}{}}) \Gamma_1
\] 
cf. \eqref{equ:gammaimg}. To obtain an element of $H(\Def(\hoe_n\to e_n)_\conn)$ one uses Lemma \ref{lem:defpresqiso} and the fact that the inclusion $e_n\to \fGraphs_{n,c}$ is a quasi-isomorphism. Concretely, these facts guarantee that there is some $X'\in \Def(\hoe_n\to \fGraphs_{n,c})_\conn$ such that 
\begin{equation}
\label{equ:XplusdX}
X + d X' \in \Def(\hoe_n\to e_n)_\conn \subset \Def(\hoe_n\to \fGraphs_n)_\conn.
\end{equation}
The cocycle $X + d X'$ is a representative in $\Def(\hoe_n\to e_n)_\conn$ of the image of the class represented by $\Gamma$ under the first map. It will be useful to know the following lemma.
\begin{lemma}
\label{lem:Xinftypart}
The element $X'$ above may be chosen so that $\pi(X') =\Gamma$ where 
\[
\pi:  \Def(\hoe_n\to \fGraphs_{n,c})_\conn \to \Def(\hoLie_{n}\to \Gra_n)_\conn
\]
is the natural projection (induced by  $\hoLie_{n}\hookrightarrow\hoe_n$ and $\fGraphs_n\to\Gra_n$).
\end{lemma}

We will postpone the proof and first use the Lemma to show Proposition \ref{prop:mapsagree}.
\begin{proof}[Proof of Proposition \ref{prop:mapsagree}.]
Let us first unravel the definition of the second map, as the connecting homomorphism in the long exact sequence of Proposition \ref{prop:homseq}. 
Given the cocycle $\Gamma\in\fGC_{n, \conn}$ as above we have to pick some element $Y$ such that 
\[
d Y \in \Def(\hoe_n\to e_n)_\conn
\]
and such that $\pi'( Y) = \Gamma$ where $\pi': \Def(\hoe_n\to \Gra_n)_\conn \to \Def(\hoLie_{n-1}\to \Gra_n)_\conn$ is the restriction map. Then the image of the cohomology class represented by $\Gamma$ under the connecting homomorphism is the class represented by $dY$ in $H(\Def(\hoe_n\to e_n)_\conn)$.

Now, due to Lemma \ref{lem:Xinftypart}, we may pick $Y=\pi''(X')$ where 
\[
\pi'':  \Def(\hoe_n\to \fGraphs_{n,c})_\conn \to \Def(\hoe_n\to \Gra_n)_\conn
\]
is the projection induced by $\fGraphs_{n,c}\to\Gra_n$.
Indeed, for this $Y$, $\pi'(Y)=\pi'(\pi''(X'))=\pi(X')=\Gamma$ and 
\[
dY = \pi''(dX')= \pi''(A - X) =\pi''(A) = A 
\]
where $X$ is as above and  
\begin{equation}
\label{equ:A}
A = X+dX' \in  \Def(\hoe_n\to e_n)_\conn.
\end{equation}
Hence the image of (the class represented by) $\Gamma$ under the second map is (the class represented by) $A$, the same as the image under the first map.
\end{proof}

Now we need to prove Lemma \ref{lem:Xinftypart}. To do this, we first show the following preliminary result.
\begin{lemma}
\label{lem:Y}
Let $\Gamma$ and $X$ be as above. Then there is an element $X''\in \Def(\hoLie_{n}\to \fGraphs_n)_\conn$ such that 
\begin{enumerate}
\item $dX''+ p(X)=0$ where 
\[
p \colon \Def(\hoe_n\to \fGraphs_{n,c})_\conn \to \Def(\hoLie_{n}\to \fGraphs_n)_\conn
\]
is the restriction.
\item $p'(X'')=\Gamma$ where 
\[
p' \colon \Def(\hoLie_{n}\to \fGraphs_{n,c})_\conn \to \Def(\hoLie_{n}\to \Gra_n)_\conn
\]
is the projection. 
\end{enumerate} 
\end{lemma} 
\begin{proof}[Proof sketch.]
%Let $\op P$ be an operad without nullary operations receiving a map from $\hoLie_{n-1}$. As noted in \ref{} $\Tw\op P(0)\cong \Def(\hoLie_{n-1}\to \op P)$.

Unravelling the definitions, one sees that 
\begin{equation}
\label{equ:otherrecall}
\Def(\hoLie_{n}\to \fGraphs_n) \cong \prod_{N\geq 0}\prod_{ M\geq 1} \Gra_n(N+M)^{S_N\times S_M}[-(N+M-1)n]
\end{equation}
and 
\begin{equation}
\label{equ:fGCrecall}
\fGC_n \cong \prod_{N'\geq 1}\Gra_n(N')^{S_{N'}}[-(N'-1)n].
\end{equation}
In the graphical language, $N$ above corresponds to the number of internal vertices while $M$ corresponds to the number of external vertices in graphs.
We will define a map of graded vector spaces
\[
f\colon \fGC_n \to \Def(\hoLie_{n}\to \fGraphs_{n,c}).
\]
To do that, we have to specify the composition with the projections $\pi_{N,M}$ to each factor in the direct product.
We will set
\[
\pi_{N,M} \circ f
=
\begin{cases}
\iota_{N,M} \circ \pi_{N+M} & \text{for $M\geq 2$} \\
0 & \text{otherwise}
\end{cases}
\]
where $\pi_{N'}$ is the projection to the $N'$-th factor in \eqref{equ:fGCrecall}, and 
\[
\iota_{N,M} \colon \Gra_n(N+M)^{S_{N+M}}[-(N+M-1)n] \hookrightarrow \Gra_n(N+M)^{S_N\times S_M}[-(N+M-1)n]
\]
is the inclusion. It is easy to see that $f$ restricts to a map of the connected parts, which we denote:
\[
f_\conn\colon \fGC_{n, \conn} \to \Def(\hoLie_{n}\to \fGraphs_{n,c})_\conn.
\]
Now we define 
\[
X'' := f_\conn(\Gamma).
\]
The projection $p'$ in the statement of the Lemma projects onto the $N=0$ factors in \eqref{equ:otherrecall}. Hence it is clear that $p'(X'')=\Gamma$.

The other assertion of the lemma is more difficult to verify. Assume (without loss of generality) that $\Gamma$ is a linear combination of graphs with exactly $N'$ vertices. Then
\[
X'' = f( \Gamma ) = \Gamma_2 + \Gamma_3 +\cdots + \Gamma_{N'}
\]
where $\Gamma_j := \iota_{N'-j, j}(\Gamma)$, silently considered as elements of \eqref{equ:otherrecall}.
The differential on $\Def(\hoLie_{n-1}\to \fGraphs_n)$ has two parts, $d=\delta+d_{\co{}{}}$, where $\delta$ is the part coming from the differential on $\fGraphs_n$ and $d_{\co{}{}}$ is the remainder. We have
\begin{equation}
\label{equ:pX}
p(X)= d_{\co{}{}} \Gamma_1.
\end{equation}
We leave for the reader to check that for all $N, M$
\[
d_{\co{}{}} \circ \iota_{N, M-1} + \delta \circ \iota_{N-1, M} = \iota_{N, M}\circ \delta
\]
where we (abusively) consider the $\iota_{i,j}$ as taking values in $\Def(\hoLie_{n-1}\to \fGraphs_{n,c})$.\footnote{It is not difficult, but lengthy to write down.}
Applying this equality to $\Gamma$ and using that $\Gamma$ is a cocycle, i.e., $\delta \Gamma=0$, we obtain
\[
d_{\co{}{}}\Gamma_j + \delta \Gamma_{j+1} = 0.
\]
Together with \eqref{equ:pX}, this shows the first assertion of the Lemma.
\end{proof}

\begin{proof}[Proof of Lemma \ref{lem:Xinftypart}.]
First we note that
\[
H(\Def(\hoLie_{n}\to \fGraphs_{n,c})_\conn)
\cong 
H(\Def(\hoLie_{n}\to e_n)_\conn)
\cong
0.
\]
Here the first equality is due to Lemma \ref{lem:defpresqiso} and Propositions \ref{prop:Graphsincl} and \ref{prop:Graphscohom} and the second equality is easily checked since $\Def(\hoLie_{n}\to e_n)_\conn$ is a 2-dimensional complex.

Now let $\Gamma$ and $X$ be as above and pick any $X'$ satisfying \eqref{equ:XplusdX} (but possibly $\pi(X')\neq \Gamma$). We will show Lemma \ref{lem:Xinftypart} by constructing an element $Z\in \Def(\hoe_n\to \fGraphs_n)_\conn$ such that $\pi(X'+dZ)= \Gamma$. (Then replacing $X'$ by $X'+dZ$ the assertions of the Lemma are satisfied.) 

Let $X''$, $p$ and $p'$ be as in Lemma \ref{lem:Y}. The element $p(X')-X''$ is closed since 
\[
dp(X')-dX''= p(dX')-dX'' = p(A-X)+dX'' = -p(X)-dX'' +p(A) = p(A)
\]
where $A$ is as in \eqref{equ:A}. But $p(A)$ lives in the 2-dimensional complex $\Def(\hoLie_{n}\to e_n)_\conn$, and it is easily verified that $p(A)\neq 0$ may only happen if $\Gamma$ contains the graph with one vertex, or the graph with two vertices and one edge. Since the latter graph is the coboundary of the former, we may assume that neither is contained in $\Gamma$ and hence $p(A)=0$.
Hence there is some $Y\in \Def(\hoLie_{n}\to \fGraphs_{n,c})_\conn$ such that $X''-p(X')=dY$.
Now pick any $Z\in \Def(\hoe_n\to \fGraphs_{n,c})_\conn$ such that $p(Z)=Y$.
Let us compute 
\[
\pi(X'+dZ) = p'\circ p (X'+dZ) = p'(p(X') + dY) = p'(X'')= \Gamma
\]
where we used Lemma \ref{lem:Y} twice. This shows Lemma \ref{lem:Xinftypart}.
%Since $\Def(\hoLie_{n-1}\to \fGraphs_n)_\conn$ is a sub-graded vector space of $\Def(\hoe_n\to \fGraphs_n)_\conn$, we may consider 
\end{proof}

\subsection{Proof of Theorem \ref{thm:GCdef}}
In section \ref{sec:seconddef} we established equation \eqref{equ:GCdef}, i.e., the first part of Theorem \ref{thm:GCdef}. It remains to show the second assertion, i. e., that $H(\GC_n)\to  H(\Der(\hoe_n))$ is a map of Lie algebras. We will in fact show the stronger assertion that 
$H(\fGC_{n,\conn}^\whl)\to  H(\Der(\hoe_n))$ is a map of Lie algebras. By Proposition \ref{prop:mapsagree} we may as well use the first description of the map $H(\fGC_{n,\conn}^\whl)\to  H(\Der(\hoe_n))$, introduced in section \ref{sec:firstdef}. But for this version the statement is the content of Proposition \ref{prop:firstmaplie}. Hence Theorem \ref{thm:GCdef} has been shown.
\hfill\qed

\begin{rem}
The proof of Theorem \ref{thm:GCdef} we have given here is not the shortest possible. One may shorten it by using only the first description of the map (from section \ref{sec:firstdef}). However, the second description of the map will still be needed below for some other results.
\end{rem} 
\section{The Grothendieck-Teichm\"uller Lie algebra}
\label{sec:grt}
In this section we collect some facts about the Grothendieck-Teichm\"uller Lie algebra.

{\bf Notation:}
Here and in the following we adopt the convention that, if we omit the subscript $n$ in $\GC_n$, $\Graphs_n$ etc., $n=2$ is implied. In particular $\GC:=\GC_2$.

% This section contains several equivalent definitions of the Grothendieck-Teichm\"uller Lie algebra $\grt$. More precisely, we will define several spaces whose cohomology is isomorphic 
\subsection{The standard definition}
\label{ssec:grtstd}
Let $\F_2=\K\langle\langle X,Y \rangle\rangle$ be the completed free associative algebra in generators $X,Y$. There is a coproduct $\Delta$ on $\F_2$ by declaring $X,Y$ to be primitive, i.e., $\Delta X=X\otimes 1 +1 \otimes X$, $\Delta Y=Y\otimes 1 +1 \otimes Y$. It is part of a topological Hopf algebra structure on $\F_2$. We call an element $\Phi\in \F_2$ \emph{group-like} if $\Delta \Phi = \Phi {\otimes} \Phi$. Equivalently, $\phi$ is group-like if $\Phi=\exp(\phi)$ with $\phi\in \hat{\F}_\Lie(X,Y) \subset \F_2$ being an element of the completed free Lie algebra generated by $X$ and $Y$.

Let $\alg{t}_n$ ($n=1,2,3,\dots$) be the Drinfeld Kohno Lie algebra. It is generated by symbols $t_{ij}=t_{ji}$, $1\leq i,j\leq n$, $i\neq j$, with relations $\co{ t_{ij} }{ t_{ik}+t_{kj} }=0$ for $\#\{i,j,k\}=3$ and $\co{ t_{ij} }{ t_{kl} }=0$ for $\#\{i,j,k,l\}=4$. 
Consider the following set of equations for group-like elements $\Phi\in \F_2$, depending on some yet unspecified parameter $\mu\in \K$.
\begin{align}
 \label{equ:pent}
\Phi(t_{12}, t_{23}+t_{24})\Phi(t_{13}+ t_{23},t_{34}) &
= \Phi(t_{23}, t_{34})\Phi(t_{12}+t_{13}, t_{24}+t_{34})\Phi(t_{12}, t_{23}) \\
 \label{equ:hex}
e^{\mu(t_{13}+t_{23})/2} &= \Phi(t_{13}, t_{12})e^{\mu t_{13}/2}\Phi(t_{13}, t_{23})^{-1} e^{\mu t_{23}/2} \Phi(t_{12}, t_{23})\\
 \label{equ:antisymm}
\Phi(x, y)&=\Phi(y,x)^{-1}
\end{align}

\begin{defi}
 The group-like solutions $\Phi\in \F_2$ of \eqref{equ:pent}, \eqref{equ:hex}, \eqref{equ:antisymm} 
are called \emph{Drinfeld associators} for $\mu\neq 0$ and elements of the \emph{Grothendieck-Teichm\"uller group} $\GRT_1$ for $\mu= 0$.
\end{defi}

H. Furusho has shown the following remarkable Theorem.
\begin{thm}[Furusho \cite{furusho}]
 Any group-like solution $\Phi\in \F_2$ of \eqref{equ:pent} automatically satisfies \eqref{equ:hex} and \eqref{equ:antisymm} for $\mu=\pm \sqrt{24c_2(\Phi)}$, where $c_2(\Phi)$ is the coeffcient of $XY$ in $\Phi$.
\end{thm}
We reprove the $\mu=0$ case of this result in Appendix \ref{sec:furusho}.

\begin{rem}
 It also follows from \eqref{equ:pent} that $\Phi$ contains no terms linear in $X,Y$, i.e., 
\[
 \Phi(X,Y) = 1+ c_2(\Phi) (XY-YX) + (\text{higher orders}) .
\]
\end{rem}

The group structure on the Grothendieck-Teichm\"uller group $\GRT_1$ is as follows.
\[
 (\Phi_1\cdot \Phi_2)(X,Y) = \Phi_1(X,Y) \Phi_2(X,\Phi_1(X,Y)^{-1}Y\Phi_1(X,Y)).
\]
We will actually be mostly interested in its Lie algebra, the Grothendieck-Teichm\"uller Lie algebra $\grt_1$. It is given by a Lie series $\phi\in \hat{\F}_\Lie(X,Y)$ such that the following hold.
\begin{align}
 \label{equ:pentl}
\phi(t_{12}, t_{23}+t_{24}) + \phi(t_{13}+ t_{23},t_{34}) &=
 \phi(t_{23}, t_{34}) + \phi(t_{12}+t_{13}, t_{24}+t_{34}) + \phi(t_{12}, t_{23}) \\
 \label{equ:hexl}
 \phi(X, Y)+\phi(Y, -X-Y)+\phi(-X-Y, X) &= 0\\
 \label{equ:antisymml}
\phi(X, Y)+\phi(Y,X) &=0
\end{align}
Again, by Furusho's result, it suffices to require \eqref{equ:pentl} and that $\phi$ contains no quadratic term (i.e., no $\co{X}{Y}$). 
%TODO: add description of Lie bracket
One may embed $\grt_1$ into $\alg t_3$ as a sub-vector space by setting $X=t_{12}$ and $Y=t_{23}$.

It is not hard to see that the Lie algebra $\grt_1$ carries a complete grading, inherited from the complete grading on $\hat{\F}_\Lie(X,Y)$ obtained by assigning the generators $X$ and $Y$ degree $+1$. One may define the Lie algebra
\[
 \grt := \K\ltimes \grt_1
\]
where elements $\lambda$ of the abelian Lie algebra $\K$ act on homogeneous elements $x\in\grt_1$ of degree $r$ by multiplictaion by $r$, i. e., $\co{\lambda}{x}=r x$.

\subsection{Definition as ``Harrison'' cohomology of \texorpdfstring{$\alg{t}$}{t}}
\label{ssec:cohomt}
Since the defining relations of $\alg{t}_n$ are homogeneous with respect to the number of $t_{ij}$'s occuring, the Lie algebras $\alg{t}_n$ are naturally graded by assigning degree one to the generators $t_{ij}$. Let us call the resulting degree of homogeneous elements the ``$t$-degree''. Distinguish the $t$-grading from the cohomological grading, with respect to which we consider $\alg t_n$ to be concentrated in degree $0$.

Let $\hat{\alg{t}}_n$ be the degree completion of $\alg{t}_n$ with respect to the $t$-degree.
The spaces $\hat{\alg{t}}_n$ in fact form an operad $\hat{\alg{t}}$ of Lie algebras. % (see \cite{}). 
Hence the Chevalley complexes $C(\hat{\alg{t}})$ form an operad of commutative coalgebras, and in particular an operad of vector spaces.
There is a map $\Com \to C(\hat{\alg{t}})$ of operads of vector spaces, sending the generator of $\Com$ to $1\in C(\hat{ \alg{t}}_2)$. Hence one also obtains a map $\Com_\infty \to C(\hat{\alg{t}})$, and can consider the deformation complex 
\[
\Def(\Com_\infty \to C(\hat{\alg{t}})) \cong \prod_{N\geq 1} \Hom_{S_N} (\Com^\vee(N), S(\hat{\alg{t}}_N[1]))
\]
where $S(\cdots)$ denotes the symmetric product space.
 The differential on $\Def(\Com_\infty \to C(\hat{\alg{t}}))$ can be written as $d_H+d_{CE}$, where $d_{CE}$ comes from the differential on $C(\hat{\alg{t}})$ (the Chevalley-Eilenberg differential). The remainder of the differential we denote by $d_H$. This part formally resembles the Harrison differential.
The deformation complex above contains a subcomplex which we denote (abusing notation) by
%By abuse of notation, we denote the degree 1 subspace of the associated graded with respect to the filtration above
\[
 \Def(\Com_\infty \to \hat{\alg{t}}[1]) 
 :=
  \prod_{N\geq 2} \Hom_{S_N} (\Com^\vee(N), \hat{\alg{t}}_N[1])
\]
as in section \ref{sec:closerlook}. This complex is endowed with the differential $d_H$.
Up to differences in grading and completion, the complex $\Def(\Com_\infty \to \hat{\alg{t}}[1])$ has been defined by V. Drinfeld, see \cite{drinfeld}, above Proposition 5.9.
\begin{lemma}[Drinfeld \cite{drinfeld}]
\label{lem:tcohom}
 $H^1(\Def(\Com_\infty \to \hat{\alg{t}}[1] ))\cong \grt_1$.
\end{lemma}
\begin{proof}
By degree reasons, $H^1(\Def(\Com_\infty \to \hat{\alg{t}}[1] ))$ is the space of closed elements in 
\[
 \Hom_{S_3}(\Com^\vee(3), \hat{\alg{t}}_3[1])
\]
modulo exact elements. From the fact that $\alg{t}_2$ is one-dimensional, one can show that there are no exact elements. Hence one needs to compute closed elements in $\hat{\alg{t}}_3$ with the correct symmetry properties. However, one checks that the required symmetries are exactly equations \eqref{equ:hexl} and \eqref{equ:antisymml}, and closedness is exactly \eqref{equ:pentl}.
\end{proof}
The cohomology of $\Def(\Com_\infty \to \hat{\alg{t}}[1] ))$ in degrees smaller than one is trivial to determine. The cohomology in degrees $>1$ on the other hand is unknown. In particular, it is a hard and long standing conjecture that the degree 2 cohomology vanishes.
\begin{lemma}
\label{lem:tharrtriv}
$H^{<0}(\Def(\Com_\infty \to \hat{\alg{t}}[1] ))=0$ and $H^{0}(\Def(\Com_\infty \to \hat{\alg{t}}[1] ))\cong \GF$.
\end{lemma}
\begin{proof}
The complex $\Def(\Com_\infty \to \hat{\alg{t}}[1] )$ has no components in degrees $<0$, and is one dimensional in degree 0, spanned by a closed element.  
\end{proof}

In \cite{pavol} P. \v Severa and the author introduced $L_\infty$ algebras $\CG(n)$, forming an operad of $L_\infty$-algebras, such that $H(\CG)=\hat{\alg{t}}$. See also section \ref{sec:graphsn} and in particular Proposition \ref{prop:ICG} of this paper. We will use here the (quasi-isomorphic) version with tadpoles $\CG^\whl$ for technical reasons.\footnote{One could proceed without the tadpoles, but then some results would attain a less elegant form.}
There is again a map $\Com \to C(\CG^\whl)$ by sending the generator to $1\in C(\CG^\whl(2))$. We consider the deformation complex 
\[
\Def(\Com_\infty \to C(\CG^\whl) )\cong \Def(\Com_\infty \to \Graphs_2^\whl ).
\]
Again, this complex contains a subcomplex which we denote (abusing notation)
\begin{equation}
\label{equ:defcomicg}
\Def(\Com_\infty \to \CG^\whl[1] )
:= 
  \prod_{N\geq 1} \Hom_{S_N} (\Com^\vee(N), \CG^\whl(N)[1]).
\end{equation}

%Using a filtration argument one can derive from the previous Proposition the following.
\begin{prop}
\label{prop}
$H(\Def(\Com_\infty \to \CG^\whl[1] ))\cong H(\Def(\Com_\infty \to \hat{\alg{t}}[1] ))$. In particular $H^1(\Def(\Com_\infty \to \CG^\whl[1] ))\cong  \grt$.
\end{prop}
\begin{proof}
Let us first recall some generalities on filtrations and spectral sequences. Suppose $A$ is some complex with a complete descending filtration. Then the associated spectral sequence may or may not converge to the cohomology $H(A)$. However, if the associated graded of $A$ is acyclic, then it does converge, i. e., $H(A)=0$. Now let $f: A\to B$ be a map of complexes, both equipped with complete descending filtrations (respected $f$). Then if the associated graded map $\gr f$ is a quasi-isomorphism, so is $f$. To see this, simply apply the previous statement to the mapping cone. 
Now let us turn to the case at hand. We recall from \cite{pavol} that there is a sub-operad of $L_\infty$ algebras $\TCG^\whl\subset \CG^\whl$ obtained by truncation, that comes equipped with quasi-isomorphisms
\begin{equation}
\label{equ:TCG}
\hat{\alg t} \twoheadleftarrow \TCG^\whl \hookrightarrow \CG^\whl.
\end{equation}
Hence we obtain maps
\begin{equation}
\label{equ:defTCG}
\Def(\Com_\infty \to \hat{\alg t}[1] )
\leftarrow
\Def(\Com_\infty \to \TCG^\whl[1] )
\to
\Def(\Com_\infty \to \CG^\whl[1] )
\end{equation}
that we claim are quasi-isomorphisms. 
Here the middle complex is defined analogously to the right hand one, by replacing all occurrences of $\CG^\whl$ by $\TCG^\whl$. All three complexes come equipped with a complete descending filtration
\[
\mF^p \Def(\Com_\infty \to \alg g[1] )
=
\prod_{N\geq p} \Hom_{S_N} (\Com^\vee(N), \alg{g}(N)[1])
\]
where $\alg g$ is either $\hat{\alg t}$, $\TCG^\whl$ or $\CG^\whl$. 
The graded versions of the morphisms in \eqref{equ:defTCG} are all quasi-isomorphisms since the morphisms in \eqref{equ:TCG} are.
By the generalities recalled at the beginning of this proof, the morphisms in \eqref{equ:defTCG} are hence quasi-isomorphisms. The Lemma hence follows from Lemma \ref{lem:tcohom}.
\end{proof}

An element of $\Def(\Com_\infty \to \CG^\whl[1] )$ is given by a collection of maps in
\[
 \Hom_{S_n}(\Com^\vee(n), \CG^\whl(n)[1]).
\]
Such maps can be depicted as internally connected graphs in $\Graphs_2(n)$, satisfying certain symmetry conditions under interchange of their external vertices. More precisely, they are required to vanish on shuffles: For any $k,l\geq 1$, $k+l=n$, one requires 
\[
 \sum_{\sigma\in ush(k,l)} \sgn(\sigma) \sigma \cdot \Gamma = 0.
\] 
Here the sum is over $(k,l)$-unshuffle permutations, and a permutation acts on graphs $\Gamma$ by interchange of the labels on the external vertices.
One can define a subcomplex $C\subset \Def(\Com_\infty \to \CG^\whl[1] )$ consisting of graphs with one external vertex and exactly one edge connecting to them. Some example graphs in $C$ are shown in Figure \ref{fig:Cgraphs}.
\begin{figure}
 \centering
\[
\begin{tikzpicture}[scale=1,
vert/.style={draw,outer sep=0,inner sep=0,minimum size=5,shape=circle,fill},
helper/.style={outer sep=0,inner sep=0,minimum size=5,shape=coordinate},
default_edge/.style={draw},
ext/.style={draw,outer sep=0,inner sep=0,minimum size=5,shape=circle},
every loop/.style={min distance=10mm,looseness=30}]

\node (v0) at (2.5,6) [ext] {};
\node (v1) at (5.5,6) [ext] {};
\node (v2) at (5,7) [vert] {};
\node (v3) at (5,8) [vert] {};
\node (v5) at (6,8) [vert] {};
\node (v7) at (6,7) [vert] {};
\node (v9) at (2.5,7) [vert] {};

\draw[default_edge] (v2) to (v3);
\draw[default_edge] (v3) to (v5);
\draw[default_edge] (v5) to (v7);
\draw[default_edge] (v2) to (v7);
\draw[default_edge] (v2) to (v5);
\draw[default_edge] (v2) to (v1);
\draw[default_edge] (v3) to (v7);
\draw[default_edge] (v9) to (v0);
\draw[default_edge,in=135, out=45,loop] (v9) to ();
\end{tikzpicture}
\]
\caption{\label{fig:Cgraphs} Some graphs in the complex $C\subset \Def(\Com_\infty \to \CG^\whl[1])$.}
\end{figure}

\begin{prop}
\label{prop:onevertgrt}
 The inclusion $C\hookrightarrow \Def(\Com_\infty \to \CG^\whl[1] )$ is a quasi-isomorphism. Hence $H^1(C)\cong \grt_1$.
\end{prop}
\begin{proof}
 It follows along the lines of the more general proof of Proposition \ref{prop:fCICG}. %There one also sees what has to be changed to circumvent using tadpoles. (To $C$ one would have to add the graph with two external vertices and one edge.)
\end{proof}

%Furthermore, note that (see \cite{pavol}) $\Graphs^\whl(n) = C(\CG^\whl(n))$ is the (completed) bar construction. Hence $\Def(\Com_\infty \to \CG^\whl[1] ) \subset \Hom(\Com_\infty \to \Graphs^\whl)$.
The same proof also shows another statement.
\begin{prop}
\label{prop:comtographs}
 The inclusion $C\hookrightarrow \Def(\Com_\infty \to \Graphs^\whl_2 )$ is a quasi-isomorphism. Hence 
\[
H^j(\Def(\Com_\infty \to \Graphs^\whl_2 ))\cong 
 \begin{cases}
 \grt_1 & \text{j=1} \\
 \GF & \text{j=0} \\
 0 & j<0
 \end{cases}.
\]
\end{prop}
\begin{proof}
 It is a copy of the proof of the more general Proposition \ref{prop:fCqiso2}.
\end{proof}

\subsection{Tamarkin's \texorpdfstring{$\grt_1$}{grt}-action (up to homotopy) on \texorpdfstring{$\hoe_2$}{hoe2} and on \texorpdfstring{$\Tpoly$}{Tpoly} }
\label{sec:tamactions}
\subsubsection{The action on \texorpdfstring{$\hoe_2$}{hoe2}}
\label{ssec:tamactionsginfty}
D. Tamarkin \cite{tamarkingrtaction} showed that there is an action (up to homotopy) of $\grt_1$ on the homotopy Gerstenhaber operad $\hoe_2$ by derivations. 
Let us recall the construction of this action.

Let $\PaP(N)$ be the category with objects the parenthesized permutations of symbols $1,\dots, n$ (for example $(13)(5(24))$ is an object of $\PaP(5)$) and with exactly one morphism between any pair of objects. Let $\PCD(N):= U\alg{t}_N\times \PaP(N)$ be the product of categories, where the completed universal enveloping algebra $U\alg{t}_N$ of $\alg{t}_N$ is considered as a category with one object.\footnote{$\PCD$ stands for ``parenthesized chord diagrams'', see \cite{barnatan}.} The $\PCD(N)$ assemble to form an operad of categories enriched over (complete) Hopf algebras, called $\PCD$. The Grothendieck-Teichm\"uller group $\GRT$ acts faithfully on $\PCD$. In particular, one has an action of the Lie algebra $\grt_1$ on $\PCD$. See \cite{barnatan} for details.
Furthermore, there is a chain of quasi-isomorphisms of operads
\[
 C\Ne\PCD \rightarrow BU\alg{t} \leftarrow e_2 \leftarrow \hoe_2.
\]
Here ``$\Ne$'' denotes the nerve and ``$C$'' the chains functor and $BU\alg{t}$ is the bar construction of $U\alg{t}$. The left arrow is obtained from the map $\PCD\to U\alg{t}$ sending $\PaP$ to the category with one object and one morphism. Since $\hoe_2$ is cofibrant, the action of $\grt_1$ on $C\Ne\PCD$ can be transferred to an action (up to homotopy) on $\hoe_2$. One hence obtains an $L_\infty$ morphism 
\[
 \grt_1 \to \Der(\hoe_2).
\]
D. Tamarkin proved the following theorem.
\begin{thm}[Tamarkin \cite{tamarkingrtaction}]
\label{thm:taminj}
 The above morphism $\grt_1 \to \Der(\hoe_2)$ is homotopy injective, i. e., the induced map $\grt_1 \to H^0(\Der(\hoe_2))$ is injective.
\end{thm}
In fact, D. Tamarkin proved this theorem by showing the following lemma.
\begin{lemma}
\label{lem:tamid}
The composition 
\[
\grt_1 \to H^0(\Der(\hoe_2)) \stackrel{\sim}{\to}
 H^1(\Def(\hoe_2\to C(\alg{t}) ))
 \to 
 H^1(\Def(\Com_\infty\to C(\alg{t}) )) \cong \grt_1
\]
is the identity. Here the first arrow is induced by the $\grt_1$ action on $\hoe_2$ and the third arrow is the restriction.
%composition from the left with the map $\hoe_2\to e_2\to \Graphs$ and the second arrow is the restriction to $\Com_\infty\subset G_\infty$.
\end{lemma} 
A proof of the Lemma and hence of Theorem \ref{thm:taminj} is sketched in Appendix \ref{sec:tamproof}.
%\begin{proof}

%\end{proof}

% \begin{rem}
%  By degree reasons, one can check that the action is trivial on the subspace $\Der(L_\infty\to G_\infty)\subset \Der(G_\infty)$.
% \end{rem}

\subsubsection{The action on \texorpdfstring{$\Tpoly$}{Tpoly}}
\label{sec:tamacttpoly}
Let next $\Tpoly=\Gamma(\R^d;\wedge T\R^d)$ be the space of polyvector fields on $\R^d$. It is naturally a Gerstenhaber algebra, hence we have maps $\hoe_2 \to e_2 \to \End(\Tpoly)$, where $\End(\Tpoly)$ is the endomorphism operad. From the $\grt_1$-action on the left one hence obtains a map
\[
 A\colon \grt_1 \to \Def(\hoe_2 \to \End(\Tpoly))[1].
\]
By rigidity of $\Tpoly$, the image of this map will be exact, cf. Proposition \ref{prop:hoengracohomfull}. So, for any $\phi\in \grt_1$ there will be some $h\in \Der(\hoe_2 \to \End(\Tpoly))$ such that $A(\phi)=d h$.
\begin{rem}
The homotopy $h$ encodes an infinitesimal $\hoe_2$ map between the two $\hoe_2$ structures on $\Tpoly$ related by $\phi$. 
\end{rem}
% TODO: explain the following better
By degree reasons, $\phi$ acts trivially on the $\hoLie_2$ part of the $\hoe_2$-structure on $\Tpoly$, i.e., it is mapped to zero under the restriction map
\[
 \Def(\hoe_2 \to \End(\Tpoly))\to \Def(\hoLie_2 \to \End(\Tpoly)).
\]
We can hence conclude that $dh'=0$ where $h'\in \Def(\hoLie_2 \to \End(\Tpoly))$ is the image of $h$ under the above map. 
In other words, $h'$ encodes a $\hoLie_2$ derivation of $\Tpoly$, or equivalently an $L_\infty$ derivation of $\Tpoly[1]$. 
%\begin{rem}
One can check that $h'$ is determined uniquely up to homotopy by $\phi$, again using rigidity of $\Tpoly$.
%\end{rem}
This action (up to homotopy) of $\grt_1$ on $\Tpoly$ by $\hoLie_2$-derivations was discovered by D. Tamarkin to the knowledge of the author. Degree 0 graph cocycles in $\GC$ can be seen as universal (independent of $d$) $\hoLie_2$ derivations of $\Tpoly$. So the map $\grt_1\to H^0(\GC)$ central to this paper is a universal version of D. Tamarkin's action of $\grt_1$ on $\Tpoly[1]$.
The universal version of Tamarkin's construction is described in more detail in Appendix \ref{sec:tammap}.

\subsection{The map from the graph cohomology to \texorpdfstring{$\grt_1$}{grt}}
By the action of the graph complex $\GC$ on the operads $\Graphs$ and $\Graphs^\whl$ described in Appendix \ref{sec:optwists} we obtain a map of complexes
\begin{equation}
\label{equ:GCtodefcom}
\GC \to \Def(\Com_\infty\to \Graphs^\whl)[1].
\end{equation}
In particular, one obtains a map 
\[
H^0(\GC) \to H^1(\Def(\Com_\infty\to \Graphs^\whl))\cong \grt_1
\]
of vector spaces between the zeroth graph cohomology and $\grt_1$. Our goal in this section and the next is to show that this map is an isomorphism of Lie algebras.

\subsubsection{Factoring through $C$}
By Proposition \ref{prop:comtographs} we know that the map of complexes \eqref{equ:GCtodefcom} must factor through the subcomplex $C$, up to homotopy. In fact, the corresponding map of complexes $\GC\to C$ has a combinatorially very simple form.

%Let $\Gamma\in \GC$ be a graph cochain (i. e., a formal linear combination of graphs). We denote by 
%Let $C$ again be the subcomplex of $\Def(G_\infty\to \Graphs)$ introduced in section \ref{ssec:cohomt}. 
We define the map of complexes
\begin{gather*}
F\colon \GC \to C[1] \\
\Gamma \mapsto 
\begin{tikzpicture}[baseline=-.65ex]
\node [ext, label=-90:1] (v) at (0,0) {};
\node [ext, label=-90:2] (w) at (.5,0) {};
\draw (v) edge (w);
\end{tikzpicture}
\circ_2
\Gamma.
\end{gather*}
Here the right hand side deserves some explanation. Both the two-vertex graph depicted and $\Gamma$ can be seen as elements of the operad $\fGraphs$, and the composition $\circ_2$ is the composition in $\fGraphs$. This yields, a priori, an element in $\fGraphs(1)$. However, it is easy to check that the element in fact lies in the subcomplex 
\[
C[1] \subset \Graphs(1) \subset \fGraphs(1).
\]
It is clear that the map $F$ respects the differentials since the operadic composition does.

\begin{rem}
By writing out the definitions explicitly, one can check that $F$ measures the failure of the map 
\begin{gather*}
\GC \to \Graphs(1) \\
\gamma \mapsto \gamma_1
\end{gather*}
(cf. \eqref{equ:x1def}) to commute with the differentials. I.e., for any graph cochain $\nu\in \GC$ one has the formula
\begin{equation}
\label{equ:Faltdef}
F(\nu) = \delta \nu_1 -(\delta \nu)_1.
\end{equation}
This could also be used as an alternative definition of $F$.
\end{rem}

\begin{prop}
Denote the maps on cohomology induced by $F$ and $G$ again by $F$, $G$. Then the following diagram commutes.
\[
 \begin{tikzpicture}
\matrix(m)[diagram]{
 H(\GC) & &  H(\Def(\Com_\infty\to \Graphs))[1] \\
& H(C)[1] & 
\\};
\draw[->]
  (m-1-1) edge node[auto]{$\scriptstyle{G}$} (m-1-3)
          edge node[auto]{$\scriptstyle{F}$} (m-2-2)
  (m-2-2) edge node[auto]{$\scriptstyle{\cong}$} (m-1-3);
\end{tikzpicture}
\]
%Furthermore, $F$ and $G$ are injective on cohomology.
\end{prop}
\begin{proof}
Let $\gamma\in \GC$ be a graph cocycle. 
Its image in $\Def(G_\infty\to \Graphs)$ is obtained by acting with $\gamma$ on the (Maurer-Cartan) element 
\[
\mu := 
(
\begin{tikzpicture}[baseline=-.65ex]
\node [ext, label=-90:1] (v) at (0,0) {};
\node [ext, label=-90:2] (w) at (.5,0) {};
\end{tikzpicture}
)
\in \Def(\Com_\infty\to \Graphs)
\]%, depicted in Figure \ref{fig:defmc}. 
Hence one needs to follow the construction of the action as in Appendix \ref{sec:optwists_graphs}. Note that the formula slightly simplifies since $\mu$ does not contain internal vertices. %Hence there are no terms (i). 
From $\gamma$ one obtains the element $\gamma_1\in \Graphs(1)$ by marking the first vertex as external. We can consider $\gamma_1$ also as an element of $\Def(G_\infty\to \Graphs)$. The action of $\gamma$ on $\mu$ is the same as the Lie bracket (in $\Def(\Com_\infty\to \Graphs)$) of $\gamma_1$ and $\mu$, $\co{\mu}{\gamma_1}$. %More explicitly, the latter term splits into two terms corresponding to the two graphs in $\mu$, i.e., $\co{\mu}{\gamma_1} = d_L\gamma_1 +d_H \gamma_1$, where we use the notation of Appendix \ref{sec:derengradiff}. Concretely, the term $d_L\gamma_1$ consists of graph with two clusters and $d_H \gamma_1$ of graphs with one cluster, see Figure \ref{fig:gamma1ex} for an example. We now have the desired element of $\Def(G_\infty\to \Graphs)_\conn$ and want to map it to $\Def(\Com_\infty\to \Graphs)_\conn$. This is easy, one just picks out those graphs with only one cluster. In our case, we are left with $d_H \gamma_1$. 
Note that the differential on $\Def(\Com_\infty\to \Graphs)$ has the form $\delta+\co{\mu}{\cdot}$, where $\delta$ comes from the differential on $\Graphs$. Hence the element $\co{\mu}{\gamma_1}$ represents the same cohomology class as $\delta \gamma_1$. But by \eqref{equ:Faltdef} and closedness of $\gamma$, we see that  $\delta \gamma_1=F(\gamma)$. This shows the proposition. 

\end{proof}

\section{\texorpdfstring{$\Der(\hoe_2)$ and $\grt_1$}{Der(hoe2) and grt}}
\label{sec:derhoe2grt}

In this section we will investigate more precisely the relation between $\grt_1$ and $\Der(\hoe_2)$ and at the end show that $H^0(\Der(\hoe_2))\cong \grt_1\oplus \GF$ as vector spaces. Theorem \ref{thm:hoe2} then easily follows in the next section, along with Theorem \ref{thm:main} as a Corollary by Theorem \ref{thm:GCdef} which has been proven in the preceding sections. We will begin by showing some auxiliary results.

As a corollary to Proposition \ref{prop:icgandt}, one obtains the following Lemma. 
\begin{lemma}
\label{lem:grandt}
\[
H^j(\gr^p \Def(\hoe_2 \to e_2)_\conn)
\cong
H^j(\gr^p \Def(\hoe_2 \to \hat{\alg{t}}[1]))
\]
where $\gr$ denotes the associated graded with respect to the filtration introduced in section \ref{sec:enfiltration}.
\end{lemma}
\begin{proof}
By Lemma \ref{lem:defpresqiso} and Propositions \ref{prop:Graphsincl} and \ref{prop:Graphscohom},
$\gr^p\Def(\hoe_2 \to e_2)_\conn$ is quasi-isomorphic to $\gr^p\Def(\hoe_2 \to \Graphs_2^\whl)_\conn$.
Note in particular that Lemma \ref{lem:defpresqiso} is applicable in this context since we may understand 
\[
 \prod_p \gr^p\Def(\hoe_2 \to e_2) \cong (\Def(\hoe_2 \to e_2), d_\wedge)
\]
as a deformation complex 
\[
 \Def(\hoe_2\stackrel{\tilde f} \longrightarrow e_2)
\]
where $\tilde f$ is the composition of the canonical map $\hoe_2\to e_2$ with the endomorphism of $e_2$ that sends the product operation to itself and the bracket to zero. See also Remark \ref{rem:tildef}.

By (the first part of) Proposition \ref{prop:fCconnqiso} the complex $\gr^p\Def(\hoe_2 \to \Graphs_2^\whl)_\conn$ is quasi-isomorphic to $(\fC_\conn, \delta)$. By Proposition \ref{prop:fCICG} in turn $H(\fC_\conn, \delta)$ is isomorphic to $H(\Def(\hoe_2 \rightarrow \ICG_2^\whl[1])_\conn, \delta+d_\wedge)$, which is isomorphic to $H(\Def(\hoe_2 \rightarrow \hat{\alg{t}}[1])_\conn)$ by Proposition \ref{prop:icgandt}.
Summarizing, we have the following zig-zag of quasi-isomorphisms of complexes.
\begin{multline*}
 \gr^p\Def(\hoe_2 \to e_2)_\conn \to \gr^p\Def(\hoe_2 \to \Graphs_2^\whl)_\conn \leftarrow 
\gr^p\fC_\conn \to
\\ \to \gr^p\Def(\hoe_2 \rightarrow \ICG_2^\whl[1])_\conn
\leftarrow
\gr^p\Def(\hoe_2 \to \hat{\alg{t}}[1])_\conn
\end{multline*}
\end{proof}

\begin{figure}
 \centering
\begin{tikzpicture}
 \draw (-1,-.5) -- coordinate (x axis mid) (4.5,-.5);
 \draw (-.5,-1) -- coordinate (y axis mid) (-.5,4.5);
\begin{scope}[yshift=-.5cm]
    	\foreach \x in {0,...,4}
     		\draw (\x,-1pt) -- (\x,-3pt)
			node[anchor=north] {\x};
\end{scope}
\begin{scope}[xshift=-.5cm]
    	\foreach \y in {0,...,4}
     		\draw (1pt,\y) -- (-3pt,\y) 
     			node[anchor=east] {\y}; 
 
\end{scope}
	%labels      
	\node[below=0.8cm] at (x axis mid) {$p$};
	\node[rotate=90, above=0.8cm] at (y axis mid) {$j-p$};
\node at (0,0) {$\K$};
\node at (0,1) {$\grt_1$};
\node at (1,0) {$\K$};
\node at (1,1) {$0$};
\node at (0,2) {$?$};
\node at (1,2) {$?$};
\node at (2,0) {$0$};
\node at (2,1) {$0$};
\node at (3,0) {$0$};
\node at (3,1) {$0$};
\node at (2,2) {$?$};
%\node at (3,1) {$0$};
%\node at (3,1) {$0$};
\begin{scope}[xshift=1.5cm]
\draw[dashed] (0,0)--(3,3);
\end{scope}
\end{tikzpicture}
\caption{\label{fig:cohomplot} Picture of the cohomology $H^j(\gr^p \Def(\hoe_2 \to e_2)_\conn)$, see Proposition \ref{prop:degrees} and Remark \ref{rem:victor}. The entries beyond the dashed line are 0. The entries above and including the question marks are unknown to the author. It is a well-known conjecture that the leftmost $?$ shown is $0$.}
\end{figure}
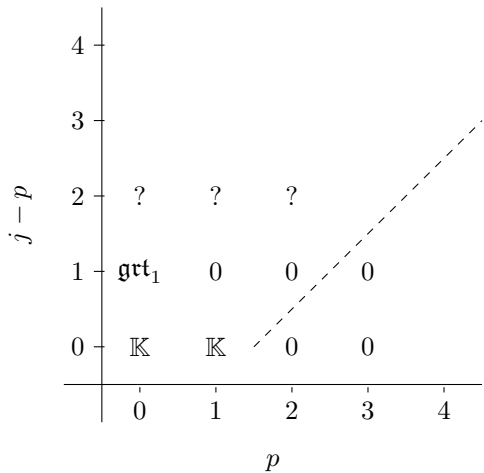

\begin{prop}
\label{prop:degrees}
\[
H^j(\gr^p \Def(\hoe_2 \to e_2)_\conn)
=
\begin{cases}
0 & \text{for $j<2p-1$ or $j \leq -1$} \\
\GF & \text{for $j=p=0$ } \\
\grt_1 & \text{for $j=1$, $p=0$ } \\
\GF & \text{for $j=p=1$ } \\
0 & \text{for $j=2$, $p=1$ } \\
\end{cases}
\]
See Figure \ref{fig:cohomplot} for an illustration of the situation.
\end{prop}

%In the cases not covered, the cohomology is not known to the author.
For the proof, recall the following fact from \cite{AT}.
\begin{prop}(Special case of Theorem 2.1 in \cite{AT})
\label{prop:AT}
Let $f(X,Y)\in \FreeLie(X,Y)$ be an element of the free Lie algebra in two generators such that 
\begin{equation}
\label{equ:fclosed}
f(Y,Z)-f(X+Y,Z)+f(X, Y+Z) - f(X,Y) = 0 \in \FreeLie(X,Y, Z).
\end{equation}
Then $f(X,Y)\propto \co{X}{Y}$. In particular, if in addition 
\begin{equation}
\label{equ:fsym}
f(X,Y)= f(Y,X)
\end{equation}
then $f(X,Y)=0$.
\end{prop}
In fact, this is one instance of a much more general (and simpler) combinatorial principle, see \cite{pavolcubical}.

\begin{proof}[Proof of Proposition \ref{prop:degrees}.]
We use Lemma \ref{lem:grandt}. An element of $\gr^p \Def(\hoe_2 \to \hat{\alg{t}}[1])_\conn$ which is in the $N$-th factor in the direct product decomposition \eqref{equ:convalgebra} has cohomological degree
\[
j = N + p -2.
\]
Since $N\geq p+1$ it follows that the cohomology vanishes for $j<2p-1$. Since the complex has no components of cohomological degree $j\leq -1$, the first line of the case distinction has been shown.
The second and third line are Lemmas \ref{lem:tharrtriv} and \ref{lem:tcohom}. 
The part of $\gr^p \Def(\hoe_2 \to \hat{\alg{t}}[1])_\conn$ of degrees $j=p=1$ is one dimensional. The generator corresponds to the graph with two vertices, each in their own cluster, and one edge. It is a cocycle, hence the fourth line follows.
For the final assertion let us again use the graphical language from Appendix \ref{sec:defengra}.
Suppose an element $x$ of cohomological degree $j=2$ and additional degree $p=1$ is given. Then, by the above formula, it has $N=3$ external vertices. They come in two clusters, one of size one and one of size two. We claim that if $x$ is closed, then $x=0$. Note that the differential of $x$ consists of two parts, one with two clusters of size two each, and one with one cluster of size one and one of size three. We claim that even if only the latter part vanishes, then $x=0$.
Concretely, $x$ is defined by an element of $T\in \hat{\alg{t}_3}$, symmetric under interchange of two indices, say of 1 and 2. It is well known that $\alg{t}_3\cong \FreeLie(t_{13},t_{23})\oplus \GF t_{12}$. By the connectedness requirement, the coefficient of $t_{12}$ in $T$ must be zero. Hence $T$ describes an element of the (completed) free Lie algebra in two symbols. The closedness condition for $x$ and the symmetry requirement translate into the conditions \eqref{equ:fclosed} and \eqref{equ:fsym} above. Hence by Proposition \ref{prop:AT} $x=0$.
\end{proof}

%In this section we will evaluate $H(\gr \Def(\hoe_n \to e_n)_\conn, d_H)$ in three different ways.
%First recall from \cite{} that there is a quasi-isomorphism of operads $e_n\to C(\alg{t}_n)$. From this it is easy to see that
%\[
%H(\gr \Def(\hoe_n \to e_n)_\conn)
%\cong
%H(\gr \Def(\hoe_n \to C(\alg{t}_n))_\conn).
%\]
%
%Next, note that $\alg{t}_n$ is quasi-isomorphic to $\ICG_n$ via maps $\alg{t}_n\leftarrow \TCG_n \rightarrow$.
%It follows that 
%\[
%H(\gr \Def(\hoe_n \to C(\alg{t}_n))_\conn)
%\cong
%H(\gr \Def(\hoe_n \to C(\ICG_n^\whl))_\conn).
%\]
%
%Analogously to Proposition \ref{prop:defengraphs}(??) the following result can be seen:
%\begin{lemma}
%\[
%H(\gr \Def(\hoe_n \to C(\ICG_n^\whl))_\conn)
%\cong
%H(\gr \Def(\hoe_n \to \ICG_n^\whl[1])_\conn).
%\]
%\end{lemma}
%\begin{proof}
%\end{proof}
%
%Furthermore as in Proposition \ref{prop:defengraphs}(??) one derives
%\begin{prop}
%\[
%H(\gr \Def(\hoe_n \to \ICG_n^\whl[1])_\conn)
%\cong
%H(\gr \Def(\hoe_n \to \alg{t}_n[1])_\conn).
%\]
%\end{prop}
%\begin{proof}
%\end{proof}
%
%As a corollary, one obtains degree bounds on the $H(\gr \Def(\hoe_n \to e_n)_\conn) $.
%\begin{cor}
%\label{cor:degrees}
%$H^j(\gr^N_{ar} \gr^p \Def(\hoe_n \to e_n)_\conn)=0$ if 
%$j-n(p-1) -N \notin (2-n)\N_0 $. Here $\gr^N_{ar}$ denotes the associated graded with respect to the arity filtration.
%\end{cor}
%\begin{proof}
%
%\end{proof}

%\subsection{The special case $n=2$}

\begin{cor}
\label{cor:hoe2e2}
\[
H^1(\Def(\hoe_2 \to e_2)_\conn)\cong \grt_1  \oplus \GF
\]
and 
\[
H^{0}(\Def(\hoe_2 \to e_2)_\conn)\cong \GF.
\]
Furthermore $H^{<0}(\Def(\hoe_2 \to e_2)_\conn)=0$.
\end{cor}
\begin{proof}
We consider the spectral sequence on the above complexes associated to the filtration of section \ref{sec:enfiltration}. It converges to cohomology by (for example) Lemma \ref{lem:findimdecomp}. 
A part of the $E^1$ term is computed in Proposition \ref{prop:degrees}.
Consider the piece $\grt_1$, living in degrees $p=0$, $j=1$ in the notation of that Proposition. Higher differentials in the spectral sequence must send it to the subspace of degrees $p\geq 1$, $j=2$. But by Proposition \ref{prop:degrees} all these vanish. The same reasoning works a fortiori for the summand $\GF$, living in degrees $p=1$, $j=1$.
The summand $\GF$, living in degrees $p=0$, $j=0$ might potentially be mapped to the summand $\GF$, living in degrees $p=1$, $j=1$ by the differential on $E^1$. However, looking at the explicit representative, which corresponds to a graph with two vertices in one cluster connected by an edge, one sees that this is not so. Finally, there is no cohomology in negative degrees even at the $E^1$ stage by Proposition \ref{prop:degrees}, hence the last assertion follows. 
\end{proof}

\begin{rem}\label{rem:victor}
 It was pointed out to the author by V. Turchin that using the above methods one can also show that 
\[
 H^2(\gr^3 \Def(\hoe_2 \to e_2)_\conn)=0.
\]
\end{rem}

\section{The proof of Theorems \ref{thm:main} and \ref{thm:hoe2}}
\label{sec:theproof}
In this section we show that $\grt_1 \cong H^0(\GC)$ as Lie algebras.\footnote{We adopt the convention that, if we omit the subscript $n$ in $\GC_n$, $\Graphs_n$ etc., $n=2$ is implied. In particular $\GC:=\GC_2$.} 
We saw in Theorem \ref{thm:GCdef} and in Corollary \ref{cor:hoe2e2} that
\[
H^0(\GC) \oplus \GF \cong H^0(\Der(\hoe_2))\cong H^1(\Def(\hoe_2 \to e_2)) \cong H^1(\Def(\hoe_2 \to e_2)_\conn)\cong \grt_1 \oplus \GF.
\] 
On both sides the class spanning $\GF$ corresponds to a graph with two vertices and the other summand to classes represented by sums of graphs with more than 2 vertices. Hence we may identify
\[
H^0(\GC)
\cong \grt_1
\]
as vector spaces. It remains to show that the isomorphism obtained is a map of Lie algebras. The inclusion $\grt_1 \to H^0(\Der(\hoe_2))$ is a map of Lie algebras by D. Tamarkin's construction, and the inclusion $H^0(GC) \to H^0(\Der(\hoe_2))$ is a map of Lie algebras by Theorem \ref{thm:GCdef}, and hence the first part of Theorem \ref{thm:main} is shown. Next, recall from Corollary \ref{cor:hoe2e2} that $H^0(\Def(\hoe_2 \to e_2))$ is one dimensional and that $H^{<0}(\Def(\hoe_2 \to e_2))=0$. Comparing with Theorem \ref{thm:GCdef}, the final assertion of Theorem \ref{thm:main}
\hfill \qed

Let us next turn to Theorem \ref{thm:hoe2}, which is almost shown by the above arguments. The only thing that remains to be verified is that the Lie bracket (in $H^0(\Der(\hoe_2))$) of $\grt_1$ elements with the additional class is as stated in the Theorem. The proof is similar to the proof of Proposition \ref{prop:firstmaplie}. The additional class corresponds to the derivation of $e_2$ given by a relative scaling of bracket and product.
Let $\xi$ be the derivation of $e_2$ that sends the product operation to zero and the bracket operation to itself. It acts by multiplication with minus the cohomological degree. Let $\Xi$ be the derivation of $\hoe_2$ that acts by multiplication with the grading on $\hoe_2$ introduced in section \ref{sec:enfiltration}. It is straightforward to check that $f\circ \Xi=\xi\circ f$ where $f: \hoe_2\to e_2$ is the projection. 
Let $\psi\in \grt_1$ be a homogeneous element of degree $r$ and let $\Psi \in \Der(\hoe_2))$ be a representative of the corresponding cohomology class in $\Der(\hoe_2))$. 
We need to show that
\[
\Xi\circ \Psi - \Psi\circ \Xi = r \Psi + (\text{exact terms}).
\]
Since $f$ is a quasi-isomorphism and by Lemma \ref{lem:defpresqiso} it suffices to show that 
\begin{equation}
\label{equ:fChiPsi}
f\circ \Xi\circ \Psi - f\circ\Psi\circ \Xi = r f\circ\Psi.
\end{equation}
Let $x$ be a generator of $\hoe_2$ of cohomological degree $d$ and of additional degree $k$. Note that automatically the arity of $x$ is fixed to be 
\[
N = 2-d-k.
\]
Apply both sides of \eqref{equ:fChiPsi} to $x$. We need to show that
\[
f\circ \Xi\circ \Psi(x) - f\circ\Psi\circ \Xi(x) = r f\circ\Psi(x).
\]
The left hand side evaluates to
\[
\xi\circ f\circ  \Psi(x) - f\circ\Psi\circ \Xi(x)
=
(- d - k) f\circ\Psi(x) = (N-2) f\circ\Psi(x).
\]

Noting that $f\circ\Psi(x)=0$ unless $r=N-2$ we are done.
\hfill \qed

\section{Explicit form of the (conjectural) generators}
\label{sec:explicitmap}
It is well known that $\grt_1$ contains a series of nontrivial elements $\sigma_3, \sigma_5,\dots$, such that $\sigma_{2j+1}$ is a Lie series in $X$, $Y$ starting with a multiple of the term $\ad_X^{2j}Y$. The Deligne-Drinfeld conjecture states that these elements freely generate $\grt_1$. In fact, it has recently been shown in \cite{brown} that these elements generate a free Lie subalgebra of $\grt_1$. However, it is still unclear whether this subalgebra is in fact the whole of $\grt_1$. 
In this section we will investigate how the graph cohomology classes corresponding to these generators look. The main result is the following.
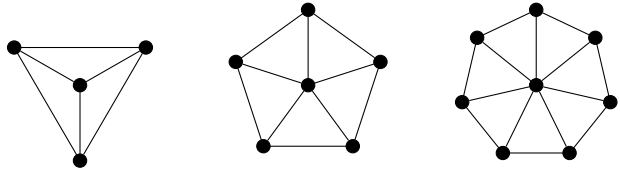
\begin{figure}
 \centering
 \[
\begin{tikzpicture}[
vert/.style={circle,draw,fill, minimum size=5pt, inner sep=0}, invvert/.style={inner sep=-1,minimum size=-1}, ext/.style={circle,draw, minimum size=5pt, inner sep=0}, scale=2 ]
\begin{scope}[shift={(-3.5,0)}]
\draw (30:0.5) node [vert] (v18) {} -- (150:0.5) node [vert] (v21) {} -- (-90:0.5) node [vert] (v19) {} -- (0:0) node [vert] (v20) {} -- (v18) -- (v19);
\draw  (v20) edge (v21);
\end{scope}
\begin{scope}[shift={(-2,0)}]
\draw (90:0.5) node [vert] (v22) {} -- (162:0.5) node [vert] (v26) {} -- (-126:0.5) node [vert] (v24) {} -- (-54:0.5) node [vert] (v1) {} -- (18:0.5) node [vert] (v23) {} -- (0:0) node [vert] (v25) {} -- (v22) -- (v23);
\draw (v24) -- (v25) -- (v26);
\end{scope}
\begin{scope}[shift={(-0.5,0)}]
\draw (90:0.5) node [vert] (v27) {} -- (141:0.5) node [vert] (v28) {} -- (-167:0.5) node [vert] (v30) {} -- (-116:0.5) node [vert] (v31) {} -- (-64:0.5) node [vert] (v32) {} -- (-13:0.5) node [vert] (v33) {} -- (39:0.5) node [vert] (v34) {} -- (v27) -- (0:0) node [vert] (v29) {} -- (v28) (v29) (v30) (v29) -- (v31) (v29) -- (v32) (v29) -- (v33) (v29) -- (v34) (v29) -- (v30);
\end{scope}
\draw (v25) -- (v1);
\end{tikzpicture}
\]
 \caption{\label{fig:wheelgraphs} Some wheel graphs. The graph cocyle corresponding to a (conjectural) generator $\sigma_{2j+1}$ of $\grt_1$ contains such a wheel with $2j+1$ spokes.}
\end{figure}

\begin{prop}
\label{prop:sigmawheel}
Under the map of Theorem \ref{thm:main}, the element $\sigma_{2j+1}\in \grt_1$ corresponds to a graph cohomology class, all of whose representatives have a nonvanishing coefficient in front of the wheel graph with $2j+1$ spokes, see Figure \ref{fig:wheelgraphs}. Moreover, the coefficient of that wheel graph equals minus the coefficient of the term $\ad_X^{2j}Y$ in the Lie series $\sigma_{2j+1}(X,Y)$, if for the wheel graph one uses the following ordering of edges to fix the sign
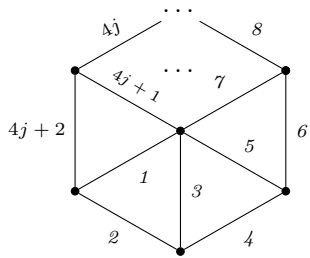
\[
 \begin{tikzpicture}
[scale=.8, baseline=-0.65ex, yshift=1.5cm
]
\node [int] (v7) at (0,0) {};
\node [int] at (30:2) {};
\node [int] (v4) at (-30:2) {};
\node [int] (v3) at (-90:2) {};
\node  (v6) at (90:2) {$\dots$};
\node  at (90:1) {$\dots$};
\node [int] (v1) at (150:2) {};
\node [int] (v5) at (30:2) {};
\node [int] (v2) at (210:2) {};
\draw[font=\scriptsize]  (v1) edge node[descr,auto, swap] {$4j+2$} (v2);
\draw[font=\scriptsize]  (v2) edge node[descr,auto,swap] {2} (v3);
\draw[font=\scriptsize]  (v3) edge node[descr,auto,swap] {4} (v4);
\draw[font=\scriptsize]  (v4) edge node[descr,auto,swap] {6} (v5);
\draw[font=\scriptsize]  (v5) edge node[descr,auto,swap] {8} (v6);
\draw[font=\scriptsize]  (v6) edge node[above, sloped] {$4j$} (v1);
\draw[font=\scriptsize]  (v7) edge node[descr, auto] {1} (v2);
\draw[font=\scriptsize]  (v7) edge node[descr,auto] {3} (v3);
\draw[font=\scriptsize]  (v7) edge node[descr,auto] {5} (v4);
\draw[font=\scriptsize]  (v7) edge node[descr,auto] {7} (v5);
\draw[font=\scriptsize]  (v7) edge node[descr, above, sloped, font=\tiny] {$4j+1$} (v1);
\end{tikzpicture}
\]
\end{prop}

Before we prove the proposition let us recollect the explicit form of the isomorphism $H^0(\GC_2)\to \grt_1$, which is scattered over the preceding sections. Let a graph cocycle $\gamma\in \GC_2$ be given. We want to find a way to read off the corresponding $\grt_1$-element from $\gamma$.
Remember that $\gamma$ is a linear combination of graphs with labeled vertices, invariant under permutations of the labels.
\\
\\
{\bf Algorithm 1:}
\begin{enumerate}
\item We assume that $\gamma$ is 1-vertex irreducible, which is possible by Proposition \ref{prop:1vi}.
 \item For each graph in $\gamma$, mark the vertex 1 as ``external''. This gives a (linear combination of) graph(s) $\gamma_{1}\in\Graphs(1)$. 
 \item Split the vertex 1 in $\gamma_1$ into two vertices, redistributing the incoming edges in all possible ways, so that both vertices are hit by at least one edge. Call this linear combination of graphs $\gamma_2' \in \CG_2(2)$.
 \item $\gamma_2'$ is closed in $\CG_2$ and has no one-edge component, hence it is the coboundary of some element $\gamma_2$. We choose $\gamma_2$ to be symmetric under interchange of the external vertices 1 and 2.
 \item Forget the non-internal-trivalent tree part of $\gamma_2$ to obtain $T_2$.
 \item \label{item:firstalgo6} For each tree $t$ occuring in $T_2$ construct a Lie word in (formal) variables $X,Y$ as follows. For each edge incident to vertex 1, cut it and make it the ``root'' edge. The resulting (internal) tree is a binary tree with leafs labelled by 1 or 2. It can be seen as a Lie tree, and one gets a Lie word $\phi_1(X,Y)$ by replacing each 1 by $X$ and 2 by $Y$. Set $\phi(X,Y)=\phi_1(X,Y)-\phi_1(Y,X)$. Summing over all such Lie words one gets a linear combination of Lie words corresponding to $\gamma$. Let us call it again $\phi_\gamma(X,Y)\in \F_\Lie(X,Y)$.
\item $\phi_\gamma$ is the desired $\grt_1$-element.
\end{enumerate}

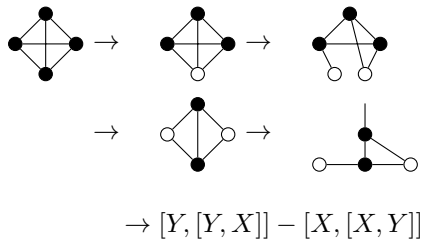
\begin{figure}
 \centering
 \[
\begin{tikzpicture}[
vert/.style={circle,draw,fill, minimum size=5pt, inner sep=0}, invvert/.style={inner sep=-1,minimum size=-1}, ext/.style={circle,draw, minimum size=5pt, inner sep=0}, scale=2 ]
\node (v0) at (0.5,1.8) [vert] {};
\node (v1) at (0.3,1.6) [vert] {};
\node (v2) at (0.5,1.4) [vert] {};
\node (v3) at (0.7,1.6) [vert] {};
\node (v4) at (1.5,1.4) [ext] {};
\node (v5) at (1.5,1.8) [vert] {};
\node (v6) at (1.3,1.6) [vert] {};
\node (v7) at (1.7,1.6) [vert] {};
\node (v8) at (0.9,1.6) [invvert] {$\to$};
\node (v9) at (1.9,1.6) [invvert] {$\to$};
\node (v10) at (2.4,1.4) [ext] {};
\node (v11) at (2.6,1.4) [ext] {};
\node (v12) at (2.3,1.6) [vert] {};
\node (v13) at (2.5,1.8) [vert] {};
\node (v14) at (2.7,1.6) [vert] {};
%\node (v15) at (3.3,1.6) [ext] {};
%\node (v16) at (3.7,1.6) [ext] {};
%\node (v17) at (3.5,1.8) [vert] {};
%\node (v18) at (3.5,1.4) [vert] {};
%\node (v19) at (2.9,1.6) [invvert] {$\to$};
\node (v24) at (1.3,1) [ext] {};
\node (v21) at (1.5,1.2) [vert] {};
\node (v22) at (1.5,0.8) [vert] {};
\node (v23) at (1.7,1) [ext] {};
\node (v25) at (0.9,1) [invvert] {$\to$};
\node (v26) at (1.9,1) [invvert] {$\to$};
\node (v27) at (2.3,0.8) [ext] {};
\node (v28) at (2.9,0.8) [ext] {};
\node (v29) at (2.6,0.8) [vert] {};
\node (v30) at (2.6,1) [vert] {};
\node (v31) at (2.6,1.2) [invvert] {};
\node (v32) at (2,0.4) [invvert] {$\to[Y,[Y,X]] - [X,[X,Y]]$};
\draw (v0)--(v2) (v2)--(v3) (v3)--(v1) (v1)--(v2) (v3)--(v0) (v0)--(v1) (v6)--(v4) (v4)--(v7) (v7)--(v5) (v5)--(v6) (v6)--(v7) (v4)--(v5) (v11)--(v14) (v14)--(v12) (v12)--(v13) (v13)--(v14) (v12)--(v10) (v13)--(v11);
% (v15)--(v16) (v17)--(v18) (v18)--(v16) (v16)--(v17) (v17)--(v15) (v15)--(v18) 
\draw (v23)--(v21) (v21)--(v22) (v22)--(v23) (v24)--(v22) (v24)--(v21) (v31)--(v30) (v30)--(v29) (v29)--(v27) (v29)--(v28) (v30)--(v28) ;
\end{tikzpicture}
\]
 \caption{\label{fig:firstalg} Illustration of the algorithm mapping a graph cohomology class to a $\grt$-element. Prefactors are omitted. The graphs depicted are (in this order) $\Gamma$, $\Gamma_1$, $\Gamma_2'$, $\Gamma_2=T$ and the Lie tree.}
\end{figure}

A graphical illustration of this algorithm for the simplest case of a three-wheel (corresponding to $\sigma_3\in \grt$) is given in Figure \ref{fig:firstalg}.

\begin{prop}
\label{prop:firstalgoagrees}
 The above Algorithm 1 produces the correct result, i.~e., $\phi_\gamma$ is indeed the $\grt_1$-element corresponding to $\gamma$ under the map of Theorem \ref{thm:main}.
\end{prop}

\begin{proof}
We proceed as in the proof of Theorem \ref{thm:main} in the previous section. There we saw that the image of $\gamma$ in $\Def(\Com_\infty\to \Graphs_2)$ is the element $d_H \gamma_1=\gamma_2'$. It can be seen as an element of $\Graphs_2$, symmetric under interchange of the external vertices. It furthermore consists entirely of graphs with only one internally connected component by 1-vertex-irreducibility of $\gamma$. Now let us compute the element of $\grt_1$, corresponding to the cohomology class of $d_H \gamma_1$. First we pick a different representative, which will be $d_H \gamma_2$. The latter can be seen as a cocycle in $\CG(3)$ satisfying some symmetry property. But any cocycle in $\CG(3)$ represents an element of $\alg{t}_3$, which can be obtained by restricting to the part consisting of trivalent internal trees. From such a tree one can recover the $\grt_1$-element as described in Appendix \ref{sec:sder}.
% To map $\Gamma$ to $\Def(\hoe_2\to \Graphs)$ we need to compute its action on the Maurer-Cartan element (see Figure \ref{fig:}). For this one computes $\Gamma_1\in \Graphs(1)$ by marking vertex 1 as external.Then one computes the bracket with the MC element, effectively splitting vertex 1 into two external vertices in two ways. Once the external vertices are connected by an edge and once they are not, but reside in the same cluster (see Appendix \ref{sec:derengra} for an explanation of these terms). We call the former term $A$ and the latter term $B$. The projection to $\Def(\Com_\infty\to \Graphs)$ kills the $A$-part and retains only $B$. Let $\Gamma_{12}$ be as in the algorithm. From Appendix \ref{sec:} one obtains the following claim.
% \\
% Claim: $\delta \Gamma_{12} = A$, where $\delta$ is the differential in $\Graphs$.
% \\
% Next let $\Gamma_{12}'$ from $\Gamma_{12}$ by (i) forgetting all parts that do not have the vertices 1 and 2 connected by an edge and (ii) remove this edge between 1 and 2 and put 1 and 2 into a cluster. Call the joint operation (i) and (ii) the ``1-2-edge-removal''.
% \\
% Claim: $\delta \Gamma_{12}' = B$
% \\
% Proof: The claim is obtained from the previous claim by noting that $A$ has vertices 1 and 2 connected by an edge. $B$ can be obtained from $A$ by 1-2-edge-removal. Since the differential commutes with 1-2-edge-removal, the claim follows.
% 
% For the final step of the algorithm, i.e., the extraction of $\phi_\Gamma(x,y)$ from $\Gamma_{12}'\in \Def(\Com_\infty\to \Graphs)$, we refer to Appendix \ref{sec:sder}.
\end{proof}

Now let us turn to the proof of Proposition \ref{prop:sigmawheel}. All we need to know about $\sigma_{2j+1}$ is that it contains a term 
\[
\ad_X^{2j}(Y)
\]
with non-trivial coefficient.
We need to trace through the above algorithm, and see that the wheel graph must be present in $\gamma$ in order to produce that term. 

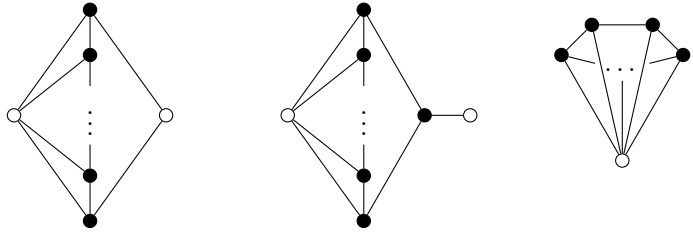
\begin{figure}
\centering
\[
\begin{tikzpicture}[
vert/.style={circle,draw,fill, minimum size=5pt, inner sep=0}, invvert/.style={inner sep=-1,minimum size=-1},ext/.style={circle,draw, minimum size=5pt, inner sep=0}, scale =2 ]
\begin{scope}[xshift=-.4cm]
\node (v0) at (0.2,1.1) [ext,label=180:{}] {};%1
\node (v1) at (1.2,1.1) [ext,label=360:{}] {};%2
\node (v2) at (0.7,1.8) [vert] {};
\node (v3) at (0.7,1.5) [vert] {};
\node (v4) at (0.7,0.4) [vert] {};
\node (v5) at (0.7,0.7) [vert] {};
\node (v6) at (0.7,1.1) {\vdots};
\end{scope}

\node (v7) at (2.5,1.1) [vert] {};
\node (v8) at (2.8,1.1) [ext,label=360:{}] {};%2
\node (v9) at (2.1,1.8) [vert] {};
\node (v10) at (3.4,1.5) [vert] {};
\node (v11) at (2.1,0.7) [vert] {};
\node (v12) at (2.1,0.4) [vert] {};
\node (v17) at (2.1,1.1) {\vdots};
\node (v14) at (3.8,0.8) [ext,label=270:{}] {}; %1
\node (v15) at (2.1,1.5) [vert] {};
\node (v16) at (1.6,1.1) [ext,label=180:{}] {};%1
\node (v18) at (3.6,1.7) [vert] {};
\node (v19) at (4,1.7) [vert] {};
\node (v20) at (4.2,1.5) [vert] {};
\node (v22) at (3.8,1.4) {\dots};
\draw (v0)--(v3) (v2)--(v3) (v0)--(v2) (v0)--(v5) (v0)--(v4) (v5)--(v4) (v5)--(v6) (v3)--(v6) (v4)--(v1) (v2)--(v1) (v16)--(v11) (v16)--(v15) (v9)--(v15) (v9)--(v16) (v16)--(v12) (v12)--(v11) (v12)--(v7) (v15)--(v17) (v11)--(v17) (v9)--(v7) (v7)--(v8) (v10)--(v14) (v10)--(v18) (v18)--(v19) (v19)--(v20) (v14)--(v20) (v18)--(v14) (v14)--(v19) (v10)--(v22) (v22)--(v20) (v22)--(v14) ;
\end{tikzpicture}
\]
\caption{\label{fig:firstproof} Three graphs occuring in the proof of Proposition \ref{prop:sigmawheel}.}
\end{figure}

\begin{proof}[Proof of Proposition \ref{prop:sigmawheel}.]
In order to produce a term $\ad_X^{2j}(Y)$ at the end, the $T_2$ in the algorithm, and hence also the $\gamma_2$, has to contain a graph as depicted in Figure \ref{fig:firstproof} (left). But this means that in $\gamma_2'$ there is a term of the form depicted in Figure \ref{fig:firstproof} (middle). But such a term can only be produced if $\gamma_1$ contains a term like that in Figure \ref{fig:firstproof} (right). But this means that $\gamma$ has to contain a wheel graph.

Next, let us track the coefficients. Suppose the wheel graph appears in the graph cocyle with coefficient 1.
Recall that by our conventions this means that the corresponding linear combination of numbered graphs is the sum over all $(2j+2)!$ numberings, divided by the order of the symmetry group $2(2j+1)$.
Hence, the graph on the right of Figure \ref{fig:firstproof} appears with coefficient 1 as well.
The graph in the middle of Figure \ref{fig:firstproof} can be produced from the splitting of vertex 1 in $2j+1$ ways, but at the same time the symmetry group is only of order 2. Hence the coefficient in front of this graph is 1 as well, as is the coefficient of the graph on the left hand side.
Finally we have to create a Lie tree by detaching an edge from one of the vertices. Here, we have to pick the lower valence vertex to obtain the desired Lie word. There are two edges to detach yielding the same contribution. However, together with the $\frac 1 2$ from the symmetry factor we obtain one contribution corresponding to a Lie tree of the form
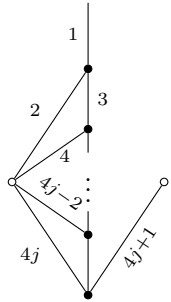
\[
\begin{tikzpicture}
\node[ext] (1) at (0,0) {};
\node (2) at (1,2.5) {};
\node[ext] (3) at (2,0) {};
\node[int] (w1) at  (1,1.5) {};
\node[int] (w2) at  (1,.7) {};
\node (wd) at  (1,0) {$\vdots$};
\node[int] (w3) at  (1,-.7) {};
\node[int] (w4) at  (1,-1.5) {};
\draw (1) edge node[auto] {$\scriptstyle 2$} (w1) 
          edge node[right] {$\scriptstyle 4$} (w2) 
          edge node[sloped, above] {$\scriptstyle 4j-2$} (w3) 
          edge node[below left] {$\scriptstyle 4j$} (w4) 
      (w1) edge node[auto] {$\scriptstyle 1$} (2) edge node[auto] {$\scriptstyle 3$} (w2) 
      (wd) edge (w2) edge (w3)
      (w4) edge node[sloped,below] {$\scriptstyle 4j+1$} (3) %
           edge node[sloped,above] {} (w3); %$\scriptstyle 4j-1$
\end{tikzpicture}%\, .
\]
which contributes to $\phi_1(Y,X)$
If we map this graph into $\alg t_3$, we obtain $\ad_X^{2j}(Y)$.
See the example in Appendix \ref{sec:sder} for the sign convention used when identifying a tree with a Lie word.
\end{proof}

For reference, let us also give an explicit algorithm for the opposite map $\grt_1\to H^0(\GC_2)$.
Let $\phi\in \grt_1$. We can assume without loss of generality that $\phi$ is homogeneous of degree $n$ with respect to the grading on $\grt_1$. I.~e., $\phi=\phi(X,Y)$ can be seen as a Lie expression in formal variables $X,Y$, with $X,Y$ occurring $n$ times. Consider the complex
\[
 D:=\Def(\hoLie_2 \to C(\alg{t}))
\]
where $C(\alg{t})$ is the Chevalley complex of the Drinfeld-Kohno (operad of) Lie algebra(s) $\alg{t}$, with trivial coefficients.
\begin{rem}
\label{rem:holietcohom}
Since $C(\alg{t})$ is quasi-isomorphic to $e_2$ (see \cite{tamarkin}), the complex $D$ has no cohomology in degree zero or lower already by degree reasons. In fact it is acyclic, as can be seen from Proposition \ref{prop:Xiqiso} in Appendix \ref{sec:endef}.
\end{rem}
We have the following algorithm, in which we ignore the overall sign:
\\
\\
{\bf Algorithm 2:}
\begin{enumerate}
\item Symmetrize (and shift in degree) the element $\phi(t_{12},t_{23})\wedge t_{12}\wedge t_{23} \in C(\alg{t}_3)$ so as to obtain an element $T_3\in D$. It is a cocycle (to be shown below).
 \item Find an element $U\in D$ with coboundary $T_3$. More concretey, split $U=U_3+U_4+\dots+U_{n+1}$ according to the grading on $D$ by arity in $\alg{t}$ (i.e., by ``number of external vertices'' if one thinks in terms of graphs). Then $d_{CE}U_3=T_3$, $d_{\co{}{}}U_3=-d_{CE}U_4$ etc. where $d_{CE}$ is the part of the differential on $D$ coming from the Chevalley-Eilenberg differential on $C(\alg{t})$ and $d_{CE}+d_{\co{}{}}$ is the full differential.
 \item By degree reasons, $U_{n+1}\in C(\alg{t})^{S_{n+1}}[-2n]\subset D$ is a linear combination of wedge products in the generators of $\alg{t}_{n+1}$ (i.e., the $t_{ij}$'s). Replacing each $t_{ij}$ by an edge between vertices $i$ and $j$ one obtains a linear combination of graphs $\gamma'$.
 \item Drop all graphs in $\gamma'$ containing vertices of valence smaller than three. This gives some element $\gamma\in \GC_2$, the result. It is closed since $d_{\co{}{}}U_{n+1}=0$ by construction.
\end{enumerate}

\begin{figure}
 \centering
 \[
\begin{tikzpicture}[
vert/.style={circle,draw,fill, minimum size=5pt, inner sep=0}, invvert/.style={inner sep=-1,minimum size=-1}, ext/.style={circle,draw, minimum size=5pt, inner sep=0}, scale=2 ]
\node (v0) at (0.3,2.4) [ext] {};
\node (v1) at (0.7,2.4) [ext] {};
\node (v2) at (1.1,2.4) [ext] {};
\node (v3) at (0.7,2.9) [vert] {};
\node (v4) at (0.7,2.7) [vert] {};
\node (v5) at (1.4,2.7) [invvert] {$\to$};
\node (v6) at (1.9,2.4) [ext] {};
\node (v7) at (2.3,2.4) [ext] {};
\node (v8) at (2.7,2.4) [ext] {};
\node (v10) at (2.3,2.7) [vert] {};
\node (v11) at (1.4,1.8) [invvert] {$\to$};
\node (v12) at (1.8,1.6) [ext] {};
\node (v13) at (2.1,1.6) [ext] {};
\node (v14) at (2.4,1.6) [ext] {};
\node (v15) at (2.7,1.6) [ext] {};
\node (v16) at (2.2,2) [vert] {};
\node (v17) at (3.2,1.6) [ext] {};
\node (v18) at (3.5,1.7) [ext] {};
\node (v19) at (3.8,1.6) [ext] {};
\node (v20) at (4.1,1.6) [ext] {};
\node (v21) at (3.6,2) [vert] {};
\node (v22) at (2.9,1.8) [invvert] {$+$};
\node (v23) at (1.4,1.2) [invvert] {$\to$};
\node (v24) at (1.8,1) [ext] {};
\node (v25) at (2.4,1) [ext] {};
\node (v26) at (2.7,1) [ext] {};
\node (v27) at (2.1,1.2) [ext] {};
\draw (v0)--(v1) (v1)--(v2);
\draw (v0)--(v3) (v4)--(v0) (v4)--(v3) (v4)--(v1) (v3)--(v2) (v6)--(v10) (v10)--(v7) (v10)--(v8) (v8)--(v7) (v7)--(v6) (v12)--(v16) (v16)--(v13) (v16)--(v15) (v12)--(v13) (v13)--(v14) (v14)--(v15) (v21)--(v18) (v17)--(v21) (v21)--(v20) (v18)--(v19) (v19)--(v20) (v17)--(v19) (v24)--(v25) (v24)--(v27) (v26)--(v27) (v27)--(v25) (v25)--(v26) ;
\draw  (v15) to [bend left] (v12)  (v24) to [bend right] (v26) (v17) to[bend right] (v20) (v6) to [bend right] (v8);
\end{tikzpicture}
\]
 \caption{\label{fig:secondalg} Illustration of Algorithm 2, mapping the $\grt$-element $\sigma_3$ to a graph cohomology class. Prefactors are omitted. The graphs depict (in this order) $T_3$, $U_3$, $d_{\co{}{}}U_3$, $U_4=\gamma$. It is explained in Appendix \ref{sec:sder} how to depict elements of $C(\alg{t})$ as graphs.}
\end{figure}
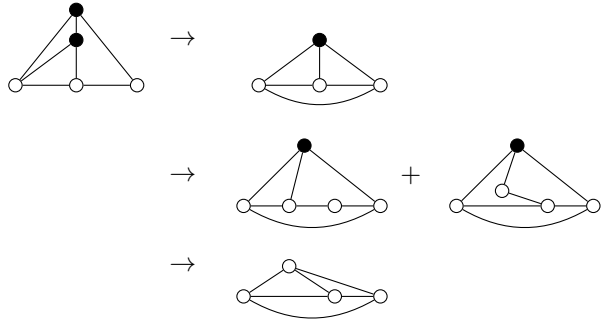

A graphical illustration of this Algorithm for $\phi=\sigma_3$ can be found in Figure \ref{fig:secondalg}.

To make the algorithm work, the following reesult has to be shown.
\begin{lemma}
The element $T_3$ constructed in the first step is indeed a cocycle.
\end{lemma}
\begin{proof}
First let us show that $T_3$ is closed under the part $d_{CE}$ of the differential.
Concretely, one has
\[
T_3\propto \phi(t_{12},t_{23})\wedge t_{12}\wedge t_{23} 
 + \phi(t_{23},t_{31})\wedge t_{23}\wedge t_{31}
 + \phi(t_{31},t_{12})\wedge t_{31}\wedge t_{12}
\]
and hence 
\begin{align*}
d_{CE}T_3 &\propto - \phi(t_{12},t_{23})\wedge \co{t_{12}}{t_{23}} + \text{(cycl.)} \\
& \quad + \co{\phi(t_{12},t_{23})}{t_{12}}\wedge t_{23} - \co{\phi(t_{12},t_{23})}{t_{23}}\wedge t_{12}+ \text{(cycl.)} \\
&=
0 +\left( \co{\phi(t_{12},t_{23})}{t_{12}} - \co{\phi(t_{23},t_{31})}{t_{31}}\right) \wedge t_{23} + \text{(cycl.)} \\
&= 0
\end{align*}
Here the term in the first line on the right is zero because of the symmetries of $\phi(t_{12},t_{23})$. In more detail, $\co{t_{12}}{t_{23}}$ is antisymmetric under the $S_3$ action on indices. Hence the symmetrization of $\phi(t_{12},t_{23})\wedge \co{t_{12}}{t_{23}}$ picks out the the antisymmetric part of $\phi(t_{12},t_{23})$, which is zero by the hexagon equation.
The second equality follows directly from the cabling relation (or equivalently, the semiclassical hexagon), which in turn follows from the hexagon equation (see \cite{barnatan} for a proof and more details on those relations).

Next we need to show that $T_3$ is also closed under $d_{\co{}{}}$. This is tedious if one writes everything out in detail.
Gathering similar terms one obtains
\begin{align*}
d_{\co{}{}}T_3 &\propto 
(A)\wedge t_{12}\wedge t_{23}\wedge t_{34}
+(B)\wedge t_{12}\wedge t_{23}\wedge t_{24}
+(\dots)
\end{align*}
where the terms $(\dots)$ can be obtained from the first two terms by permutations of indices, so as to make whole expression symmetric. One calculates 
\[
A=-\phi^{2,3,4}-\phi^{1,2,3}+\phi^{12,3,4}-\phi^{1,23,4}+\phi^{1,2,34}=0
\]
by the pentagon equation. Here we use the usual notation $\phi^{12,3,4}=\phi(t_{13}+t_{23},t_{34})$ etc. Similarly one computes
\begin{align*}
B&= -\phi^{1,2,3}-\phi^{3,2,4}+\phi^{1,2,4}+\phi^{1,24,3}-\phi^{1,23,4}+\phi^{3,12,4} \\
&=-\phi^{1,2,3}-\phi^{3,2,4}-\phi^{2,4,3}+\phi^{12,3,4}-\phi^{1,23,4}+\phi^{1,2,34} \\
&= -\phi^{3,2,4}-\phi^{2,4,3}+\phi^{2,3,4} =0.
\end{align*}
Here we used twice the pentagon and once the hexagon equation. Hence $d_{\co{}{}}T_3=0$.

\end{proof}

Now that we know that $T_3$ is closed, we can construct $U$, which exists and is unique up to exact terms by the remark preceding Algorithm 2. % \ref{rem:holietcohom}.
One can verify (see Appendix \ref{sec:tamalgodoagree}) that the output of Algorithm 2 agrees with the Tamarkin map from $\grt_1$ to $H^0(\GC_2)$ (see Appendix \ref{sec:tamalgoagree}) and hence also with the map of Theorem \ref{thm:main}.

\section{Applications}
\label{sec:appdefq}
In this section we will discuss applications of the results obtained in this paper.
We suppose that the reader is already familiar with the basic objects and questions of deformation quantization. If not, we refer to M. Kontsevich's seminal paper \cite{K1}.
\subsection{The action of the graph complex on formality morphisms}\label{sec:grcactionfm}
Let $T_{\rm poly}^\bullet=\Gamma(\R^d; \wedge^\bullet T\R^d)$ be the space of multivector fields on $\R^d$. There is an action of the operad $\Gra_2$ on $T_{\rm poly}$, given by \eqref{equ:Graaction}.
%It is a Gerstenhaber algebra with the wedge product and the Schouten bracket. Furthermore, there is an action of the operad $\Gra_2$ on $T_{\rm poly}$, considered as a vector space. Concretely, for a graph $\Gamma\in \Gra_2(N)$ and $\gamma_1,\dots, \gamma_N\in T_{\rm poly}$ the action reads
%\[
% \Gamma(\gamma_1,\dots,\gamma_N) = \mu\circ \left( \prod_{(i,j)}\sum_{k=1}^d \pd{}{x_k^{(j)}} \pd{}{\xi_k^{(i)}} 
%+ \pd{}{x_k^{(i)}} \pd{}{\xi_k^{(j)}}
%\right)\left(\gamma_1\otimes\cdots\otimes\gamma_n \right).
%\]
%Here $\mu$ is the operation of multiplication of $n$ multivector fields and the product runs over all edges $(i,j)$ in $\Gamma$, in the order given by the numbering of edges. The notation $\pd{}{x_k^{(j)}}$ means that the partial derivative is to be applied to the $j$-th factor of the tensor product, and similarly for $\pd{}{\xi_k^{(i)}}$.
Hence $\Tpoly$ is also a Gerstenhaber algebra via the operad map \eqref{equ:enGramapdef}.
It also follows that there is an action of the dg Lie algebra $\Def(\hoLie_2\to \Gra_2)$ on $T_{\rm poly}[1]$ by pre-$L_\infty$-derivations. Here the closed degree zero elements act by true $L_\infty$-derivations. In particular, the closed degree zero elements of $\GC_2\subset \Def(\hoLie_2\to \Gra_2)$ act in this way. By the identification $\grt_1\cong H^0(\GC_2)$ there is also an action of $\grt_1$ on $T_{\rm poly}$, defined up to homotopy. 

Any such action on $\Tpoly$ also yields an action on the space of its Maurer-Cartan elements, i.~e., on Poisson structures, see also \cite{K3}. 

Furthermore, let $\Dpoly$ be the space of polydifferential operators on $\R^d$. The main result of deformation quantization is M. Kontsevich's Formality Theorem, stating the existence of an $L_\infty$ quasi-isomorphism 
\[
 \Tpoly[1] \to \Dpoly[1].
\]
In the following, such a morphism will be called a formality morphism. By composition, closed degree zero elements of $\Def(\hoLie_2\to \Gra_2)$ and $\GC_2$ and $\grt_1$ act on the space of formality morphisms. Notice that in this manner one obtains a \emph{right} action on formality morphisms, and one has to flip the sign to obtain a left action.

\subsection{Recollection: D. Tamarkin's proof of the Formality Theorem and the \texorpdfstring{$\GRT_1$}{GRT1} action}\label{sec:tamconstr}
M. Kontsevich originally proved his Formality Theorem by writing down an explicit formula, involving certain hard-to-compute configuration space integrals. D. Tamarkin later gave a different proof \cite{tamdefq, hinich} as follows. First, one endows $\Dpoly$ with a $\hoe_2$ structure, whose $\hoLie_2$ part agrees with the usual one. This step depends on the choice of a Drinfeld associator. Secondly, by homotopy transfer one knows that $\Dpoly$ with that $\hoe_2$ structure is $\hoe_2$-quasi-isomorphic to its cohomology $\Tpoly$, with some non-standard $\hoe_2$ structure. Call this latter space $\Tpoly'$ to distinguish it from $\Tpoly$ with the standard $\hoe_2$ structure. Finally, D. Tamarkin shows a rigidity result, which states that there are no obstructions for the existence of a $\hoe_2$-quasi-isomorphism $\Tpoly\to \Tpoly'$.
The $\hoLie_2$ part of the composition $\Tpoly\to \Tpoly'\to \Dpoly$ is the sought after formality morphism.

There are multiple ways of endowing $\Dpoly$ with a $\hoe_2$ structure. Contrary to \cite{tamdefq} we will use a formality morphism of the operad of chains of the little disks operad as in \cite{tamarkin}. The latter operad is quasi-isomorphic to the braces operad $\Br$, which naturally acts on $\Dpoly$ \cite{KS1}. One may write down a chain of quasi-isomorphisms of operads
\[
\Br \to \cdots \to C\Ne \PaCD \to BU\alg{t} \leftarrow %C(\alg{t}) 
e_2 \leftarrow \hoe_2.
\]
The operads occurring in this chain are explained in more detail in \cite{tamarkin,graphthings,pavol} and are unimportant here.
Lifting (up to homotopy) yields a map $\hoe_2\to \Br$ such that the restriction $\hoLie_2\to \Br$ is the usual map and hence one obtains the desired $\hoe_2$ structure on $\Dpoly$. A Drinfeld associator (cf. Section \ref{sec:grt}) is used for constructing the map to $C\Ne \PaCD$ in the above chain, see \cite{tamarkin} for details.

The Grothendieck-Teichm\"uller group $\GRT_1$ (and its Lie algebra $\grt_1$) act on $\PaCD$ and hence on $C\Ne \PaCD$ in the above chain. Equivalently, $\GRT_1$ acts on the set of Drinfeld associators, one element of which is chosen to produce the map into $C\Ne \PaCD$. Still equivalently (up to homotopy) $\GRT_1$ and $\grt_1$ act on $\hoe_2$ up to homotopy by transfer. The latter action is the one from \cite{tamarkingrtaction}, recalled in section \ref{ssec:tamactionsginfty}.

Now suppose we have fixed a Drinfeld associator $\Phi_1$ yielding a map $F_1:\hoe_2\to\Br$. Suppose we obtain another associator $\Phi_2$ by the action of $g\in \GRT_1$. Denote the automorphism of $\hoe_2$ corresponding to $g$ by $g$ as well, abusing notation.
From $F_1$ and $F_2:=F_2\circ g$ we obtain two $\hoe_2$ structures on $\Dpoly$, say $\Dpoly$ and $\Dpoly'$, and two essentially unique homotopy Gerstenhaber quasi-isomorphisms $\Tpoly\to \Dpoly$, $\Tpoly\to \Dpoly'$.
The map $\Tpoly\to \Dpoly'$ may be chosen to fit into a commutative diagram
\[
\begin{tikzpicture}
 \matrix[diagram](m) {\Tpoly & \Dpoly \\ \Tpoly' & \Dpoly' \\};
 \draw[->] (m-1-1) edge node[above] {$f_1$} (m-1-2) edge node[left] {$g'$} (m-2-1) edge node[above]{$f_2$} (m-2-2) (m-2-1) edge node[above] {$f_2'$} (m-2-2);
\end{tikzpicture}
\]
where $\Tpoly'$ is obtained from $\Tpoly$ by acting on the standard $\hoe_2$ structure with $g$, and $g':\Tpoly \to \Tpoly'$ is the essentially unique homotopy Gerstenhaber isomorphism.

Of course, as $\hoLie_2$ algebras $\Tpoly=\Tpoly'$ and $\Dpoly=\Dpoly'$.
Restricting to the $\hoLie_2$ parts of $f_1$ and $f_2$ one obtains $L_\infty$ quasi-isomorphisms
\[
\tilde f_1, \tilde f_2 \colon \Tpoly[1] \to \Dpoly[1].
\] 
If we denote the $L_\infty$ parts of $g'$ and $f_2'$ by $\tilde g'$ and $\tilde f_2'$, one can then check that possibly up to homotopy $\tilde f_2'=\tilde f_1$ and hence
\[
 \tilde f_1 = \tilde f_1\circ \tilde g'.
\]
The infinitesimal version of this action is the one described in section \ref{sec:tamacttpoly}.

\begin{rem}
 It has been pointed out to the author by V. Dolgushev that the above statement is not obvious because the morphisms considered are homotopy or ($\infty$-)morphisms and not strict morphisms of $\hoe_2$ algebras. If $A\to B$ is a strict quasi-isomorphism of $\hoe_2$ algebras and $A'$, $B'$ are the same objects but with $\hoe_2$ structure altered by $g$ as above, then clearly the same map induces a strict morphism $A'\to B'$. In particular, if the $\hoLie_2$ structure is unaltered then the resulting strict morphisms of $\hoLie_2$ algebras $A\to B$ trivially agree.
The analog statement for homotopy morphisms instead of strict ones is not obvious. However, we may understand each such quasi-isomorphism $A\to B$ as a zigzag of strict quasi-isomorphisms $A \leftarrow X\to B$. Then, changing all $\hoe_2$ structures by $g$, we obtain a zigzag of strict quasi-isomorphisms $A' \leftarrow X'\to B'$.
If the $\hoLie_2$ part of the structure is unaltered, one again obtains trivially the same homotopy $\hoLie_2$ morphism.
In our setting the quasi-isomorphism $f_2':\Tpoly'\to \Dpoly'$ is essentially unique (up to homotopy) and hence the induced $L_\infty$ morphism $\tilde f_2'$ must agree (up to homotopy) with $\tilde f_1$. An alternative approach is described in \cite{paljug}.
\end{rem}

\subsection{Relation of the Tamarkin and Kontsevich formality morphisms}
% Let $\Phi$ be a Drinfeld associator, i.e., a power series in some variables $X,Y$, which is group-like and satisfies the pentagon, hexagon and antisymmetry conditions as in Section \ref{sec:grt}. 
% To $\Phi$ one can associate a formality morphism
% \[
%  \mU_\Phi \colon \Tpoly \to \Dpoly.
% \]
% In fact, there are multiple, conjecturally equivalent ways to associate the morphism $\mU_\Phi$ to the Drinfeld associator $\Phi$. Let us choose here D. Tamarkin's construction using formality of the little disks operad. In any way, one obtains many inequivalent formality morphisms, one for eeach associator. On the other hand, M. Kontsevich's original proof featured a very different construction of one particular formality morphism $\mU_K$, using graphical techniques. This construction contains essentially no ``free parameters''. It has been an open question for some time how D. Tamarkin's and M. Kontsevich's constructions of Formality Morphisms fit together.

D. Tamarkin's proof of the Kontsevich Formality Theorem depends on the choice of a Drinfeld associator. Different choices of associator yield different formality quasi-isomorphisms. On the other hand, M. Kontsevich's original proof featured a very different construction of one particular formality morphism $\mU_K$, using graphical techniques. This construction contains essentially no ``free parameters''. It has been an open question for some time how D. Tamarkin's and M. Kontsevich's constructions of formality morphisms fit together.

I think I can now give an answer, to some extend. First it is shown in \cite{graphthings} that D. Tamarkin's formality morphism (more precisely the version of section \ref{sec:tamconstr}) is homotopic to the one constructed by M. Kontsevich \cite{K1}, if one chooses for the associator the Alekseev-Torossian associator $\Phi_{AT}$ \cite{pavol, ATassoc}. 
% Secondly, the is an action of $\grt_1$ on the set of Drinfeld associators according to the formula
% \[
%  (\psi \cdot \Phi)(X,Y) =  \left. \frac{d}{dt}\right|_{t=0} \exp(t\phi) \Phi(X,\exp(-t\phi)Y\exp(t\phi))
% \]
% for $\phi\in \grt_1$. This action 
 Secondly, the action of the Grothendieck-Teichm\"uller group on the set of Drinfeld associators is free and transitive. Hence any Drinfeld associator can be obtained from a particular one, say $\Phi_{AT}$, by integrating the flow on the space of Drinfeld associators generated by some $\grt_1$-element.
By this action $\grt_1$ also acts on $L_\infty$ formality morphisms obtained by the Tamarkin construction as detailed in the preceding subsection.
As seen there this action is up to homotopy the same as the one obtained by first mapping $\grt_1$ to $H^0(\GC_2)$ by the map of Theorem \ref{thm:main}, and then acting with a representing degree zero graph cohomology class as in section \ref{sec:grcactionfm}.
%$\psi$ and $\Gamma$ on $\Tpoly$ agreed, up to homotopy. 
From these arguments, one sees that any Tamarkin type formality morphism can be obtained, up to homotopy, by integrating the action of some degree 0 cocycle $\Gamma\in \GC_2$ on $\Tpoly$, followed by applying Kontsevich's formality morphism $\mU_K$. 
Since $\mU_K$ is already constructed by ``graphical means'', one can incorporate this action into the construction of $\mU_K$. One just has to modify the way weights are associated to graphs in M. Kontsevich's definition of $\mU_K$.
This shows how one can obtain, using ``purely graphical techniques'', a big class of formality morphisms.% More details and a discussion of a related construction of Drinfeld associators using configuration space integrals will be contained in a forthcoming paper.

\subsection{An answer to a question of B. Tsygan}
Let $\Phi$ be a Drinfeld associator.
To $\Phi$ one can assign a formal odd power series 
\[
 \tilde f_\Phi(x) = \sum_{j\geq 1} \tilde f_{2j+1} x^{2j+1}
\]
where the number $\tilde f_{2j+1}$ is the coefficient of $X^{2j}Y$ in the series $\Phi(X,Y)$.
Now let $\psi$ be a $\grt_1$-element, and use it to act on $\Phi$. It is easy to check that this action changes the coefficient of $X^{2j}Y$ by the coefficient of $\ad_X^{2j}Y$ in $\psi(X,Y)$. Let us call this coefficient $\tilde s_{2j+1}(\psi)$. The $\tilde s_{2j+1}$ ($j=1,2,\dots$) can be checked to vanish on $\co{\grt_1}{\grt_1}$ and hence form Lie algebra cocycles of $\grt_1$. 

Next let $\mU_\Phi$ be the formality morphism associated to $\Phi$ by the Tamarkin construction as in section \ref{sec:tamconstr}. One can use it to construct a proof of Duflo's Theorem along the lines of \cite[section 8]{K1}. The version of the Duflo morphism obtained then has the form
\[
 J\circ e^{\sum_{j\geq 1} f_{2j+1}
 \tr(\ad_\p^{2j+1}) }
\]
where $J$ is the usual Duflo morphism and the $f_{2j+1}$ are numbers, depending on the formality morphism. Hence there is a natural way to define another odd formal power series for $\Phi$, namely
\[
 f_\Phi = \sum_{j\geq 1} {f}_{2j+1}  x^{2j+1},
\]
B. Tsygan asked the following question: Is $f_\Phi=\tilde f_\Phi$?
\begin{lemma}
 The answer to B. Tsygan's question is yes, i.~e., $f_\Phi=\tilde f_\Phi$.
\end{lemma}
Since the action of $\grt_1$ on Drinfeld associators is transitive, it will be sufficient to prove the following two statements.
%Any associator $\Phi$ can be obtained by acting with some element of the Grothendieck Teichm\"uller group on $\Phi_{AT}$. Hence, to prove the Lemma, it will be sufficient to show the following two statements.
\begin{enumerate}
 \item For the Alekseev-Torossian associator $\Phi_{AT}$, $f_{\Phi_{AT}}=\tilde{f}_{\Phi_{AT}}$.
 \item For $\psi\in \grt_1$ and any associator $\Phi$,
%\footnote{For this equation we extend the definition of $f_{(\cdot)}$ to arbitrary power series in $x,y$ in the obvious way.}
\[
f_{\psi\cdot\Phi}=\tilde{f}_{\psi\cdot\Phi}
\]
where $f_{\psi\cdot\Phi}$ should be understood as ``derivative of $f_{(\cdot)}$ along $\psi\cdot$ at $\Phi$''.
\end{enumerate}

Let us begin with the first statement. The Alekseev-Torossian associator can be shown to be even, i.~e., it contains only monomials in $X,Y$ of even degree. Hence $f_{\Phi_{AT}}=0$. From the previous subsection on the other hand we know that $\mU_{\Phi_{AT}}$ is homotopic to Kontsevich's morphism $\mU_K$. It is been shown in \cite{K1} and \cite{shoiwheels} that from $\mU_K$ one obtains the original Duflo morphism. Hence we conclude $\tilde f_{\Phi_{AT}}=0=f_{\Phi_{AT}}$.

For the second statement, fix some $\psi\in \grt_1$. The right hand side of the equation we want to show is $\tilde f_{\psi\cdot\Phi}=\sum_j \tilde s_{2j+1}(\psi)x^{2j+1}$. Let $\Gamma$ be graph cocycle, whose cohomology class corresponds to $\psi$. From the identification of Lie algebras $\grt_1\cong H^0(\GC_2)$ 
one (of course) also obtains an identification of the Lie algebra cohomology classes. The class represented by $\tilde s_{2j+1}$ corresponds to the graph homology class represented by a wheel with $2j+1$ spokes by Proposition \ref{prop:sigmawheel}. If we call the graph cycle given by this wheel graph $s_{2j+1}\in \GC_2^*$, with sign convention as in Proposition \ref{prop:sigmawheel}, we can hence write
\[
 \tilde f_{\psi\cdot\Phi} = -\sum_j s_{2j+1}(\Gamma)x^{2j+1}.
\]
Next we should compute ${f}_{\psi\cdot\Phi}$. Suppose we change our formality morphism by precomposing it with some $L_\infty$-automorphism $\mV$ of $\Tpoly$. Then the associated Duflo morphism changes by precomposing with the automorphism $\mV_1^\pi$ of $S\alg{g}$. Here $\alg{g}$ is the Lie algebra for which we write down the Duflo morphism and $\mV_1^\pi$ is the first Taylor component of the $L_\infty$-morphism obtained by twisting $\mV$ by the Poisson structure $\pi$ on $\alg{g}^*$.
Suppose now that that $\mV$ is obtained by integrating the action of the graph cocycle $-\Gamma$, i.~e., $\mV = \exp(-t\Gamma\cdot)$.\footnote{Note that the degree zero graph cocyles carry a left action on $\Tpoly$, hence to obtain a left action on formality morphisms one needs to introduce an extra sign.} We can assume without loss of generality that $\Gamma$ is 1-vertex irreducible. Then using the linearity of $\pi$ it can be checked that the only terms in $\mV_1^\pi$ that do not vanish on $S\alg{g}$ come from the wheel graphs. More precisely
\[
 \mV_1^\pi\mid_{S\alg{g}} = e^{-t \sum_j s_{2j+1}(\Gamma) \tr(\ad_\p^{2j+1})}.
\]
Taking the derivative at $t=0$ we can hence conclude that ${f}_{\psi\cdot\Phi}= -\sum_j s_{2j+1}(\Gamma)X^{2j+1}=\tilde f_{\psi\cdot\Phi}$.
\hfill \qed

\subsection{Globalization and a proof of a result announced by M. Kontsevich}
 Degree zero cocycles of $\GC_2$ can be naturally interpreted as $L_\infty$-derivations of the polyvector fields on $\R^d$.
 Using some globalization methods as in \cite{dolgushev, dolgushev-2005, calaquerossibook}, one can obtain $L_\infty$-derivations also on the space of polyvector fields on any smooth manifold, or on the sheaf cohomology of the sheaf of holomorphic polyvector fields on a complex manifold. This has been used in \cite{vasilychrisme} together with Proposition \ref{prop:sigmawheel} above, to show that the cocycles corresponding to the $\sigma_{2j+1}$ act on the sheaf cohomology of the sheaf of holomorphic polyvector fields by contraction with the odd Chern characters, thus showing an earlier claim by M. Kontsevich \cite{K2}.

\appendix
%%\input appgrcohom.tex
%\subsection{Harrison complex of the cofree coalgebra}
\section{Harrison complex of the cofree coalgebra}
\label{sec:harrison}
Let us recall some general (and well known) facts from homological algebra. Let 
\[
 \F_{\coCom}(X_1,\dots, X_k)
\]
be the free cocommutative coalgebra cogenerated by symbols $X_1,\dots, X_k$. Since this coalgebra is cofree, its reduced Harrison complex $\Harr(\F_{\coCom}(X_1,\dots, X_k))$ has cohomology
\[
H^j(\Harr(\F_{\coCom}(X_1,\dots, X_k)))
\cong
\begin{cases}
\GF X_1\oplus \dots\oplus \GF X_k & \text{for $j=1$} \\
0 & \text{otherwise}
\end{cases}.
\]
Here $\GF X_1\oplus \dots\oplus \GF X_k\cong \GF^k$ shall denote the $k$-dimensional vector space generated by the set $X_1,\dots, X_k$. Since the differential on $\Harr(\cdots)$ cannot create or annihilate any of the formal variables
\[
\Harr(\F_{\coCom}(X_1,\dots, X_k)) 
\]
inherits a $\Z^k$ grading. 
Of particular interest to us is the subcomplex of degree $(1,1,\dots, 1)$, 
\[
\Harr^{(1,\dots, 1)}(\F_{\coCom}(X_1,\dots, X_k)).
\]
It immediately follows from the formula for the cohomology of  $\Harr(\cdots)$ above that
\[
H^j(\Harr^{(1,\dots, 1)}(\F_{\coCom}(X_1,\dots, X_k)))
\cong
\begin{cases}
\GF & \text{for $j=1$ and $k=1$} \\
0 & \text{otherwise}
\end{cases}.
\]
Let us denote the $p$-fold symmetric tensor product of a dg vector space $V$ by $S^p V$, i.e., 
\[
S^pV=(V^{\otimes p})^{S_p}.
\]
$S^p \Harr(\F_{\coCom}(X_1,\dots, X_k))$ inherits a $\Z^k$ grading. 
The subcomplex of degree $(1,\dots, 1)$ is of special importance to us and we will abbreviate
\[
V_{p,k,n} := (S^p(\Harr(\F_{\coCom}(X_1,\dots, X_k))[1-n]))^{(1,\dots, 1)}.
\]
Because taking invariants with respect to finite group actions commutes with taking cohomology,
\[
H^j(V_{p,k,n})
\cong
\begin{cases}
\GF & \text{for $j=kn$ and $p=k$} \\
0 & \text{otherwise}
\end{cases}.
\]
Note that $V_{p,k,n}$ carries a natural action of the symmetric group $S_k$ by permuting the indices of the variables $X_j$. Again because taking invariants with respect to finite group actions commutes with taking cohomology, the following Lemma is evident.
\begin{lemma}
\label{lem:harrisoncohom}
Let $p, k,n\in \mathbb{N}$, let $G\subset S_k$ be a subgroup, and let $M$ be some $G$-module. Then
\[
H^j((V_{p,k,n}\otimes M)^{G})
\cong 
\begin{cases}
(M\otimes \sgn^{\otimes n})^{G} & \text{for $j=kn$ and $p=k$} \\
0 & \text{otherwise}
\end{cases}.
\]
Here $V_{p,k,n}$ is considered a $G$-module, and the $G$-action on the tensor product is the diagonal action. 
\end{lemma}
 
\section{The deformation complex of \texorpdfstring{$n$}{n}-algebras}
\label{sec:endef}
\subsection{The (co)operads \texorpdfstring{$e_n$}{en} and \texorpdfstring{$e_n^\vee$}{env}}
The operad $e_n$ is the operad governing $n$-algebras. An $n$-algebra is a (graded) vector space $V$ with binary operations $\wedge$ of degree 0 and $\co{\cdot}{\cdot}$ of degree $1-n$ satisfying the following relations:
\begin{enumerate}
 \item $(V,\wedge)$ is a graded commutative algebra.
 \item $(V[n-1],\co{\cdot}{\cdot})$ is a graded Lie algebra.
 \item For all (homogeneous) $v\in V$, the unary operation $\co{v}{\cdot}$ is a derivation of degree $|v|+1-n$ on $(V,\wedge)$.
\end{enumerate}

Elements of $e_n(N)$ can be written as linear combinations of expressions of the form 
\begin{equation}
\label{equ:enform}
 L_1(X_{1},\cdots,X_{N})\wedge \cdots \wedge L_k(X_{1},\cdots,X_{N})
\end{equation}
where $X_1,\dots,X_N$ are formal variables, $L_j$ are $\Lie\{n-1\}$ words and each $X_i$ occurs exactly once in the expression. The action of the symmetric group $S_N$ on $e_n(N)$ is given by permuting the labels on the $X_1,\dots,X_N$.
From the length of the individual Lie words one can derive various filtrations on $e_n$.  We denote by by 
\[
 e_n(N)_{k_1,k_2,\dots}
\]
the subspace spanned by elements of the form \eqref{equ:enform} with $k_j$ the number of Lie words of length $j$, $j=1,2\dots$. Clearly $\sum_j jk_j =N$ and furthermore the cohomological degree is $N-\sum_j k_j$.

\begin{ex}
For example, for the expression
\[
 X_1\wedge\co{X_4}{X_3}\wedge X_2
\]
the number of Lie words of length 1 is $k_1=2$.
\end{ex}

The cooperad $e_n^\vee$ is the Koszul dual cooperad to $e_n$. One can show that $e_n^\vee=e_n^*\{n\}$. Dualizing the direct sum decomposition of $e_n(N)$ into  subspaces $e_n(N)_{k_1,k_2,\dots}$ we may write 
\begin{equation}
\label{equ:enveesplitting}
 e_n^\vee(N) = \bigoplus_{k_1,k_2,\dots} e_n^\vee(N)_{k_1,k_2,\dots}
\end{equation}
where $e_n^\vee(N)_{k_1,k_2,\dots}$ is dual to $(e_n\{n\})(N)_{k_1,k_2,\dots}$.

\begin{rem}
 One may represent $e_n^\vee$ using graphs. Then $e_n^\vee(n)_{k_1,k_2,\dots}$ is the space of graphs with $k_1$ isolated vertices, $k_2$ connected components of size 2, $k_3$ connected components of size 3 etc.
\end{rem}

%We will be interested in the filtration coming from the number $k_1$ of Lie words of length 1. For example, for the expression
%\[
% X_1\wedge\co{X_4}{X_3}\wedge X_2
%\]
%the number of Lie words of length 1 is $k_1=2$. 
%
%
%
% Hence from the filtration by $k_1$ on $e_n$ one obtains a filtration on $e_n^\vee$.

%Elements of $e_n^\vee(N)$ can be written as linear comvbinations of expressions of the form 

\subsection{\texorpdfstring{$\hoe_n$}{hoen} and a(nother) filtration on the deformation complex}
Recall that $\hoe_n=\Omega(e_n^\vee)$ is the operadic cobar construction of the cooperad $e_n^\vee$. 
Inserting the direct sum decomposition above into \eqref{equ:convalgebra} we see that 
\[
\Def(\hoe_n\to e_n) = \prod_{k_1,k_2,\dots} 
\Hom_{S_N}(e_n^\vee(N)_{k_1,k_2,\dots}, e_n(N))
\]
where we abbreviate $N=\sum_j jk_j$ within the product on the right hand side.
The differential on $\Def(\hoe_n\to e_n)$ contains one part, say $d_+$, that raises $k_1$ by one, and one part that leaves $k_1$ constant. Let us put a filtration on  $\Def(\hoe_n\to e_n)$ so that $d_+$ is the differential on the associated graded.
Concretely, 
\[
\mF^p \Def(\hoe_n\to e_n) = 
\prod_{k_1,k_2,\dots } 
\Hom_{S_N}(e_n^\vee(N)_{k_1,k_2,\dots }, e_n(N))_{ \geq k_1-p }
\]
where the final subscript shall indicate that only the subspace of cohomological degrees $\geq k_1-p$ is taken.

\subsection{A more concrete description of the differential \texorpdfstring{$d_+$}{d+}}
Let 
\[
e_n^\vee(N)' := \bigoplus_{k_2, k_3,\dots} e_n^\vee(N)_{0, k_2, k_3,\dots}.
\]
Then we may write
\begin{align*}
\gr \Def(\hoe_n\to e_n)
&\cong 
 \prod_N \prod_{k_1\leq N} 
\Hom_{S_{k_1}\times S_{N-k_1}}( \sgn^n_{k_1}\otimes e_n^\vee(N-k_1)',e_n(N))
\\ &\cong
\prod_N \prod_{k_1\leq N}
\Hom_{S_{N-k_1}}( e_n^\vee(N-k_1)',\Hom_{S_{k_1}}(\sgn^n_{k_1}, e_n(N))).
\end{align*}
The differential $d_+$ on this complex may be seen as induced by a differential on
\[
\Hom_{S_{k_1}}(\sgn^n_{k_1}, e_n(N) ) \cong \left( \sgn^n_{k_1}\otimes e_n(N) \right)^{S_{k_1}}
\]
which we abusively also denote by $d_+$.
%Up to permutations of labels, any element of $e_n(N)$ can be written as a sum of expressions of the form
%\[
% X_1\wedge X_2 \wedge \cdot \wedge X_{k_1} \wedge L_1(X_{1},\cdots,X_{N})\wedge \cdots \wedge L_k(X_{1},\cdots,X_{N})
%\]
%where all $\Lie\{n-1\}$-words $L_1,\dots,L_n$ have length at least 2.
%Hence one can decompose
%\[
% e_n(N) = \bigoplus_{k_1=0}^N \Ind_{S_{k_1}\times S_{N-k_1}}^{S_N} (\sgn^n_{k_1}\otimes e_n(N-k_1)')
%\]
%where $e_n(N-k_1)'\subset e_n(N-k_1)'$ is the subspace generated by expressions as above with all Lie words of length at least 2.
%A similar splitting holds for the $e_n^\vee$, and allows us to write
%\begin{multline*}
%  \Hom_{S_N}(e_n^\vee(N), e_n)_{k_1} \cong 
%\Hom_{S_{k_1}\times S_{N-k_1}}( \sgn^n_{k_1}\otimes e_n^\vee(N-k_1)',e_n(N))
%\cong \\ \cong
%\Hom_{S_{N-k_1}}( e_n^\vee(N-k_1)',\Hom_{S_{k_1}}(\sgn^n_{k_1}, e_n(N))).
%\end{multline*}
%where $\Hom_{S_N}(e_n^\vee(N), e_n)_{k_1}\subset \Hom_{S_N}(e_n^\vee(N), e_n)$ is the subspace of fixed $k_1$. 
%The differential $d$ on $\Der(\hoe_n \to e_n)$ has some part, say $d_+$, that increases $k_1$ by 1 and some part that leaves $k_1$ constant or decreases it. We will need an explicit formula for $d_+$. 
Elements of $\left( \sgn^n_{k_1}\otimes e_n(N) \right)^{S_{k_1}}$
may be understood as expressions $P(X_1,\dots,X_{k_1}, A_1,\dots,A_m)$ in formal variables $X_1,\dots, X_{k_1}, A_1,\dots,A_m$, with $m=N-k_1$, of the form \eqref{equ:enform}, invariant under permutations (with signs) of $X_1,\dots, X_{k_1}$. 
%is spanned by expressions $P(X_1,\dots,X_{k_1}, A_1,\dots,A_m)$, where $m=N-k_1$ and $P$ is a linear combination of expressions as in \eqref{equ:enform} that is invariant under permutations (with signs) of $X_1,\dots, X_{k_1}$. 
Using this notation, the formula for the differential $d_+$ is:
\begin{align*}
&\pm (d_+ P)(X_1,\dots,X_{k_1+1},A_1,\dots A_m) = \\
&= \sum_{i=1}^{k_1+1}(-1)^{n (i+1)} \co{X_i}{P(X_1,\dots,\hat{X}_i,\dots,X_{k_1+1},A_1,\dots ,A_m)}
\\ &\quad\quad -
\sum_{1\leq i<j\leq k_1+1}(\pm)(-1)^{n (i+j+1)} P(\co{X_i}{X_j},X_1,\dots,\hat{X}_i,\dots,\hat{X}_j,\dots,X_{k_1+1},A_1,\dots ,A_m)
\\ &\quad\quad -
\sum_{i=1}^{k_1+1}\sum_{j=1}^m (\pm) (-1)^{n (i+1)} P(X_1,\dots,\hat{X}_i,\dots,X_{k_1+1},A_1,\dots,\co{X_i}{A_j}, \dots,A_m)
\\ &=
\sum_{1\leq i<j\leq k_1+1}(\pm)(-1)^{n (i+j+1)} P(\co{X_i}{X_j},X_1,\dots,\hat{X}_i,\dots,\hat{X}_j,\dots,X_{k_1+1},A_1,\dots ,A_m).
\end{align*}
Here the signs $(\pm)$ occur only in the case of even $n$ and are a bit tricky, owed to the oddness of the Lie bracket.\footnote{Strictly speaking, our notation here is almost wrong since different summands of $P$ may pick up different signs.} For example, in the third line, one should write down a term of $P$, and then compute the number of brackets to the left of the first argument. This gives the sign.
The last equality follows since $\co{X_i}{\cdot}$ is a derivation with respect to the product and the Lie bracket. 

\begin{rem}
\label{rem:dpluszero}
 In particular, note that $d_+=0$ on the subspace with $k_1=0$.
\end{rem}

%Let us describe in more detail the $S_{N-k_1}$-module $\Hom_{S_{k_1}}(\sgn^n_{k_1}, e_n(N))$. We claim that
%\[
% \Prod_{k_1\geq 0} \Hom_{S_{k_1}}(\sgn^n_{k_1}, e_n(m+k_1))
%=
%\]

\subsection{The cohomology is concentrated in degree \texorpdfstring{$k_1=0$}{k1=0}}
Let 
\[
\Xi= \prod_{k_2,k_3, \dots} 
\Hom_{S_N}(e_n^\vee(N)_{0,k_2,k_3,\dots}, e_n(N)) 
\subset
\Der(\hoe_n \to e_n)
\]
be the subspace of $\Der(\hoe_n \to e_n)$ corresponding to $k_1=0$.
The main result of this section is the following.
\begin{prop}
\label{prop:Xiqiso}
 $\Xi \subset \Der(\hoe_n \to e_n)$ is a subcomplex. The inclusion is a quasi-isomorphism.
\end{prop}

The first statement follows directly from Remark \ref{rem:dpluszero}. For the second statement, let us compute the spectral sequence associated to the filtration by $k_1$ introduced above. Because $\Xi$ is a subcomplex, it will be sufficient to show that the cohomology of the associated graded $\gr \Def(\hoe_n\to e_n)$ is $\Xi$. Given the description of the differential on the associated graded above, it is sufficient to show that the complexes
\[
 C_m := \left( \bigoplus_{k_1\geq 1} \left( \sgn^n_{k_1}\otimes e_n(m+k_1) \right)^{S_{k_1}}, d_+ \right)
\]
are acyclic. Note that clearly $\Lie\{ n-1 \} (m+k_1)\subset e_n(m+k_1)$. One can make a first reduction on the problem.
\begin{lemma}
 If the complexes 
\[
 L_m := \left( \bigoplus_{k_1\geq 1} \left( \sgn^n_{k_1}\otimes \Lie\{ n-1 \}(m+k_1) \right)^{S_{k_1}}, d_+ \right)
\]
are acyclic, then so are the complexes $C_m$.
\end{lemma}
\begin{proof}
Consider a complex $\tilde C_m$ that is (like $C_m$) spanned by products 
\[
P(X_1,\dots,X_{k_1}, A_1,\dots,A_m) 
\]
of Lie words, but with the $A_1,\dots,A_m$ allowed to occur with repetitions. It splits as a sum of subcomplexes according to the number of occurences of $A_1, A_2$ etc. $C_m\subset\tilde C_m$ is simply the subcomplex in which each $A_1,A_2,\dots$ occurs exactly once. It is sufficient to prove that the cohomology of $\tilde C_m$ is precisely its $k_1=0$-part.
However, $\tilde C_m$ is a cofree cocommutative coalgebra, the coproduct being the deconcatenation of Lie words. The cogenerators are those expressions $P(X_1,\dots,X_{k_1}, A_1,\dots,A_m)$ containing only a single Lie word. If $L_m$ is acyclic, the cohomology of the space of generators is precisely its $k_1=0$-part.
%\[
% C_m = \left( \bigoplus_{k_1\geq 1} \left( \sgn^n_{k_1}\otimes e_n(m+k_1) \right)^{S_{k_1}} \right)
%\]

 %Recall from the previous subsection that $C_m$ is spanned by wedge products of Lie words $P(X_1,\dots,X_{k_1}, A_1,\dots,A_m)$. In particular, 
\end{proof}

Hence we are left with showing that the complexes $L_m$ are acylic. For $L_0$ this is easy to see ($L_0$ is only two-dimensional by the Jacobi identity). Now suppose that $m\geq 1$. Then the space of Lie words in $X_1,\dots,X_{k_1}, A_1,\dots,A_m$, each symbol occuring exactly once, can be identified with the space of ordinary (associative) words in symbols $X_1,\dots,X_{k_1}, A_2,\dots,A_m$, each occuring exactly once. For example (for $m=3$, $k_1=2$)
\[
 \co{X_1}{\co{A_3}{\co{X_2}{\co{A_2}{A_1}}}}
\leftrightarrow
X_1A_3X_2A_2.
\]
A basis for the (anti-)symmetric part in the $X_j$ can be given by symbols like
\[
 XA_3XA_2.\leftrightarrow X_1A_3X_2A_2 \pm X_2A_3X_1A_2.
\]
In this basis, the differential is given by ``doubling'' $X$'s, i.e., 
\[
 d_+(XA_3XA_2) = XXA_3XA_2 \pm XA_3XXA_2.
\]
It is well known that this complex is acyclic. An explicit homotopy is ``contracting'' pairs of $X$'s. Hence the proposition is proven. \hfil \qed

%Here the differential $d_+$ is the one from the previous subsection. One can see that the above complex splits into subcomplexes, according to the 
%\begin{align*}
%\end{align*}

\section{\texorpdfstring{$\Def(\hoe_n\to\Graphs_n)$ and $\Def(\hoe_n\to\Gra_n)$}{Deformation complexes}}
\label{sec:defengra}
At several points of this paper we use combinatorial arguments to compute the cohomology of $\Def(\hoe_n\to\Graphs_n^\whl)$ or some subcomplex or quotient. The combinatorics are much easier to understand and to describe in words if one identifies elements of $\Def(\hoe_n\to\Graphs_n^\whl)$ with certain linear combinations of graphs. 
Using the splitting \eqref{equ:enveesplitting} again, we may write
\[
\Def(\hoe_n\to \Graphs_n^\whl) = \prod_{k_1,k_2,\dots} 
\Hom_{S_N}(e_n^\vee(N)_{k_1,k_2,\dots}, \Graphs_n^\whl(N))
\]
where again $N=\sum_j jk_j$ within the product. 
%In this section we give a graphical description of the complexes $\Def(\hoe_n\to\Graphs_n)$ and its quotient $\Def(\hoe_n\to\Gra_n)$. Let us start with the former.
%We have 
%\[
%\Def(\hoe_n\to\Graphs_n) \cong \prod_N 
% \Hom_{S_N}(e_n^\vee(N), \Graphs_n(N)).
%\]
%As in Appendix \ref{sec:endef}, this space is a graded Lie algebra and the differential is given by the bracket with a Maurer-Cartan element. Here the Maurer-Cartan element is the image $\tilde{f}=\tilde f_{\wedge}+\tilde f_{\bracket}$ under $\Def(\hoe_n\to e_n)\to \Def(\hoe_n\to\Graphs_n)$ of the map $f=f_{\wedge}+f_{\bracket}$ of Appendix \ref{sec:endef}. More precisely, $\tilde{f}$ maps the coproduct in $e_n^\vee(2)$ to the graph in $\Graphs_n(2)$ with two vertices and one edge connecting them, and maps the cobracket to the ``empty'' graph with two vertices and no edge.

%For later proofs, it is helpful to have a combinatorial description of $\Def(\hoe_n\to\Graphs_n)$ as a complex spanned by special graphs. Elements of $e_n^\vee(N)$ can be written as linear combinations of expressions like
%\[
%(324)\wedge (51) \wedge (6)
%\] 
%(for $n=6$), where the order of the brackets is arbitrary (up to sign) and the order of numbers within the brackets is defined only modulo shuffles. 
%One can decompose any element of $e_n^\vee(N)$ according to the number of brackets with $j$ numbers inside. Let 
%\[
% e_n^\vee(N)_{k_1,k_2,\dots}
%\]
%Be the space of elements of $e_n^\vee(N)$ with $k_1$ brackets of length 1, $k_2$ brackets of length 2 etc. Here $N=\sum_j j k_j$.

The space $\Hom_{S_N}(e_n^\vee(N)_{k_1,k_2,\dots}, \Graphs_n^\whl(N))$ 
may now (-up to a degree shift-) be understood as the subspace of elements $\Gamma \in \Graphs_n^\whl(N)$ with the following symmetry properties:
\begin{enumerate}
 \item Consider the $N$ external vertices to be organized into clusters, with $k_1$ clusters of $1$ vertices, followed by $k_2$ clusters of 2 vertices, etc. Then the (linear combination of) graphs $\Gamma$ must be invariant under interchange (with sign) of clusters of the same size. The sign is $(-1)^{j+n+1}$ for the interchange of clusters of size $j$.
 \item $\Gamma$ must ``vanish on shuffles'' in any cluster. This means the following. Fix some cluster of length $j$, and fix $j_1,j_2\geq 1$ s.t. $j=j_1+j_2$. Let $\mathit{USh}_{j_1,j_2}$ be the set of $(j_1,j_2)$-unshuffle permutations, which we consider acting on $\Graphs_n^\whl(N)$ by permuting the vertices in the cluster under consideration. Then we require
\[
 \sum_{\sigma \in \mathit{USh}_{j_1,j_2}} \pm \sigma\cdot \Gamma = 0.
\]
\end{enumerate}

%\begin{lemma}
% The space $\Hom_{S_N}(e_n^\vee(N)_{k_1,k_2,\dots}, \Graphs_n(N))$ is isomorphic to the space of linear combinations of graphs in $\Graphs_n(N)$ that have the following symmetry properties. 
%\begin{enumerate}
% \item Consider the input vertices as organized into clusters, with $k_1$ clusters of $1$ vertices, $k_2$ clusters of 2 vertices, etc. Then the (linear combination of) graphs $\Gamma\in\Graphs_n(N)$ must be invariant under interchange (with sign) of clusters of the same size. The sign is $(-1)^{j+n+1}$ for the interchange of clusters of size $j$.
% \item We fix a linear ordering on the vertices of each cluster. Then $\Gamma$ must ``vanish on shuffles'' in any cluster. This means the following. Fix some cluster of length $j$, and $j_1,j_2\geq 1$ s.t. $j=j_1+j_2$. Let $\mathit{USh}_{j_1,j_2}$ be the set of $(j_1,j_2)$-unshuffle permutations, which we consider acting on $\Graphs(n)$ by permuting the vertices in the cluster under consideration. Then we require
%\[
% \sum_{\sigma \in \mathit{USh}_{j_1,j_2}} \pm \sigma\cdot \Gamma = 0.
%\]
%\end{enumerate}
%\end{lemma}

The complex $\Def(\hoe_n\to\Gra_n^\whl)$ has a similar interpretation in terms of graphs. Indeed, it is obtained from $\Def(\hoe_n\to\Graphs_n^\whl)$ by sending all graphs with internal vertices to zero. The complexes $\Def(\hoe_n\to\Gra_n)$ and $\Def(\hoe_n\to\Graphs_n)$ can similarly be given a graphical description, by just dis-allowing tadpoles in the graphs.

\begin{rem}
 The degree of a graph can be computed as $n\cdot(\#vertices-1)-(n-1)\cdot(\#edges)$, when we count a cluster with $k$ vertices as $k$ vertices (of course) and $k-1$ edges.
\end{rem}

\subsection{A Graphical description of the differential}
\label{sec:derengradiff}
\begin{figure}
\centering
\[
\begin{tikzpicture}[scale=1,
 vert/.style={draw,outer sep=0,inner sep=0,minimum size=5,shape=circle,fill},
 helper/.style={outer sep=0,inner sep=0,minimum size=0,shape=coordinate},
 default_edge/.style={draw},
 ext/.style={draw, outer sep = 0, inner sep=0,minimum size=5,shape=circle}]
 
 \node (v0) at (0.300000011920929,3.09999996423721) [vert] {};
 \node (v1) at (1.80000007152557,3.30000001192093) [vert] {};
 \node (v2) at (1.80000007152557,2.89999997615814) [vert] {};
 \node (v3) at (4.09999990463257,3.30000001192093) [vert] {};
 \node (v4) at (4.09999990463257,2.89999997615814) [ext] {};
 \node (v5) at (2.60000014305115,3.09999996423721) [ext] {};
 \node (v6) at (3.89999985694885,3.70000001043081) [helper] {};
 \node (v7) at (4.09999990463257,3.70000001043081) [helper] {};
 \node (v8) at (4.29999995231628,3.70000001043081) [helper] {};
 \node (v9) at (3.89999985694885,2.60000002384186) [helper] {};
 \node (v10) at (4.09999990463257,2.60000002384186) [helper] {};
 \node (v11) at (4.29999995231628,2.60000002384186) [helper] {};
 \node (v12) at (2.40000009536743,2.79999995231628) [helper] {};
 \node (v13) at (2.60000014305115,2.79999995231628) [helper] {};
 \node (v14) at (2.79999995231628,2.79999995231628) [helper] {};
 \node (v15) at (2.40000009536743,3.39999997615814) [helper] {};
 \node (v16) at (2.60000014305115,3.39999997615814) [helper] {};
 \node (v17) at (2.79999995231628,3.39999997615814) [helper] {};
 \node (v18) at (1.6000000834465,2.60000002384186) [helper] {};
 \node (v19) at (1.80000007152557,2.60000002384186) [helper] {};
 \node (v20) at (2.00000011920929,2.60000002384186) [helper] {};
 \node (v21) at (1.6000000834465,3.6000000089407) [helper] {};
 \node (v22) at (1.80000007152557,3.6000000089407) [helper] {};
 \node (v23) at (2.00000011920929,3.6000000089407) [helper] {};
 \node (v24) at (0.100000001490116,3.39999997615814) [helper] {};
 \node (v25) at (0.300000011920929,3.39999997615814) [helper] {};
 \node (v26) at (0.5,3.39999997615814) [helper] {};
 \node (v27) at (0.100000001490116,2.79999995231628) [helper] {};
 \node (v28) at (0.300000011920929,2.79999995231628) [helper] {};
 \node (v29) at (0.5,2.79999995231628) [helper] {};
 \node (v58) at (0.400000005960464,3.09999996423721)
 [helper,label=0:$\mapsto \sum\pm$] {};
 \node (v57) at (2.70000004768372,3.09999996423721)
 [helper,label=0:$\mapsto \sum\pm$] {};
 \node (v59) at (0.29999977350235,1.89999985694885) [ext,label=180:(] {};
 \node (v60) at (0.699999749660492,1.89999985694885) [ext] {};
 \node (v61) at (0.899999797344208,1.89999985694885) [ext] {};
 \node (v62) at (1.09999978542328,1.89999985694885) [ext] {};
 \node (v63) at (1.29999983310699,1.89999985694885) [ext,label=0:)] {};
 \node (v64) at (0.499999761581421,1.89999985694885) [ext] {};
 \node (v65) at (3.10000014305115,1.59999990463257) [ext,label=180:(] {};
 \node (v66) at (3.29999995231628,1.59999990463257) [ext] {};
 \node (v67) at (3.49999976158142,1.59999990463257) [ext,label=0:)] {};
 \node (v68) at (3.5,2) [ext,label=180:(] {};
 \node (v69) at (3.70000004768372,2) [ext] {};
 \node (v70) at (4.09999990463257,2) [ext,label=0:)] {};
 \node (v71) at (3.70000004768372,1.99999976158142) [ext] {};
 \node (v72) at (3.89999985694885,2) [ext] {};
 \node (v73) at (1.60000002384186,1.89999985694885)
 [helper,label=0:$\mapsto \sum\pm$] {};
 \node (v87) at (0.500000059604645,0.899999856948853) [ext,label=180:(] {};
 \node (v88) at (0.900000035762787,0.899999856948853) [ext] {};
 \node (v90) at (1.09999978542328,0.899999856948853) [ext] {};
 \node (v91) at (1.29999983310699,0.899999856948853) [ext,label=0:)] {};
 \node (v92) at (0.700000047683716,0.899999856948853) [ext] {};
 \node (v93) at (1.60000002384186,0.899999856948853)
 [helper,label=0:$\mapsto \sum\pm$] {};
 \node (v94) at (3.10000002384186,0.899999856948853) [ext,label=180:(] {};
 \node (v95) at (3.5,0.899999856948853) [ext] {};
 \node (v96) at (3.70000004768372,0.899999856948853) [ext] {};
 \node (v97) at (3.90000003576279,0.899999856948853) [ext] {};
 \node (v98) at (4.1000000834465,0.899999856948853) [ext,label=0:)] {};
 \node (v99) at (3.30000001192093,0.899999856948853) [ext] {};
 
 \draw[default_edge] (v3)--(v4);
 \draw[default_edge] (v4)--(v9);
 \draw[default_edge] (v10)--(v4);
 \draw[default_edge] (v24)--(v0);
 \draw[default_edge] (v0)--(v25);
 \draw[default_edge] (v26)--(v0);
 \draw[default_edge] (v3)--(v6);
 \draw[default_edge] (v7)--(v3);
 \draw[default_edge] (v3)--(v8);
 \draw[default_edge] (v4)--(v11);
 \draw[default_edge] (v5)--(v12);
 \draw[default_edge] (v13)--(v5);
 \draw[default_edge] (v5)--(v14);
 \draw[default_edge] (v5)--(v17);
 \draw[default_edge] (v16)--(v5);
 \draw[default_edge] (v5)--(v15);
 \draw[default_edge] (v23)--(v1);
 \draw[default_edge] (v1)--(v22);
 \draw[default_edge] (v21)--(v1);
 \draw[default_edge] (v1)--(v2);
 \draw[default_edge] (v2)--(v18);
 \draw[default_edge] (v19)--(v2);
 \draw[default_edge] (v2)--(v20);
 \draw[default_edge] (v29)--(v0);
 \draw[default_edge] (v0)--(v28);
 \draw[default_edge] (v27)--(v0);
 \draw[default_edge] (v68)--(v67);
 \end{tikzpicture}
\]
 \caption{\label{fig:defdiff} A schematic illustration of the various parts of the differential on the complex $\Def(\hoe_n\to\Graphs_n)$. Using the splitting of the differential as in \eqref{equ:diffsplit2}, the first line corresponds to $\delta$, the second to $d_\bracket$ and the third to $d_\wedge$.
For a detailed description, see the text below.}
\end{figure}

Let us describe the differential on $\Def(\hoe_n\to\Graphs_n^\whl)$ combinatorially, using the graphical language from the last subsection.
The differential has four parts:
\begin{enumerate}
 \item The first part splits an internal vertex into two internal vertices. It comes from the differential on $\Graphs_n^\whl$.
 \item The second part splits an external vertex into an external and an internal vertex. It also comes from the differential on $\Graphs_n^\whl$.
 \item The third part, $d_L$, splits a cluster of external vertices into two clusters, by splitting one external vertex in that cluster. It creates an edge between the two vertices that the original vertex was split up into (see Figure \ref{fig:defdiff}). This part of the differential is denoted by $d_\bracket$ in \eqref{equ:dwedge}.
 \item The fourth part, denoted $d_\wedge$ in \eqref{equ:dwedge},  also creates an external vertex, but does not split the cluster and does not introduce a new edge. Hence it maps a cluster of length $j$ to one of length $j+1$. (For a picture, see again Figure \ref{fig:defdiff}.)
\end{enumerate}

In the complex $\Def(\hoe_n\to\Gra_n)$ the first two terms are absent.

\section{Sketch of proof of Lemma \ref{lem:prekazhdan}}
\label{sec:prekazhdanproof}
Our goal is to show that the double complex
\[
 0\to \Xi_\conn \to \Def(\hoe_n\to \Gra_n^\whl)_\conn \to \Def(\hoLie_n\to \Gra_n^\whl)_\conn\to 0
\]
is acyclic. Equivalently, we have to show that the inclusion
\[
 \Xi_\conn \to \Def'(\hoe_n\to (\Gra_n^\whl))_\conn
\]
is a quasi-isomorphism, where the $'$ on the right hand side indicates that one restricts to those derivations with vanishing $\hoLie_n$-part.

To go further, we need to use notation from Appendix \ref{sec:defengra}, in particular the graphical representation of $\Def(\hoe_n\to \Gra_n^\whl)$.
For each graph occurring in $\Def(\hoe_n\to \Gra_n^\whl)_\conn$ one can associate a number $k_1$, the number of ``clusters'' (notation as in Appendix \ref{sec:defengra}) of length 1 in that graph. 
This yields a filtration on the graded vector space $\Def(\hoe_n\to \Gra_n^\whl)_\conn$, which descends to a filtration on $\Def'(\hoe_n\to (\Gra_n^\whl))_\conn$. One takes an appropriate spectral sequence such that the first differential increases $k_1$ by exactly one. 
More concretely, this differential, say $d_+$, acts by spliting off a length one cluster from any vertex. Let us introduce new terminology. Let us call vertices in clusters of length 1 ``internal'' and all others external. Then a computation very much similar to the computation of $H(\Graphs_n^\whl)\cong e_n$ shows that the cohomology of $\Def'(\hoe_n\to (\Gra_n^\whl))_\conn$ is given by closed graphs without internal vertices. In other words, graphs with all clusters of length $\geq 2$, and which actually lie in $\Def(\hoe_n\to e_n)_\conn\subset \Def(\hoe_n\to \Gra_n^\whl)_\conn$. But that space is the space we called $\Xi_\conn$.
\hfill \qed

\section{(Re-)Derivation of Furusho's result}
\label{sec:furusho}
In a remarkable paper \cite{furusho}, H. Furusho showed that the hexagon equation \eqref{equ:hex} is a consequence of the pentagon equation \eqref{equ:pent}, if one requires that $\phi\in \FreeLie(X,Y)$ (as in those equations) does not contain the term $\co{X}{Y}$. 
Rephrasing this result in operadic language, it reads as follows.
\begin{prop}
\[
 H^1(\Def(\Ass_\infty\to t[1])) \cong   H^1(\Def(\Com_\infty\to t[1])) \oplus \GF[-1].
\] 
\end{prop}
Here the $\GF[-1]$ corresponds to the cohomology class represented by $\phi=\co{X}{Y}$.
\begin{proof}
We will actually show that $H^1(\Def(\Ass_\infty\to \CG^\whl[1])) \cong   H^1(\Def(\Com_\infty\to \CG^\whl[1])) \oplus \GF[-1]$.
Take a spectral sequence as in the proof of Lemma \ref{lem:fCqiso}. The difference to the situation there is that we now have to compute the Hochschild instead of the Harrison cohomology of a cofree cocommutative coalgebra. But it is well known that the Hochschild cohomology of such an algebra is the Koszul dual algebra, i.~e., an (anti-)commutative free algebra. Translated into graphical language, it means that the cohomology is given by graphs, whose external vertices are connected by exactly one edge each, and which are antisymmetric under interchange of external vertices.
The space of such graphs with $q$ external vertices forms a subcomplex $C_q\subset \Def(\Ass_\infty\to \CG^\whl[1])$. In the language of \cite{LambrechtsTurchin} this subcomplex computes the part of the cohomology of Hodge degree $q$.
Let us compute the degree 1 cohomology of $C_q$ for various $q$.
We already saw that $H(C_1)\cong H^1(\Def(\Com_\infty\to t[1]))\cong \grt_1$.

Next, suppose we are given a closed linear combination of graphs $x_q \in C_q$ of cohomological degree 1. 
We split the differential on $\Def(\Ass_\infty\to \CG^\whl[1])$ into a part $d_s$ creating an external vertex, and a part $\delta$ not creating one. In fact, $\delta$ is the part of the differential coming from the differential on $\CG^\whl$.
If $q\neq 3$ we can write $x_q=-\delta y_q$ for some $y_q$ since $H(\CG^\whl)\cong \alg t$ is concentrated in degree 0.
Hence our cohomology class is also represented by $d_s y_q=:x_{q+1}$. If $q\geq 4$ we can continue in this manner (i.~e., $x_{q+1}=-\delta y_{q+1}$, etc...) indefinitely by degree reasons and see that the cohomology class represented by $x_q$ is trivial.
Hence we have 
\[
 H^1(\Def(\Ass_\infty\to t[1])) \cong H^1(C_1)\oplus H^1(C_2) \oplus H^1(C_3).
\]
If $q=3$, then $x_q$ describes a non-trivial cohomology class only if it describes an element of $\alg t_3$ which can be recovered by projecting to the internal-trivalent-tree part of the graphs. However, since the external vertices are of valence 1, the only possible tree with all external vertices of valence 1 is the graph, say $T_3$, with one internal vertex, see Figure \ref{fig:graphwoneivertex}. This graph corresponds to $\co{t_{12}}{t_{23}}\in \alg{t}_3$. Hence $H^1(C_3)\cong \K$, the one cohomology class represented by $T_3$.

The statement of the proposition is thus reduced to showing that $H^1(C_2)=0$. So suppose $x_2\in \C_2$ is a degree 1 cocycle. We decompose $x_2=x_{2a}+x_{2b}$ where $x_{2a}$ contains those graphs for which both external vertices are connected to the same internal vertex and $x_{2b}$ the remainder. One can see that necessarily $x_{2a}=\delta y_{2a}+(\cdots)$ where $y_{2a}\in C_2$ and $(\cdots)$ is a linear combination of graphs with the two (univalent) external vertices connected to different internal vertices. Hence we may assume that in fact $x_{2a}=0$ from the start. In this case $x_2=x_{2b}=\delta y_2$ for a $y_2$ obtained by contracting one of the external edges. One checks that combinatorially $d_s y_2$ cannot contain any trivalent-tree-part: Such a tree must necessarily have (at least) two internal vertices connected to different external vertices. Since one external vertex of every graph in $y$ has valence one, it means that one of the internal vertices must have two edges connected to the same external vertex, which is a contradiction (or rather, the graph is 0). Hence $d_s y$ is $\delta$-exact again and hence by the same reasoning as in the $q>3$-case, the cohomology class of $x_q$ vanishes. Hence $H^1(C_2)=0$ and we are done.

\end{proof}

\begin{rem}
 Note that the calculation of $H^1(\Def(\Ass_\infty\to t[1]))$ above is similar to calculations of $H(\gr^p \Def(\hoe_2 \to e_2)_\conn)$ performed in section \ref{sec:derhoe2grt}, up to some sign and degree differences.
The Hodge degree $q$ corresponds to the parameter $p+1$ in section \ref{sec:derhoe2grt}.
\end{rem}

\begin{figure}
\centering
\begin{tikzpicture}[scale=1,
vert/.style={draw,outer sep=0,inner sep=0,minimum size=5,shape=circle,fill},
helper/.style={outer sep=0,inner sep=0,minimum size=5,shape=coordinate},
default edge/.style={draw},
ext/.style={draw,outer sep=0,inner sep=0,minimum size=5,shape=circle}]

\node (v0) at (0.300000011920929,1.19999998807907) [ext,label=-90:1] {};
\node (v1) at (0.900000035762787,1.19999998807907) [ext,label=-90:2] {};
\node (v2) at (1.5,1.19999998807907) [ext,label=-90:3] {};
\node (v3) at (0.900000035762787,1.79999999701977) [vert] {};

\draw[default edge] (v0)--(v3);
\draw[default edge] (v3)--(v1);
\draw[default edge] (v2)--(v3);
\end{tikzpicture}
 \caption{\label{fig:graphwoneivertex} The graph corresponding to the $\co{t_{12}}{t_{23}}\in \alg{t}_3$.}
\end{figure}

\section{The one vertex irreducible part of \texorpdfstring{$\GC_n$}{GCn} is quasi-isomorphic to \texorpdfstring{$\GC_n$}{GCn}}
\label{sec:onevi}
%Probably everything in this subsection has been shown by M. Kontsevich, but has not been written up.
Let $\GC_{n}^{1vi}\subset \GC_n$ be the subspace of 1-vertex irreducible graphs, i.e., those graphs that remain connected after deleting one vertex.
\begin{lemma}
 $\GC_{1vi}\subset \GC$ is a sub-dg Lie algebra.
\end{lemma}
\begin{proof}
 It is clear.
\end{proof}

The following Proposition has been shown by Conant, Gerlits and Vogtman.
\begin{prop}[\cite{conant_cut_2005}]
\label{prop:1vi}
 $\GC_n^{1vi}\hookrightarrow \GC$ is a quasi-isomorphism.
\end{prop}
We nevertheless give a different short sketch of proof for completeness, following the idea of Lambrechts and Volic \cite{LV}.
\begin{proof}[(very sketchy) Sketch of proof]
 We need to show that $\GC_n/\GC_n^{1vi}$ is acyclic. Any non-1vi graph can be written using the following data:
(i) a family of 1-vertex-irreducible graphs (``1vi components'') (ii) a tree (iii) for each vertex in the tree, a subset of vertices in the irrucible components. All vertices in that subset are glued together to form the graph. The differential can be decomposed into two parts: One part that changes one of the 1vi components, and one that does not. One can set up a spectral sequence, such that its first term is the latter part of the differential (leaving invariant the 1vi components).\footnote{Note that the spectral sequence converges because the graph complex splits into finite dimensional subcomplexes.} The resulting complex splits into subcomplexes according to the family of 1vi components. Fix one such subcomplex, say $C$, and fix one vertex $v$ in one of the $1vi$ components, that belongs to the subset associated to vertex $t$ of the tree. 
There is a filtration $C\supset C_1$, where $C_1$ is the subspace containing graphs such that the number of vertices in the subset of $t$ is one and $t$ has only one incident edge. Take the spectral sequence. Its first term contains a differential mapping 
\[
 C/C_1 \to C_1
\]
which can easily be seen to be an isomorphism.
\end{proof}

\section{A note on the convergence of spectral sequences}
\label{sec:specs}
In this paper, we use spectral sequences to compute the cohomology of graph complexes at various places, in particular for $\GC_n$ and $\fullGC_n$. We claim that these spectral sequences indeed converge to the cohomology.

For $\GC_n$ this is very easy to see. $\GC_n$ splits into a direct sum of finite dimensional subcomplexes, for fixed values of the difference
\[
 \Delta=e-v:=\#\text{edges}-\#\text{vertices}.
\]
Indeed, since each vertex is at least trivalent, one has $e\geq \frac{3}{2} v$ and hence $v\leq 2\Delta$ is bounded within each subcomplex. But there are only finitely many graphs with fixed $\Delta$ and bounded $v$, and hence each subcomplex is finite dimensional.

For $\fullGC_n$ the argument is a bit more subtle. Let $\mathcal{F}$ be a filtration on $\fullGC_n$, compatible with the differential. Let $\fullGC_n'$ be the same complex as $\fullGC_n$, but with the degrees shifted by $(n-\frac{1}{2})\Delta$, so that the new degree of a graph is
\[
 n(v-1)-(n-1)e +(n-\frac{1}{2})\Delta=\frac{1}{2}(v+e)-n.
\]
It is clear that (i) each grading component of $\fullGC_n'$ is finite dimensional, that (ii) the cohomology of $\fullGC_n'$ is the same as that of $\fullGC_n$ up to degree shifts and that (iii) one has a filtration $\mathcal{F}'$ on $\fullGC_n'$. By (i) the filtration $\mathcal{F}'$ is bounded and hence the associated spectral sequence converges to the cohomology of $\fullGC_n'$. But the spectral sequence associated to $\mathcal{F}$ is the same as that associated to $\mathcal{F}'$, up to some degree shifts. Hence it  converges to the degree shifted cohomology of $\fullGC_n'$, which by (ii) is the cohomology of $\fullGC_n$.

\section{\texorpdfstring{$\alg{t}$, $\grt_1$, $\sder$}{t, grt, sder}}
\label{sec:sder}
Let $\sder$ be the operad of Lie algebras of special derivations of free Lie algebras (see \cite{AT}, or \cite{drinfeld}). Elements of $\sder$ can be seen as internal trivalent trees in $\CG$, modulo the Jacobi identity.
As noted in \cite{pavol} there is a spectral sequence coming from the filtration on $\CG$ by internal loops, whose first term contains $\sder$. 

In particular, $\alg{t}$ is a sub-operad of Lie algebras of $\sder$. The map $\alg{t}\to \sder$ sends a generator $t_{ij}$ to the graph with a single edge between external vertices $i$ and $j$. 

\begin{prop}
 Let $\phi \in \alg{t}_3$ and $\Gamma$ its image in $\sder(3)$. 
%Assume $\phi$ has no $t_{13}$ component. 
Then $\phi$ can be recovered from $\Gamma$ as follows:
\begin{enumerate}
 \item Forget all graphs in $\Gamma$ that have more than one edge incident to vertex 2.
 \item The coeffiecient $c$ of the (single possible) graph with no vertex at 2 yields the coefficient of $t_{13}$ in $\phi$.
 \item Interpret each remaining tree as a Lie tree rooted at 2, corresponding to a Lie expression $\psi\in F_\Lie(X,Y)$.
 \item Then $\phi = c\cdot  t_{13} + \psi(t_{12},t_{13})$.
\end{enumerate}

\end{prop}
\begin{proof}
 By induction on the degree.
\end{proof}

When interpreting the tree as a Lie tree, we use the following sign convention.
In the ordering of the edges, the root edge of the tree must come first, then all vertices of its left subtree, and then all vertices of the right subtree. For each subtree, we apply the same convention recursively. For example, consider the following tree with the indicated ordering of the edges.
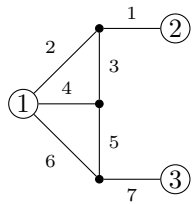
\[
\begin{tikzpicture}
\node[ext] (1) at (0,0) {$1$};
\node[ext] (2) at (2,1) {$2$};
\node[ext] (3) at (2,-1) {$3$};
\node[int] (w1) at  (1,1) {};
\node[int] (w2) at  (1,0) {};
\node[int] (w3) at  (1,-1) {};
\draw (1) edge node[auto] {$\scriptstyle 2$} (w1) edge node[auto] {$\scriptstyle 4$} (w2) edge node[below left] {$\scriptstyle 6$} (w3) 
      (w1) edge node[auto] {$\scriptstyle 1$} (2) edge node[auto] {$\scriptstyle 3$} (w2) 
      (w3) edge node[below] {$\scriptstyle 7$} (3) edge node[right] {$\scriptstyle 5$} (w2);
\end{tikzpicture}
\]
This tree is to be interpreted as the Lie word
\[
 [X,[X,[X,Y]]]\,.
\]

%TODO: add grt-part
%$\grt$ can be realized in \cite{AT} as a subalgebra, say $\grt'\cong \grt$, of $\sder_2$. The map $\grt'\to \grt$

\section{Twisting of operads}
\label{sec:optwists}
\subsection{A Foreword}
Any operad describes a certain kind of algebraic object. Often the algebraic object (the representation of the operad) is easier to describe and comprehend than the operad itself. So let us describe first what algebraic situation we want to consider by defining a ``twisted operad''.

Suppose we are given some operad $\op P$ together with a $\op P$-algebra $A$. Suppose further that we have a map $\hoLie\to \op P$. In particular, this means that $A$ is also a $\hoLie$-algebra. Sweeping convergence issues under the rug, it makes sense to talk about Maurer-Cartan elements in $A$, which are simply Maurer-Cartan (MC) elements in the $\hoLie$-algebra $A$. Fix such an $MC$ element $m$. Twisting the $\hoLie$ structure on $A$ by $m$, one can in particular endow $A$ with a new differential. Furthermore one can construct new operations on $A$ by inserting $m$'s into the $\op P$-operations on $A$. The twisted operad $\Tw\op P$ is defined such that the algebra $A$ with twisted differential and with this extended set of operations is a $\Tw \op P$-algebra.

The important claim for this paper is that there is an action of the deformation complex $\Def(\hoLie\to \op P)$ on the operad $\Tw\op P$. Concretely, it is defined as follows:
Suppose that we have a closed degree zero element $x\in \Def(\hoLie\to \op P)$. Such an element gives a $\hoLie$ derivation (i.e., infinitesimal automorphism) of the $\hoLie$-algebra $A$. Having an MC element $m\in A$, one can twist this derivation by $m$. One obtains in particular an (infinitesimally) different MC element $m'\in A$ %=m+\epsilon x\cdot m$
and hence also a new $\Tw \op P$-structure on $A$. The action of $x$ on the operad $\Tw\op P$ is defined such that it induces that change of $\Tw \op P$-structure.

\begin{rem}
 It is important that we define $\Tw \op P$ such that $A$ as above \emph{with differential twisted by $m$} is a $\Tw \op P$-algebra. 
\end{rem}

%  and (ii) an infinitesimal $\Lie_\infty$-morphism from the $\Lie_\infty$-algebra $A$, twisted by $m$, to $A$, twisted by $m'$.
% In particular, the first Taylor component intertwines the differential on $A$ twisted by $m$ with the differential on $A$ twisted by $m'$.
% The action of $x$ on the operad $\Tw\op P$ is given such that (i) the MC element $m$ is changed to $m'$ and (ii) an (infinitesimal) automorphism is applied to $A$, namely the first Taylor component of the twisted $\Lie_\infty$-morphism from before.

\subsection{The construction}
Let $\op P$ be any (dg) operad. We will consider it here as a contravariant functor from the category of finite sets (with bijections as morphisms) to the category of dg vector spaces. For a finite set $S$, the space $\op P(S)$ can be seen as the space of $\# S$-ary operations, with inputs labelled by elements of $S$. We will write for short $\op P(n):= \op P(\{1,\dots, n\})$. For some operation $a\in \op P(n)$ and for some symbols $s_1,\dots,s_n$ we will write 
\[
 a(s_1,\dots,s_n) \in \op P(\{s_1,\dots,s_n\})
\]
for the image of $a$ under the map $\op P(f)$, where $f: \{s_1,\dots,s_n\} \to [n]$ is the bijection sending $s_j\mapsto j$. Similarly, for $a\in \op P(n)$, $b\in \op P(m)$, and symbols $s_1,\dots,s_{m+n-1}$ we will write
\[
 a(s_1,\dots, s_{j-1},b(s_{j},\dots,s_{j+m-1}), s_{j+m},\dots,s_{m+n-1})
=(a\circ_{j}b)(s_1,\dots,s_{m+n-1})
\]
for the operadic composition.

Let $\hoLie_{k+1}:=(\Lie\{k\})_\infty$ be the minimal resolution of the degree shifted Lie operad. Concretely, the Lie bracket here has degree $-k$. An operad map $\hoLie_{k+1}\to \op P$ is described by a Maurer-Cartan element $\mu$ in 
\[
 \alg{g} = \Def(\hoLie_{k+1}\to \op P)\cong \prod_j (\op P(j)\otimes (\R[-k-1])^{\otimes j})^{S_j}[k+1].
\]
On the right, the permutation group $S_j$ acts on $\op P(j)$ as usual and on the tensor product of $\R[-k-1]$'s is by permutation, with appropriate signs.
Up to degree shift, elements of $\alg{g}$ are sums of symmetric or antisymmetric elements of $\op P$.

\begin{rem}
Fixing a Maurer-Cartan element $\mu\in \alg{g}$ fixes the notion of ``Maurer-Cartan element'' in any (nilpotent) $\op P$-algebra $A$. Concretely, the latter is just a Maurer-Cartan element in the $\hoLie_{k+1}$-algebra $A$. 
\end{rem}

We next want to define an operad $\Tw\op P$ that governs $\op P$-algebras ``twisted by'' a Maurer-Cartan element in the sense of the Remark. The underlying functor is defined on a set $S$ as
\[
 \Tw\op P(S) = \prod_{j\geq 0} (\op P(S\sqcup \{\bar 1,\dots \bar j\})\otimes (\R[k+1])^{\otimes j})^{S_j}
\]
where the symmetric group $S_j$ acts by permutation on the symbols $\bar 1,\dots, \bar j$ and be appropriate permutation (with signs) on the $\R[k+1]$. An element of the $j$-th factor in the product should be seen as an operation in $\op P$ invariant under permutations of the last $j$ slots. For an element $a$ in the $j$-th factor of $\Tw \op P(n)$ and some symbol set $S=\{s_1,\dots , s_{n+j}\}$ the expression
\[
 a(s_1,\dots,s_n,s_{n+1},\dots s_{n+j}) \in \op P(S)[j(k+1)]
\]
accordingly makes sense.
The operadic compositions are defined for homogeneous (wrt. both the grading by $j$ and the degree) elements $a\in\Tw\op P(m), b\in \Tw\op P(n)$
\begin{multline*}
 (a\circ_l b)(s_1,\dots, s_{m+n-1}, \bar 1,\dots \overline{j_1+j_2})
= \\ = 
\sum_{I\sqcup J=[j_1+j_2]}
\sgn(I,J)^{k+1}(-1)^{|b|j_1 (k+1)}
a(s_1,\dots,s_{l-1},b(s_l,\dots,s_{l+n-1},\bar J), \dots,s_{m+n-1}, \bar I)
\\
\in (\op P(S\sqcup \{\bar 1,\dots \overline{j_1+j_2} \})\otimes (\R[k+1])^{\otimes (j_1+j_2)})^{S_{j_1+j_2}}
\end{multline*}
where $\bar I$ is shorthand for $\bar i_1, \bar i_2,\dots$ with $i_1<i_2<\dots$ being the members of the set $I$, and similarly for $\bar J$. 
The sign is the sign of the shuffle permutation bringing $i_1, \dots ,j_1, \dots$ in the correct order.

\begin{rem}
An element of $\Tw\op P$ should be thought of as an operation in $\op P$, with Maurer-Cartan elements inserted into some of its slots.  
\end{rem}

Next we want to define a differential on $\Tw \op P$. First, denote by $\widetilde{\Tw}\op P$ the operad $\Tw \op P$ as defined so far, with the differential solely the one coming from the differential on $\op P$. On $\widetilde{\Tw}\op P$ we have a right action of the (dg) Lie algebra $\alg{g}$. Concretely, for homogeneous $x\in \alg{g}$, $a\in \widetilde{\Tw}\op P(n)$ we have 
\[
 (a\cdot x)(1,\dots,n,\bar 1,\dots \overline{j_1+j_2})
=
\sum_{I\sqcup J=[j_1+j_2]}
\sgn(I,J)^{k+1}
a(1,\dots,n, \bar I, x( \bar J)).
\]
\begin{lemma}
\label{lem:rightaction}
This formula describes a right action by operadic derivations. 
\end{lemma}
\begin{proof}
 For homogeneous elements we compute
\begin{align*}
& ((a\circ_l b)\cdot x ) (s_1,\dots,s_{m+n-1}, \bar 1,\dots, \overline{j_1+j_2+j_3})
\\
&=
\sum_{I\sqcup J\sqcup K=[j_1+j_2+j_3]}
\sgn(I\cup J,K)^{k+1}
\sgn(I,J)^{k+1}(-1)^{|b|j_1 (k+1)}
\\&\quad\quad\quad\quad\quad\quad
a(s_1,\dots,s_{l-1},b(s_l,\dots,s_{l+n-1},\bar J), \dots,s_{m+n-1}, \bar I, x(\bar K))
\\& \quad\quad\quad\quad\quad\quad +
\sgn(I\cup J,K)^{k+1}
\sgn(I,J)^{k+1}(-1)^{(|b|+|x|)j_1 (k+1)}
\\& \quad\quad\quad\quad\quad\quad
a(s_1,\dots,s_{l-1},b(s_l,\dots,s_{l+n-1},\bar J, x(\bar K)), \dots,s_{m+n-1}, \bar I)
\\&=
((-1)^{|b||x|}(a\cdot x)\circ_l b + a\circ_l (b\cdot x)) (s_1,\dots,s_{m+n-1}, \bar 1,\dots, \overline{j_1+j_2+j_3})
\end{align*}
For the last equality we used that 
\[
\sgn(I\cup J,K)\sgn(I,J)=\sgn(I, J\cup K)\sgn(J,K)=\sgn(I\cup K, J)\sgn(I,K)(-1)^{|J||K|}.
\]
\end{proof}
Of course, multiplying by a sign, one can make this right action into a left action.

For any operad $\op Q$, the unary operations $\op Q(1)$ form an algebra, hence in particular a Lie algebra, which acts on $\op Q$ by operadic derivations. Concretely, for $q\in \op Q(1), a\in \op Q(n)$ the formula is 
\begin{equation}
\label{equ:cdotdef}
 q\cdot a = q\circ_1 a -(-1)^{|a||c|}\sum_{j=1}^n a\circ_j q.
\end{equation}
Suppose that in addition some Lie algebra $\alg{h}$ acts from the left on $\op Q$ by operadic derivations. Then also the Lie algebra 
$\alg{h}\ltimes \op Q(1)$ acts on $\op Q$ by operadic derivations. 
Applying this to our case, we see that the Lie algebra 
\[
 \hat{\alg{g}}=\alg{g}\ltimes \widetilde{\Tw}\op P(1)
\]
acts on $\widetilde{\Tw}\op P(1)$ by operadic derivations. Given some (homogeneous) element $x \in \alg{g}$, we construct an element $x_1\in\widetilde{\Tw}\op P(1)$ by the formula
\begin{equation}
\label{equ:x1def}
 x_1(1,\bar 1,\dots , \bar j) = x(1,\bar 1,\dots , \bar j).
\end{equation}
%where $S$ is the degree shift operator.
\begin{lemma}
 If $\mu \in \alg{g}$ is a Maurer-Cartan element, then $\hat{\mu}=\mu -\mu_1$ is a Maurer-Cartan element in $\hat{\alg{g}}=\alg{g}\ltimes \widetilde{\Tw}\op P(1)$. % (TODO: check ``-'')
\end{lemma}
\begin{proof}
 Let us compute this for homogeneous $\mu$, the general case is analogous.
\begin{align*}
 \frac{1}{2}\co{\mu_1}{\mu_1}
(s_1, \bar 1, \dots , \overline{2m-2})
&=
\sum_{I\sqcup J = [2m-2]}
\sgn(I,J)^{k+1}
\mu(\mu(s_1, \bar I), \bar J)
\\
&=
-
\sum_{I\sqcup J = [2m-2]}
\sgn(I,J)^{k+1}
(-1)^{m(k+1)}
\mu(\mu(\bar I), s_1, \bar J)
\\ 
&=
-
\sum_{I\sqcup J = [2m-2]}
\sgn(I,J)^{k+1}
(-1)^{k+1}
\mu(s_1, \bar J, \mu(\bar I))
\\ 
&=
-
\sum_{I\sqcup J = [2m-2]}
\sgn(I,J)^{k+1}
(-1)^{m(m-1)(k+1)}
\mu(s_1, \bar I, \mu(\bar J))
\\ 
&=
-\mu_1\cdot \mu (s_1, \bar 1,\dots, \overline{2m-2})
=
\mu\cdot \mu_1(s_1, \bar 1,\dots, \overline{2m-2}) 
\end{align*}
\end{proof}
Given a Maurer-Cartan element in $\hat{\alg{g}}$ we can twist the differential on $\widetilde{\Tw}\op P$.
\begin{defi}
 Let $\op P$ be any operad and $k\in \mathbb{Z}$ be an integer. Let $\mu\in \Def(\Lie_\infty^{(k)}\to \op P)$ be a Maurer-Cartan element. Then we define the $\mu$-\emph{twisted operad} $\Tw \op P$ as the operad constructed above with differential
\[
 d_{\op P} + \hat{\mu}\cdot.
\]
\end{defi}

By construction, one has an action on $\Tw \op P$ by the $\mu$-twisted version of the Lie algebra $\alg{g}$, i.e., the Lie algebra $\alg{g}$ with the term $\co{\mu}{\cdot}$ added to the differential.

\begin{rem}
Let $\op P, k, \mu$ be as above and suppose that a $\op P$-algebra $A$ is given. Let $\alg{n}$ be a commutative nilpotent or pro-nilpotent algebra, e.g., $\alg{n}=\epsilon\R[[\epsilon]]$. Let $m\in A\otimes \alg{n}$ be a Maurer-Cartan element, i.e., a Maurer-Catan element in the $\hoLie_k$-algebra $A\otimes \alg{n}$. Then $A\otimes \alg{n}$ is an algebra over the operad $\Tw \op P$ by the formula
\[
 p(x_1,\dots, x_n) = \frac{1}{j!} p(x_1,\dots,x_n,m,\dots,m)
\]
for $p\in \Tw \op P$ homogeneous wrt. the degree in $j$ and $a_1,\dots,a_n\in A$.\footnote{I admit that the notation here is suboptimal. The $p$ on the left means the element of $\Tw \op P$, while the $p$ on the right is the underlying operation in $\op P$, which is (anti-)symmetric in its last $j$ slots.}
\end{rem}

\subsection{More explicit description of the action on the Graphs operad}
\label{sec:optwists_graphs}
Let us specialize the above constructions to the case $\op P = \Gra_n^\whl$. Then $\alg{g}=\Def(\hoLie_{n-1}\to \op P)=\fullGC_n^\whl$ is the ``full'' graph complex, containing all possible graphs, possibly with multiple connected components, tadpoles or multiple edges.
Elements of $\alg{g}$ should be considered as (possibly infinite) linear combinations of graphs with numbered vertices, invariant under permutation of vertex labels. In pictures we draw a graph with black unlabelled vertices. This should be understood as the sum of all possible numberings of the vertices, divided by the order of the automorphism group. Note that the picture is still inaccurate since we do not specify the overall sign. 
An explicit Maurer-Cartan element $\mu\in \alg{g}$ is given by the graph with two vertices and one edge.

By twisting one obtains the operad $\Tw\op P=\fullGraphs_n^\whl$. The $N$-ary operations of $\Tw\op P$ are (possibly infinite) linear combinations of graphs, having two kinds of numbered vertices, ``internal'' and ``external''. It is required that there are exactly $N$ external vertices and that the linear combination is invariant under interchange of the labels on the internal vertices. In pictures, we draw the internal vertices black without labels, with the convention that one should sum up over all possible numberings, and divide by the order of the symmetry group.
From the Maurer-Cartan element $\mu$ one obtains the element $\mu_1\in \Tw\op P(1)$. It is given by the graph with one external and one internal vertex, and an edge between them. The differential on $\Tw\op P$ has three terms: (i) There is a term coming from the action of $\mu$ as in Lemma \ref{lem:rightaction}. Concretely, this amounts to splitting each internal vertex into two and reconnecting the incoming edges. (ii) There is a term $\sum_{j} (\cdot)\circ_j \mu_1$. This amounts to splitting off from each external vertex one internal vertex, and reconnecting the incoming edges. (iii) There is a term $\mu_1\circ_1 (\cdot)$. This term adds a new internal vertex and connects it to every other vertex (but one at a time). Note that if all internal vertices are at least bivalent, the terms (iii) precisely cancel those terms from (i) and (ii) that contain graphs with univalent internal vertices.

Next consider more generally the action of an arbitrary element $\gamma\in \alg{g}$ on some $\Gamma\in \Tw\op P$. It again contains three terms:
(i) There is a term coming from the action as in Lemma \ref{lem:rightaction}. It amounts to inserting $\gamma$ at the internal vertices of $\Gamma$ and reconnecting the incoming edges. (ii) Build the element $\gamma_1\in \Tw\op P(1)$ by marking the vertex 1 in $\gamma$ as external. Then there is a term in the action stemming from  $\sum_{j} \Gamma \circ_j \mu_1$. This amounts to inserting $\gamma_1$ at all external vertices. (iii) There is the term $\gamma_1\circ_1 \Gamma$. This amount to inserting $\Gamma$ at the external vertex of $\gamma_1$.

\begin{rem}
\label{rem:Hopfaction}
A Hopf operad is an operad in the symmetric monoidal category of counital coalgebras (see \cite[section 5.3.5]{lodayval}). 
The operad $\Graphs_n$ is a (cocommutative) Hopf operad. Concretely, $\Graphs_n(N)=S(\ICG_n(N)[1])$ may be identified with a symmetric (co)product space (more precisely the the Chevalley-Eilenberg complex) of the $L_\infty$ algebra of internally connected graphs $\ICG_n(N)$, cf \eqref{equ:GraphsICG}.
The operad $e_n$ is a (quasi-isomorphic) sub-Hopf operad of $\Graphs_n$.
Combinatorially, the action of $\GC_n$ on $\Graphs_n$ always merges zero or more internally connected components into one new internally connected component. Hence it is compatible with (i.e., a derivation with respect to) the coproduct. Furthermore it respects the counit. Thus $\GC_n$ acts on $e_n$ by Hopf operad derivations, up to homotopy.
\end{rem}
  \section{The Tamarkin map \texorpdfstring{$\grt_1\to H^0(\GC_2)$}{grt to H(GC)} and Algorithm 2 of section \ref{sec:explicitmap}. }
\label{sec:tamalgoagree}
Part of the arguments in this subsection came out of a discussion with Pavol \v Severa.

\subsection{The map}\label{sec:tammap}
Let us first describe the map $\grt_1\to H^0(\GC_2)$, which is the universal version of the action of $\grt_1$ on $\Tpoly$ discussed in section \ref{sec:tamacttpoly}. As in section \ref{ssec:tamactionsginfty}, we have a chain of quasi-isomorphisms of operads 
\begin{equation}
\label{equ:cnepcdchain}
 C\Ne\PCD \rightarrow BU\alg{t} \leftarrow e_2 \leftarrow \hoe_2.
\end{equation}
Fix a lift up to homotopy $F:\hoe_2 \to C\Ne\PCD$ for now. It exists because $\hoe_2$ is cofibrant.

The Lie algebra $\grt_1$ acts on $\PCD$. Let some $\phi\in \grt_1$ be given. 
Its action on $F$ determines an element $\phi'\in \Def(\hoe_2 \to C\Ne\PCD)$.
Since the map $F$ is a quasi-isomorphism, we may decompose $\phi'$ into a lift $l\in \Der(\hoe_2)$ up to a homotopy $h\in \Def(\hoe_2 \to C\Ne\PCD)$. The lift may furthermore be chosen so that it has trivial $\hoLie_2$ part by degree reasons.
So we have the following diagram. 
\[
\begin{tikzpicture}
 \matrix(m)[diagram]{
 \hoLie_2  & &  \hoe_2 & &\\
 \hoe_2 & C\Ne\PCD & C\Ne\PCD& \Gra & \End(\Tpoly)
\\};
\draw[->]
  (m-1-1) edge (m-1-3)
  (m-1-1) edge (m-2-1)
  (m-2-1) edge (m-2-2)
  (m-1-3) edge (m-2-3)
  (m-2-2) edge node[below] {$\mathit{id}+\epsilon \phi\cdot$} (m-2-3)
  (m-2-1) edge[dashed]  node[auto]{$\id+\epsilon l$} (m-1-3)
  (m-2-3) edge (m-2-4)
  (m-2-4) edge (m-2-5);
  %(m-2-2) edge node[auto]{$\scriptstyle{\int}$} (m-2-3); 
\draw (m-2-3) +(-1.2,.7) node {$\stackrel{h}{\Longrightarrow}$};
\end{tikzpicture}
\]
The (infinitesimal) automorphism $l$ of $\hoe_2$ describes the action of the $\grt_1$-element $\phi$ on $\hoe_2$. 
The composition of the homotopy $h$ with the maps $\hoLie_2 \to \hoe_2$ from the left and $C\Ne\PCD \to \Gra$ from the right yields a degree 0 cocycle 
\[
 \xi\in \Def(\hoLie_2 \to \Gra) \cong \fGC
\]
which represents a cohomology class in the (full) graph complex. 
Since $H^0(\fGC)\cong H^0(\GC)$ this class corresponds to a cohomology class of M. Kontsevich's graph complex $\GC$. A representative in $\GC$ may be obtained by dropping from $\xi$ all graphs that are not connected or have a less then trivalent vertex. Composing $\xi$ with the map  $\Gra\to \End(\Tpoly)$ we recover D. Tamarkin's action of $\phi$ on $\Tpoly$.

%, to which $x$ is mapped under the map $\grt \to \GC$. 

\begin{rem}
Concretely, the map $C\Ne\PCD \to \Gra$ occuring above is constructed as follows.
First there is a map $C\Ne\PCD\to BU(\alg t)$ to the bar construction of the completed universal enveloping algebra of $\alg t$, by forgetting the parenthesization. Next one projects to the abelianization,  $ BU(\alg t)\to BU(\alg t/\co{\alg t}{ \alg t})$. For an abelian Lie algebra there is a canonical map from the bar construction of the universal enveloping algebra to the Chevalley complex, $BU(\alg t/\co{\alg t}{ \alg t})\to C(\alg t/\co{\alg t}{ \alg t})$. Finally one realizes that $C(\alg t/\co{\alg t}{ \alg t}) \cong \Gra$.
\end{rem}

\begin{rem}
 Note that the graph cohomology class $[\xi]$ we associate to the element $l\in \Der(\hoe_2)$ by the above recipe is the same as the one obtained using the map of Theorem \ref{thm:GCdef}: From $l$ we obtain an element in $\Def(\hoe_2 \to e_2)$ by composition with $\hoe_2\to e_2$ and  degree 1 element $l'\in \Def(\hoe_2 \to \Gra)$ by further composition with $e_2\to \Gra$.
Furthermore from the homotopy $h$ we obtain a degree zero element $h'\in \Def(\hoe_2 \to \Gra)$, whose coboundary is $l'$. Using these data to compute the connecting homomorphism in the long exact sequence from Proposition \ref{prop:homseq}, we see that graph cocycle obtained is just the $\hoLie_2$ part of $h'$, i. e., $\xi$.
 (Note that $l'$ is not guaranteed to be in the subcomplex $\Xi_\conn$, but one may remedy by adding some exact elements with trivial $\hoLie_2$ part.)
\end{rem}

\subsection{The lift \texorpdfstring{$\hoe_2 \to C\Ne\PCD$}{hoe2 to CNPaCD}}
To obtain more information one has to consider the lift $F:\hoe_2 \to C\Ne\PCD$ of $\hoe_2 \stackrel{f}\to BU(\alg t) \stackrel{g} \twoheadleftarrow C\Ne\PCD$. To determine $F$, it suffices to specify the image of the generators, i.~e., of $e_2^\vee(N)[1]\cong e_2^*\{2\}(N)[1]$ for $N=2,3,\dots$ (respecting $S_N$ equivariance).
These images are obtained by recursively (in $N$) solving the equations
\begin{align*}
 dF(x) &= F(dx) \\
 g(F(x)) &= f(x) 
\end{align*}
for $x\in e_2^\vee(N)[1]$, where we denote the differentials on $\hoe_2$ and $C\Ne\PCD$ both by $d$, abusing notation. Solutions exists since the maps $f$ and $g$ are quasi-isomorphisms, but they are not unique.

Note furthermore that there are additional gradings on the objects involved: On $C\Ne\PCD$ and $BU(\alg t)$ there is the grading by $t$-degree, while on $\hoe_2$ there is the grading of section \ref{sec:enfiltration}. The maps $f$, $g$ and the differential $d$ respect these gradings and hence we may choose $F$ in such a way that it respects the additional grading as well.

We will be interested in the restrictions to $\hoCom\subset \hoe_2$ and $\hoLie_2\subset \hoe_2$.
Denote by $m_2, m_3,\dots$ generators of $\hoCom$ and by $\mu_2,\mu_3,\dots$ the generators of $\hoLie_2$.
The gradings are as follows. $m_j$ has cohomological degree $2-j$, while the additional degree (of section \ref{sec:enfiltration}) is $0$. $\mu_j$ has cohomological degree $3-2n$ and additional degree $n-1$.

\begin{ex}
For $N=2$ we we may pick
\begin{align*}
 F(m_2) &= \frac 1 2 (1_{12}+1_{21}) & F(\mu_2) &= \frac 1 2 \left(t_{12}\otimes \id_{12} + t_{12}\otimes \id_{21} \right)
\end{align*}
where $1_{12}$ and $1_{21}$ denote the two zero-chains ``over'' the two objects $(12)$ and and $(21)$ of $\PCD(2)$, and $t_{12}\otimes \id_{12}$ and $t_{12}\otimes \id_{21}$ are one-chains made out of the identity morphisms of $(12)$ and $(21)$ and $t_{12}\in \alg t_2\subset U(\alg t_2)$.

Note also that the $t$ degree of $F(m_2)$ is zero and that of $F(\mu_2)$ is $1$, so that $F$ respects also the additional degrees.
\end{ex}

\subsection{Sketch of proof of Lemma \ref{lem:tamid} and hence of Theorem \ref{thm:taminj}}
\label{sec:tamproof}
Defining $m_3$ appropriately (there is a choice) one finds that $F(m_3)$ has to satisfy the equation
\[
 dF(m_3) = F(m_2)\circ_1 F(m_2) - F(m_2)\circ_2 F(m_2)
=
\frac 1 4 ( 1_{(12)3} +1_{(21)3} + 1_{3(12)} + 1_{3(21)} 
- 1_{1(23)}  - 1_{1(32)} - 1_{(23)1}  - 1_{(32)1}).
\]
where we used the notation $1_{(12)3}$ for the zero chain above object $(12)3$ etc.
As remarked before, we may pick all $F(m_j)$ to have $t$-degree 0.

Let us show Lemma \ref{lem:tamid}. Pick some $\grt_1$ element $\phi$. To obtain the corresponding element of $\Def(\hoe_2\to BU(\alg t) )$ we have to let it act on the lift map $F:\hoe_2$ and $C\Ne\PCD$ and then compose with the projection $C\Ne\PCD\to BU(\alg t)$. Note that the action on $F(m_j)$ produces a $j-2$-chain, in which all of the $j-2$ involved morphisms except for one have $t$-degree 0.
Since $BU(\alg t)$ is the \emph{normalized} bar construction, the image in $BU(\alg t)$ hence vanishes unless $j=3$.

\begin{rem} 
Let $h$ be the unique morphism in $\PCD(3)$ from any of the four objects $(ij)k$, $k(ij)$, $(ji)k$ or $k(ji)$ to any of the objects $i(jk)$, $(jk)i$, $(kj)i$ or $i(kj)$, where $i,j,k$ is a permutation of $1,2,3$.
Then acting with $\phi$ on the 1-chain described by $h$ and projecting the result to $BU(\alg t)$ yields the degree $-1$ element $\phi(t_{ij}, t_{jk})\in BU(\alg t)$.
\end{rem}

From the remark and the hexagon identity it is not hard to see that acting with $\phi$ on $F(m_3)$ and projecting the result $BU(\alg t)$ we obtain the element $\phi(t_{12}, t_{23})$. Note that explicit knowledge of $F(m_3)$ is not required here, only knowledge of $dF(m_3)$ as above.
But $\phi(t_{12}, t_{23})$ is the $\grt_1$-element we started with, considering $\grt_1$ as a subset of $\alg t_3$. This shows Lemma \ref{lem:tamid} and hence D. Tamarkin's Theorem \ref{thm:taminj}.

\subsection{Comparison to Algorithm 2 of section \ref{sec:explicitmap}}
\label{sec:tamalgodoagree}
% Also note that there are two gradings on $C\Ne\PCD$, one by the (non-positive) cohomological degree and one coming from the grading by number of $t_{ij}$'s on $\alg{t}$. We will call the latter grading the $t$-grading. The differential on $C\Ne\PCD$, call it $d$, has degree 1 with respect to the cohomological grading and degree 0 with respect to the $t$-grading.
% 
% Under $\hoe_2 \to C\Ne\PCD$ the generator $m_2$ has to be sent to a 0-chain in $C\Ne\PCD(2)$ of $t$-degree zero. Since $\PCD(2)$ has exactly two objects (namely $(12)$ and $(21)$) the space of such chains is two dimensional, generated by, say $1_{12}$ and $1_{21}$. By symmetry and the normalization condition that the $F(m_2)$ must be mapped to $1\in BU\alg t$ under $C\Ne\PCD\to BU\alg t$ we see that
% \[
%  F(m_2) = \frac 1 2 (1_{12}+1_{21}).
% \]
% Up to a prefactor, which we may choose to be $+1$ by defining $m_3$ appropriately, $F(m_3)$ then has to satisfy the equation 
% \[
%  dF(m_3) = F(m_2)\circ_1 F(m_2) - F(m_2)\circ_2 F(m_2)
% =
% \frac 1 4 ( 1_{(12)3}-1_{1(23)} + 1_{(12)3}-1_{1(23)} + )
% \]

% The goal of this section is to show the following result:
% \begin{prop}
% \label{prop:tamalgo2agree}
%  The cocycle $\xi$ describes the same graph cohomology class as the image of 
%  $x$ under the map $\grt_1\to H^0(\GC_2)\subset H^0(\Def(L_\infty \to \Gra))$ from Theorem \ref{thm:main}.
%  %cocycle $\Gamma$ produced by the Algorithm 2 in Section \ref{sec:explicitmap}.
% \end{prop}
% \begin{proof}
% Let us recall how the map of Theorem \ref{thm:main} was defined.
% First, 
% 
% \end{proof}

We will show that the cochain $\gamma'$ produced in Algorithm 2 of section \ref{sec:explicitmap} agrees with $\xi$ as in section \ref{sec:tammap} up to exact terms. Restrict the map $F:\hoe_2 \to C\Ne\PCD$ we picked to $\hoLie_2\subset \hoe_2$ to obtain a map $\hoLie_2\to C\Ne\PCD$. It is defined by specifying the images of the generators $\mu_2,\mu_3,\dots$. They have to satisfy equations of the form
\begin{equation}
\label{equ:muninduct}
d F(\mu_n) =(\text{linear combination of $\mu_j\circ \mu_{k}$ for $j+k=n+1$.}). 
\end{equation}
The cohomological degree of $F(\mu_n)$ must be $3-2n$. This means that $\mu_n$ is some linear combination of chains of morphisms of $\PCD$ of length $2n-3$. 
Furthermore, as before, we may pick $F$ such that $F(\mu_n)$ has $t$-degree $n-1$.

\begin{rem}
 In particular, this means that a chain of morphisms of $\PCD$ occuring in $\mu_n$ must contain at least $2n-3-(n-1)=n-2$ morphisms of $t$-degree zero. In particular, mapping $\mu_n$ along $C\Ne\PCD\to \Gra$ yields zero, except for $n=2$.
\end{rem}

Next consider the action of $\phi\in \grt_1$. It produces some 
\[
 a \in \Def(\hoe_2 \to C\Ne\PCD)
\]
of degree 1. Precompose with $\hoLie_2\to \hoe_2$ and compose with $C\Ne\PCD\to BU\alg{t}$ to obtain some element
\[
 a' \in \Def(\hoLie_2 \to BU\alg{t}).
\]
By the last Remark one sees that $a'$ is determined solely by the image of $\mu_3$ under the action of $x$. By an explicit calculation one can see that $a'$ in fact agrees with the element $T_3\in \Def(\hoLie_2 \to C\alg{t})\subset \Def(\hoLie_2 \to BU\alg{t})$ from the second algorithm in section \ref{ssec:tamactionsginfty}. Consider next the homotopy $h$. Precompose it with $\hoLie_2\to \hoe_2$ and compose with $C\Ne\PCD\to BU\alg{t}$ so as to obtain some degree 0 element
\[
 h' \in \Def(\hoLie_2 \to BU\alg{t}).
\]
Because $h$ was the homotopy making the lower right cell in the commutative diagram above commute, one has
\[
 a'=dh'
\]
where $d$ is now the differential in $\Def(\hoLie_2 \to BU\alg{t})$.
Since $a'=T_3$ actually lives in $\Def(\hoLie_2 \to C\alg{t})$ and $C\alg{t}\to BU\alg{t}$ is a quasi-isomorphism, there is a closed $c\in \Def(\hoLie_2 \to BU\alg{t})$, such that $U:=h'+c\in \Def(\hoLie_2 \to C\alg{t})$. Since $H^0(\Def(\hoLie_2 \to C\alg{t}))=0$, $c$ is in fact exact.
Taking $U$ for the $U$ in Algorithm 2 of section \ref{sec:explicitmap}, we see that the output $\gamma'$ of that algorithm is the image of $U$ after composition with $BU\alg{t}\to \Gra$. The image of $h'$ after composition with $BU\alg{t}\to \Gra$ is $\xi$, hence 
\[
 \gamma' = \xi +(\text{exact terms}).
\]

\section{A short note on the directed version of the graph complex.}
One may define a directed version of the graph complex $\GC_n$ by (i) taking directed instead of undirected edges, (ii) allowing bivalent vertices and (iii) requiring that a graph has at least one trivalent vertex.\footnote{This last condition is to remove the loops as in Figure \ref{fig:loops} and otherwise unnecessary.} Call the resulting graph complexes $\dGC_n$.
\begin{prop}
 \[
H(\dGC_n) \cong H(\GC_n).
 \]
\end{prop}
There is even an explicit quasi-isomorphism of dg Lie algebras
\[
 \GC_n \to \dGC_n
\]
sending an undirected graph to a sum of directed graphs, obtained by interpreting each edge as the sum of edges in both directions.
\begin{proof}[Proof of the Proposition]
Set up a spectral sequence on the number of bivalent vertices. The first differential produces one bivalent vertex. As in the undirected case, consider a graph with bivalent vertices as one with trivalent vertices and labelled edges, see Figure \ref{fig:frombivtotrivlab}. The first differential just changes the labels. The resulting subcomplex is (essentially) a tensor product the complex for one edge. The cohomology of the latter can be seen to have a single nontrivial cohomology class, represented by the sum of edges in both direcions (i.e., $\leftarrow + \rightarrow$). Hence the first convergent in the spectral sequence is $\GC_n$ and hence $\GC_n\to \dGC_n$ is a quasi-isomorphism.
\end{proof}

\begin{figure}
 \centering
 \[
 \begin{tikzpicture}[scale=1,
vert/.style={draw,outer sep=0,inner sep=0,minimum size=5,shape=circle,fill},
helper/.style={outer sep=0,inner sep=0,minimum size=5,shape=coordinate},
default_edge/.style={draw},
diredge/.style={-triangle 60},
ext/.style={draw,outer sep=0,inner sep=0,minimum size=5,shape=circle},
helper2/.style={outer sep=0,inner sep=0,minimum size=5,shape=coordinate},
every loop/.style={}]

\node (v0) at (2.5,8.5) [vert] {};
\node (v1) at (1,8.5) [vert] {};
\node (v3) at (1.5,7) [vert] {};
\node (v5) at (3,7.5) [vert] {};
\node (v7) at (4.5,7.5) [vert] {};
\node (v9) at (6,7) [vert] {};
\node (v11) at (4.5,6) [vert] {};
\node (v13) at (2.5,6) [vert] {};
\node (v17) at (6,5.5) [vert] {};
\node (v22) at (9.5,7.5) [vert] {};
\node (v23) at (12,7.5) [vert] {};
\node (v25) at (11,6) [vert] {};
\node (v31) at (9.5,7.9) [helper2,label=90:{"$<>>$"}] {};
\node (v32) at (10.5,7.8) [helper2,label=0:{"$>>>$"}] {};
\node (v33) at (9.3,6.5) [helper2,label=0:{"$><$"}] {};
\node (v34) at (11.8,6.5) [helper2,label=0:{"$><$"}] {};
\node (v35) at (10.8,6.9) [helper2,label=0:{"$>$"}] {};

\draw[diredge] (v1) to (v0);
\draw[diredge] (v1) to (v3);
\draw[diredge] (v3) to (v5);
\draw[diredge] (v5) to (v7);
\draw[diredge] (v7) to (v9);
\draw[diredge] (v9) to (v11);
\draw[diredge] (v11) to (v13);
\draw[diredge] (v0) to (v3);
\draw[diredge] (v3) to (v13);
\draw[diredge] (v9) to (v17);
\draw[diredge] (v11) to (v17);
\draw[default_edge] (v22) to (v23);
\draw[default_edge,in=-10, out=-90] (v23) to (v25);
\draw[default_edge] (v22) to (v25);
\draw[default_edge,in=135, out=45,loop] (v22) to ();
\draw[default_edge] (v25) to (v23);
\end{tikzpicture}
\]
\caption{\label{fig:frombivtotrivlab} Passing from a graph with bivalent vertices to one with trivalent vertices only, but labelled edges. Note that we cheat a little, since one has to assign some directions, by some covention, to the edges on the right so as to interpret the labels correctly. However, for the argument in the proof this does not matter.}
\end{figure}

\section{The automorphism group of \texorpdfstring{$\hoe_n$}{hoen}}
\label{app:unipotent}

Let $\op C$ be a coaugmented cooperad with zero differential, finite dimensional in each arity. We will consider the cobar construction $\Omega(\op C)$ and its automorphism group $\Aut(\Omega(\op C))$. First note that by functoriality of the cobar construction we have a map from the automorphism group of the coaugmented cooperad $\op C$
\[
\phi: \Aut(\op C) \to \Aut(\Omega(\op C)).
\]
Recall from \cite[section 6.5]{lodayval} that elements of $\Omega(\op C)$ may be understood as linear combinations of certain trees, whose vertices are decorated by elements of $\overline{\op C}$. 
The operadic composition is grafting of trees. Hence there is a filtration by the number of vertices in a tree on the operad $\Omega(\op C)$. Concretely, $\mF^p\Omega(\op C)$ is spanned by trees with $p$ or more vertices. By quasi-freeness any automorphism of $\Omega(\op C)$ preserves this filtration.
Clearly 
\[
\overline{\op C} \cong (\Omega(\op C)/\mF^2 \Omega(\op C))[1]
\]
and hence we have a map 
\[
\psi: \Aut(\Omega(\op C))\to \Aut_{\S-\rm{mod}}(\op C)
\]
where $\Aut_{\S-\rm{mod}}(\op C)$ is the group of automorphisms of the $\S$-module $\op C$.
It is clear that $\psi\circ\phi$ agrees with the inclusion $\Aut(\op C)\subset \Aut_{\S-\rm{mod}}(\op C)$.
We claim that the image of $\psi$ actually lands in $\Aut(\op C)$. Indeed any automorphism $f\in \Aut(\Omega(\op C))$ has to respect the differential. On generators, The differential is (see \cite[6.5.5]{lodayval}) just the infinitesimal decomposition (see \cite[6.1.7]{lodayval}) on $\op C$, up to degree shifts.\footnote{Here we use that there is no differential on $\op C$.} It follows that $\psi(f)$ has to respect this infinitesimal decomposition. But the latter generates all cocompositions in $\op C$ and hence $\psi(f)$ has to respect the cooperad structure. (The coaugmentation is also respected by construction.)
Thus we have have a decomposition
\[
\Aut(\Omega(\op C)) = \Aut(\op C) \ltimes  \Aut_1(\Omega(\op C))
\]
where we abbreviate 
\[
\Aut_1(\Omega(\op C)) = \ker \phi.
\]

\begin{lemma}
$\Aut_1(\Omega(\op C))$ is a pro-unipotent group.
\end{lemma}
\begin{proof}
Let $\op P_n\subset \Omega(\op C)$ be the suboperad spanned by trees all of whose vertices have $\leq n$ children. Equivalently this sub-operad is generated by an $\op S$ module obtained from $\op C$ by setting to zero all $\op C(N)$ for $N>n$.
By quasi-freeness all automorphisms of $\Omega(\op C)$ have to respect $\op P_n$. Similarly all automorphisms of $\op P_n$ have to respect $\op P_{n-1}\subset \op P_n$ etc.
We hence have maps
\[
\Aut(\op P_{n-1}) \leftarrow \Aut(\op P_n).
\]
Clearly $\lim_\leftarrow \Aut(\op P_n) \cong \Aut(\Omega(\op C))$.
The filtration $\mF$ from above induces filtrations on each $\op P_n$.
Denote by $\Aut_1(\op P_n)\subset \Aut(\op P_n)$ the subgroup of automorphisms fixing all generators modulo elements of $\mF^2$. Then, in the same manner as above we have arrows 
\[
\Aut_1(\op P_{n-1}) \leftarrow \Aut_1(\op P_n)
\]
and $\lim_\leftarrow \Aut_1(\op P_n) \cong \Aut_1(\Omega(\op C))$.

We claim that each $\Aut_1(\op P_n)$ is a unipotent algebraic group. Indeed by the finiteness assumption on $\op C$ it is an algebraic subgroup of some $GL(\op P_n(n))$. Furthermore, if $g\in \Aut_1(\op P_n)$ then $(g-\id)$ maps $\mF^p\op P_n$ into $\mF^{p-1}\op P_n$ and hence is a nilpotent element.
\end{proof}

One may also verify that the Lie algebra of $\Aut_1(\Omega(\op C))$ is given by the closed degree zero elements in $\Der'(\Omega(\op C))$ (as defined in section \ref{sec:defcomplexes}).

\begin{rem}[The special case $\hoe_n$]
Let us specialize to the case $\op C=e_n^\vee$ relevant to this paper.
In this case $\Aut(\op C)=\Aut(e_n^\vee)\cong \K^\times \times \K^\times$, the two factors acting by rescaling the two cogenerators. No (non-trivial) such rescaling is homotopy trivial.
\end{rem}

% to be compatible with the differential it has to respect the infinitesimal coproduct 
%
%compatibility with the differential is equivalent to  that the infinitesimal cocompositions 
%
%
%
%There is an operadic ideal $I\subset \hoe_n$ spanned by trees with at least two levels. 
%
%
%
%with zero differential, finite dimensional in each arity
%
%
%The goal of this section is to show that the exponential map 
%\[
%\exp : \Der(\hoe_n)' \to \Aut( \hoe_n ) 
%\]
%exists and is onto. 
%First recall from \cite[section ]{lodayval} that elements of $\hoe_n=\Omega(\hoe_n^\vee)$ may be understood as linear combinations of certain trees, decorated by elements of $\overline{\hoe}_n^\vee$. 
%The operadic composition is grafting of trees. There is an operadic ideal $I\subset \hoe_n$ spanned by trees with at least two levels. 

\nocite{Turchin1, Turchin2, Turchin3, LambrechtsTurchin}
\bibliographystyle{plain}
\bibliography{biblio} 

%\begin{thebibliography}{9}
%\bibitem{AT1} Alekseev Torossian KV
%\bibitem{AT2} Alekseev Torossian graphs
%\bibitem{K1}Kontsevich, quantization
%\bibitem{K2} Kontsevich, operads
%\bibitem{LV} Lambrechts Volic
%\end{thebibliography}

\end{document}